\documentclass[12pt]{article}
\def\date{18 January 2012}  %This will appear on page 1 as "revision date"
\usepackage{amsfonts}
\usepackage{amssymb}
\usepackage{amsmath}
\usepackage{graphicx}
\usepackage{subfigure}
\usepackage{theorem}

\newtheorem{proposition}{Proposition}[section]
\newtheorem{lemma}[proposition]{Lemma}
\newtheorem{corollary}[proposition]{Corollary}
\newtheorem{theorem}[proposition]{Theorem}
\newtheorem{claim}[proposition]{Claim}

\newcommand{\qed}{$\square$\bigskip}

\newcommand{\proof}{{\noindent\bf Proof. }}
{\theorembodyfont{\rmfamily}

}
\def\mytextindent#1{\indent\indent\llap{\rm#1\enspace}\ignorespaces}
\def\fattextindent#1{\indent\indent\llap{#1\enspace}\ignorespaces}
\def\myitem{\par\hangindent0pt\mytextindent}
\def\myitemitem{\par\hangindent\parindent\fattextindent}
\def\claim#1#2{\bigskip\noindent\rlap{\rm(#1)}\ignorespaces
   \hangindent=30pt\hskip30pt{\sl#2}\bigskip}

\def\myclaim#1#2\par{{\medbreak\noindent\rlap{\rm(#1)}\ignorespaces
 \rightskip20pt
 \hangindent=20pt\hskip20pt{\ignorespaces\sl#2}\smallskip}}
\def\junk#1{}

\def\mylabel#1{{\label{#1}}}
\font\smallrm=cmr8
\font\smallrm=cmr8
\newcommand{\tengirl}{L_1}
\newcommand{\tenguy}{L_2}

\newcommand{\numbernine}{L_4}
\newcommand{\elevengirl}{L_5}
\newcommand{\elevenguy}{L_6}

\begin{document}

\phantom{a}\vskip .25in
\centerline{\bf FIVE-COLORING GRAPHS ON THE KLEIN BOTTLE}

\vskip.5in
%authors
\centerline{Nathan Chenette%
  \footnote{\smallrm Partially supported by NSF under
             Grant No.~DMS-0200595.}
}
\centerline{Luke Postle$^{1,3}$
}
\centerline{Noah Streib$^1$
}
\centerline{Robin Thomas%
  \footnote{\smallrm Partially supported by NSF under
             Grants No.~DMS-0200595,  DMS-0354742, and DMS-0701077. }
}
\centerline{Carl Yerger%
  \footnote{\smallrm Partially supported by NSF Graduate Fellowship.}
}
\bigskip
\centerline{School of Mathematics}
\centerline{Georgia Institute of Technology}
\centerline{Atlanta, Georgia  30332-0160, USA}

\vskip 0.5in 
\centerline{\bf ABSTRACT}
\bigskip
\parshape=1.5truein 5.5truein
We exhibit an explicit list of nine graphs such that a graph
drawn in the Klein bottle is $5$-colorable if and only if it
has no subgraph isomorphic to a member of the list.

\vfill \noindent 31 August 2007, revised \date.
%\hfil\break\noindent This is an incomplete draft not intended for distribution.
\vfil\eject
%End of title page
%%%%%%%%%%%%%%%%%%%%%%%%%%%%%%%%%%%%%%%%%%%%%%%%%%%%%%%%%%%%%%%%%%%%%%%%%%%%

\section{Introduction}

All graphs in this paper are finite, undirected and simple.
By a {\em surface} we mean a compact, connected $2$-dimensional 
manifold with empty boundary.
The classification theorem of surfaces 
(see e.g.~\cite{Massey}) states that each surface
is homeomorphic to either $S_g$, the surface obtained from the sphere by
adding $g$ handles, or $N_k$, the surface obtained from the sphere
by adding $k$ cross-caps.
Thus $S_0=N_0$ is the sphere, $S_1$ is the torus, $N_1$ is the
projective plane and $N_2$ is the Klein bottle.

In this paper we study a specific instance of the following more
general question: Given a surface $\Sigma$ and an integer $t\ge0$,
which graphs drawn in $\Sigma$ are $t$-colorable?

Heawood~\cite{Heawood} proved that if $\Sigma$ is not the
sphere, then every graph in $\Sigma$
is $t$-colorable as long as 
$t\ge H(\Sigma):=\lfloor(7+\sqrt{24\gamma+1})/2\rfloor$,
where $\gamma$ is the {\em Euler genus of $\Sigma$}, 
defined as $\gamma=2g$ when $\Sigma=S_g$ and $\gamma=k$ when $\Sigma=N_k$.
Incidentally, the assertion holds for the sphere as well, by the
Four-Color Theorem~\cite{AppHak1,AppHakKoc,AppHak89,RobSanSeyTho4CT}.
Ringel and Youngs (see~\cite{Ringel}) proved that the bound is
best possible for all surfaces except the Klein bottle.
Dirac~\cite{Dirmap} and Albertson and Hutchinson~\cite{AlbHut}
improved Heawood's result by showing that
every graph in $\Sigma$ is actually $(H(\Sigma)-1)$-colorable,
unless it has a subgraph isomorphic to the complete graph
on $H(\Sigma)$ vertices.

We say that a graph is {\em $(t+1)$-critical} if it is not $t$-colorable,
but every proper subgraph is.
Dirac~\cite{Dircritical} also proved that for every $t\ge8$ and every
surface $\Sigma$ there
are only finitely many $t$-critical graphs on $\Sigma$.
Using a result of Gallai~\cite{Galcritical} it is easy to
extend this to $t=7$. In fact, the result extends to $t=6$ by the
following deep theorem of Thomassen~\cite{ThoCritical}.

\begin{theorem}
For every surface $\Sigma$ there are only finitely many $6$-critical
graphs in $\Sigma$.
\end{theorem}

Thus for every $t\ge5$ and every surface $\Sigma$ there exists a
polynomial-time algorithm to test whether a graph in $\Sigma$
is $t$-colorable. What about $t=3$ and $t=4$? For $t=3$ the
$t$-coloring decision problem is NP-hard even when $\Sigma$ is
the sphere~\cite{GarJoh}, and therefore we do not expect to be
able to say much.
By the Four-Color Theorem the $4$-coloring decision problem is trivial 
when $\Sigma$ is the sphere, but it is open for all
other surfaces. 
A result of Fisk~\cite{Fis} can be used to construct infinitely
many $5$-critical graphs on any any surface other than the sphere,
but the structure of $5$-critical graphs on surfaces 
appears complicated~\cite[Section~8.4]{MohTho}.
\nocite{MohTho}

\begin{figure}
\centering
\includegraphics[scale=.8]{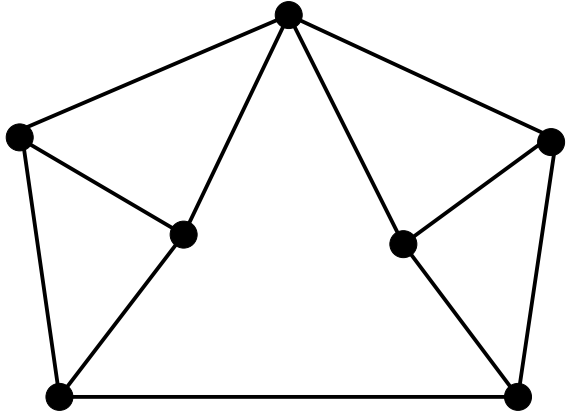}
\caption{The graph $H_7$}
\label{fig:h7}
\end{figure}

Thus the most interesting value of $t$ for the $t$-colorability
problem on a fixed surface seems to be $t=5$.
By the Four-Color Theorem
every graph in the sphere is $4$-colorable, but on every other
surface there are graphs that cannot be $5$-colored.
Albertson and Hutchinson~\cite{AlbHut} proved that a graph in
the projective plane is $5$-colorable if and only if it has no
subgraph isomorphic to $K_6$, the complete graph on six vertices.
Thomassen~\cite{Tho5torus}
 proved the analogous (and much harder) result for the
torus, as follows.
If $K,L$ are graphs, then by $K+L$ we denote the graph obtained
from the union of a copy of $K$ with a disjoint copy of $L$ by
adding all edges between $K$ and $L$.
The graph $H_7$ is depicted in Figure~\ref{fig:h7} and the graph $T_{11}$ is
obtained from a cycle of length $11$ by adding edges joining
all pairs of vertices at distance two or three.

\begin{theorem}
\mylabel{thomassentorus}
A graph in the torus is $5$-colorable if and only if it has no
subgraph isomorphic to $K_6$, $C_3+C_5$, $K_2+H_7$, or $T_{11}$.
\end{theorem}

\begin{figure}
\centering
\includegraphics[scale=1.2]{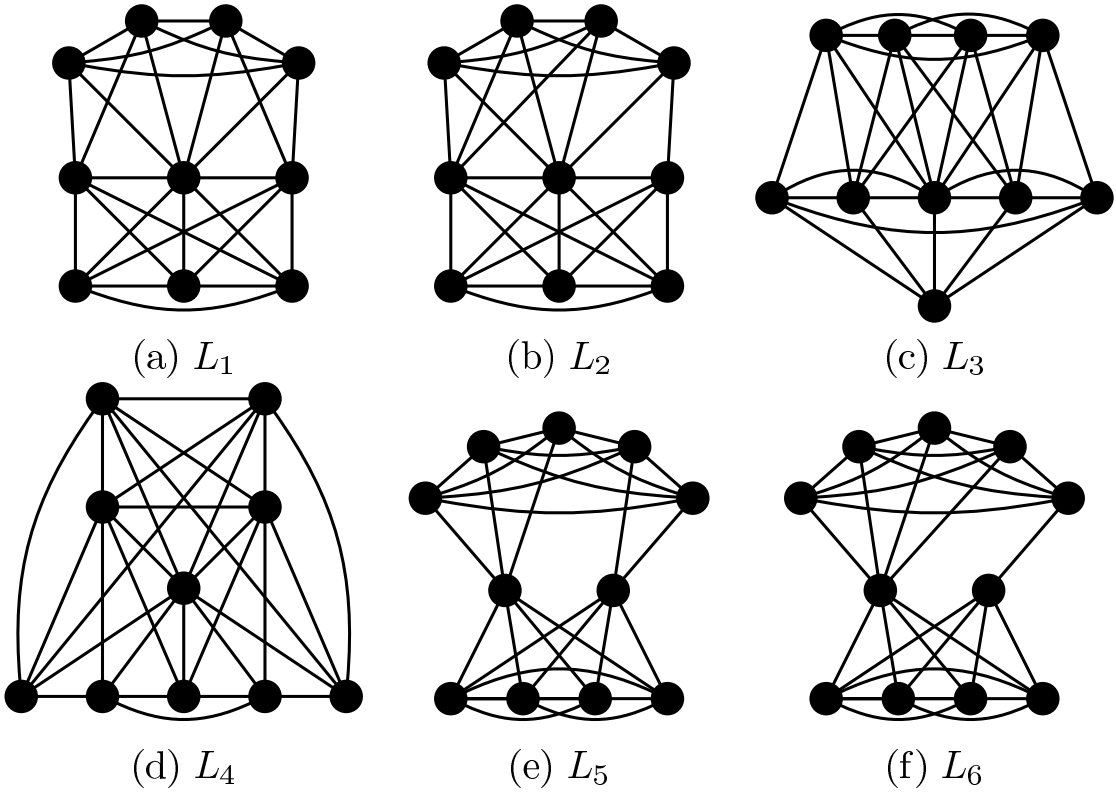}
\caption{The graphs $L_1,L_2,\ldots,L_6$}
\label{fig:allels}
\end{figure}

Our objective is to prove the analogous result for the Klein bottle,
stated in the following theorem. 
The graphs $L_1,L_2,\ldots,L_6$ are defined in Figure~\ref{fig:allels}.
Lemma~\ref{K5MINUS} explains how most of these graphs arise in the proof.

\begin{theorem}
\mylabel{main}
A graph in the Klein bottle is $5$-colorable if and only if it has no
subgraph isomorphic to $K_6$, $C_3+C_5$, $K_2+H_7$, or any of the
graphs $L_1,L_2,\ldots,L_6$.
\end{theorem}

Theorem~\ref{main} settles a problem of 
Thomassen~\cite[Problem~3]{ThoCritical}.
It also implies that in order to test $5$-colorability of a graph $G$ drawn in the
Klein bottle it suffices to test subgraph isomorphism to one of the
graphs listed in Theorem~\ref{main}.
Using the algorithms of~\cite{Eppsubiso} and~\cite{Mohlinear}
we obtain the following corollary.

\begin{corollary}
There exists an explicit linear-time algorithm to decide whether
an input graph embeddable in the Klein bottle is $5$-colorable.
\end{corollary}

It is not hard to see that with the sole exception of $K_6$, none of
the graphs listed in Theorem~\ref{main} can be a subgraph of an Eulerian
triangulation of the Klein bottle. Thus we deduce the following
theorem of
Kr\'al', Mohar, Nakamoto, Pangr\'ac and Suzuki~\cite{KraMohNakPanSuz}.

\begin{corollary}
\label{kbtriang}
An Eulerian triangulation of the Klein bottle is $5$-colorable if and
only if it has no subgraph isomorphic to $K_6$.
\end{corollary}

It follows by inspection that each of the graphs from Theorem~\ref{main}
has a subgraph isomorphic to a subdivision of $K_6$. Thus we deduce the
following corollary.

\begin{corollary}
\label{K6subd}
If a graph in the Klein bottle is not $5$-colorable, then it has a subgraph
isomorphic to a subdivision of $K_6$.
\end{corollary}

This is related to Haj\'os' conjecture, which states that for every
integer $k\ge1$, if a graph $G$ is not $k$-colorable, then it has 
a subgraph isomorphic to a subdivision $K_{k+1}$.
Haj\'os' conjecture is known to be true for $k=1,2,3$ and false
for all $k\ge6$. The cases $k=4$ and $k=5$ remain open.
In~\cite[Conjecture~6.3]{ThoHajos} Thomassen conjectured that Haj\'os'
\nocite{ThoHajos}
conjecture holds for  every graph in the projective plane or the torus.
His results~\cite{Tho5torus} imply that it suffices to prove this
conjecture for $k=4$, but that is still open.
Likewise, one might be tempted to extend Thomassen's conjecture to
graphs in the Klein bottle; Corollary~\ref{K6subd} then implies that
it would suffice to prove this extended conjecture for $k=4$.

Thomassen proposed yet another related 
conjecture~\cite[Conjecture~6.2]{ThoHajos} stating that every graph
which triangulates some surface satisfies Haj\'os' conjecture.
He also pointed out that this holds for $k\le4$ for every surface
by a deep theorem of Mader~\cite{MadDirac}, and that it holds
for the projective plane and the torus by~\cite{Tho5torus}.
Thus Corollary~\ref{K6subd} implies that Thomassen's second conjecture
holds for graphs in the Klein bottle.
For general surfaces the conjecture was disproved by Mohar~\cite{MohHajos}.
Qualitatively stronger counterexamples were found by R\"odl and Zich~\cite{RodZich}.

Our proof of Theorem~\ref{main} follows closely the argument of~\cite{Tho5torus},
and therefore we assume familiarity with that paper.
We proceed as follows. 
The result of Sasanuma~\cite{Sas6reg}
%First we show, using the description of all $6$-regular graphs
%in the Klein bottle, 
that every $6$-regular graph in the Klein bottle is $5$-colorable
(which follows from the description of all $6$-regular graphs on the
Klein bottle)
allows us to select a minimal counterexample 
$G_0$ and a suitable vertex $v_0\in V(G_0)$ of degree five.
If every two neighbors of $v_0$ are adjacent, then $G_0$
has a $K_6$ subgraph and the result holds.
We may therefore select two non-adjacent neighbors $x$ and $y$
of $v_0$. Let $G_{xy}$ be the graph obtained from $G_0$ by deleting
$v_0$, identifying $x$ and $y$ and deleting all resulting parallel
edges. If $G_{xy}$ is $5$-colorable, then so is $G_0$, as is easily
seen. Thus we may assume that $G_{xy}$ has a subgraph isomorphic
to one of the nine graphs on our list, and it remains to show
that either $G_0$ can be $5$-colored, or it has a subgraph
isomorphic to one of the nine graphs on the list.
That  occupies most of the paper.

We would like to acknowledge that Theorem~\ref{main}
was independently obtained by Kawara\-bayashi, Kr\'al', Kyn\v{c}l, and 
Lidick\'y~\cite{KawKraKynLid}.
Their method relies on a computer search.
The result of this paper forms part of the doctoral dissertation~\cite{YerPhD} 
of the last author.

\section{Lemmas}
%%%%%%%%%%%%%%%%%%%%%%%%%%%%%%%%%%%%%%%%%%%%%%%%%%%%%%%%%%%%%%%%%%%%%%%%%%

Our first lemma is an adaptation of~\cite[Theorem~6.1, Claim~(8)]{Tho5torus}.
%A graph is {\em $(t+1)$-vertex-critical} if it is not $t$-colorable,
%but every proper induced subgraph is.

\begin{lemma}
\mylabel{smallcrit}
Let $G$ be a graph in the Klein bottle that is not $5$-colorable
and has no subgraph isomorphic to $K_6$, $C_3 + C_5$, or $K_2 + H_7$.
Then $G$ has at least $10$ vertices, and if it has exactly $10$, then
it has a vertex of degree nine.
\end{lemma}

\proof
We follow the argument of~\cite[Theorem~6.1, Claim~(8)]{Tho5torus}.
%also serves as a proof of this lemma.~\qed
Let $G$ be as stated, and let it have at most ten vertices.
We may assume, by replacing $G$ by a suitable subgraph, 
that $G$ is $6$-critical.
%We now employ a
%result of Gallai~\cite{Galcritical}, which states that if $G$ is a
%$k$-critical graph with at most $2k - 2$ vertices, then $G$ is of
%the form $G = G_1 + G_2$, where $G_i$ is $k_i$-critical for $i = 1,
%2$, and $k_1 + k_2 = k$.  
By a result of Gallai~\cite{Galcritical} it follows that $G$ is
of the form $H_1 + H_2$, where $H_i$ is $k_i$-critical, $k_1
\leq k_2$, and $k_1 + k_2 = 6$.  If $k_1 = k_2 = 3$, then we obtain
that $G$ is isomorphic to
either $K_6$ or $C_3 + C_5$, a contradiction.  So $k_1
\leq 2$ and therefore $G$ has a vertex adjacent to all other
vertices.  Now, suppose for purposes of contradiction that $|V(G)|\leq 9$.
If $k_1 = 1$, then $|V(H_2)| \leq 8$ and so $H_2$ is of
the form $H_2' + H_2''$, where $H_2' = K_2$ or $K_1$.
Thus we may assume that  $k_1 = 2$
and that $H_2$ is 4-critical.  By the results of
\cite{Galcritical} and \cite{Toft}, the only 4-critical
graphs with at most seven vertices are $K_4, K_1 + C_5, H_7$ and
$M_7$, where $M_7$ is obtained from a 6-cycle, $x_1x_2\cdots x_6x_1$
by adding an additional vertex $v$ and edges $x_1x_3, x_3x_5,
x_5x_1, vx_2, vx_4, vx_6$.  However,
$G$ has no subgraph isomorphic to $K_2+K_4=K_6$,
$K_2 + (K_1+ C_5) = C_3 + C_5$, or $K_2 + H_7$.  This
implies that $G$ is isomorphic to $K_2 + M_7$.
The latter graph  has nine
vertices and  27 edges, and so 
%$G$ is isomorphic to $K_2 + M_7$ and 
triangulates the Klein bottle. However, $K_2 + M_7$ has
a vertex whose neighborhood is not Hamiltonian,
a contradiction.~\qed

Our next lemma is an extension of~\cite[Lemma~4.1]{Tho5torus},
which proves the same result for cycles of length at most six.
If $C$ is a subgraph of a graph $G$ and $c$ is a coloring of $C$,
then we say that a vertex $v\in V(G)-V(C)$ {\em sees a color $\alpha$ on $C$}
if $v$ has a neighbor $u\in V(C)$ such that $c(u)=\alpha$.

\begin{figure}
\centering
\includegraphics[scale=1.1]{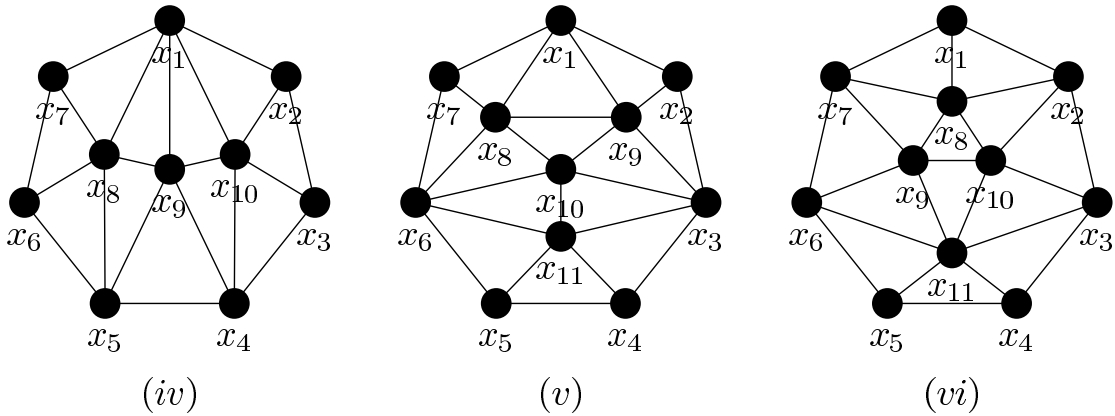}
\caption{Graphs that have non-extendable colorings}
\label{fig:criticalgraphs}
\end{figure}

\begin{lemma} 
\mylabel{7cycle-5coloring}
Let $G$ be a plane graph with an outer cycle $C$ of length $k\le7$, and
%$S: x_1 x_2 \cdots x_k x_1$, $k \le 7$. 
let $c$ be a $5$-coloring of $G[V(C)]$. Then $c$ cannot 
be extended to a 5-coloring of $G$ if and only if $k\ge5$ and
the vertices of $C$ can be numbered $x_1,x_2,\ldots,x_k$ in order
such that one of the following conditions hold:
\myitem{(i)}some vertex of $G-V(C)$ sees five distinct colors on $C$,
\myitem{(ii)}$G - V(C)$ has two adjacent vertices that both see the
same four colors on $C$,
\myitem{(iii)} $G - V(C)$ has three pairwise adjacent vertices that each
see the same three colors on $C$,
\myitem{(iv)} $G$ has a subgraph isomorphic to the first graph shown
in Figure~\ref{fig:criticalgraphs}, 
and the only pairs of vertices of $C$ colored the same  are 
either $\{x_5,x_2\}$ or $\{x_5,x_3\}$, and 
either $\{x_4, x_6\}$ or $\{x_4, x_7\}$,
\myitem{(v)}$G$ has a subgraph isomorphic to the second graph shown in 
Figure~\ref{fig:criticalgraphs}, 
and the only pairs of vertices of $C$ colored the same  are 
exactly $\{x_2,x_6\}$ and $\{x_3,x_7\}$,
\myitem{(vi)}$G$ has a subgraph isomorphic to the third graph shown in 
Figure~\ref{fig:criticalgraphs}, 
and the only pairs of vertices of $C$ colored the same  are
exactly $\{x_2,x_6\}$ and $\{x_3,x_7\}$.
\end{lemma}

\proof
Clearly, if one of (i)--(vi) holds, then $c$ cannot be extended to
a $5$-coloring of $G$. To prove the converse we will show, by induction on
$|V(G)|$, that if none of (i)--(vi) holds, then $c$ can be extended to
a $5$-coloring of $G$. Since $c$ extends if 
$|V(G)| \le 4$, we assume that
$|V(G)| \ge 5$, and that the lemma holds for all graphs on fewer
vertices. We may also assume that $V(G)\ne V(C)$, and that every
vertex of  $G-V(C)$ has degree at least five, for  we
can delete a vertex of $G-V(C)$ of degree at most four and
proceed by induction.
Likewise, we may assume that 

\claim{$*$}{%
  the graph $G$ has no cycle of length at most four whose 
  removal disconnects $G$.
}

\noindent This is because if a cycle $C'$ of length at most four 
separates $G$, then we first
delete all vertices and edges drawn in the open disk bounded by $C'$
and extend $c$ to that graph by induction. Then, by another application
of the induction hypothesis we extend the resulting coloring of $C'$ to
a coloring of the entire graph $G$.
Thus we may assume ($*$).
\medskip

Let $v$ be a vertex of $G-V(C)$ joined to $m$ vertices of $C$, where $m$ is
as large as possible. 
Then we may assume that $m\ge3$, for otherwise $c$ extends to 
a $5$-coloring of $G$ by the Theorem of~\cite{Tho5choose}.

%We claim that if $m\le2$, then the lemma holds. 
%To prove the claim assume 
%that every vertex of $G-V(C)$ has at most two neighbors
%in $C$. We deduce that some vertex of $G-V(C)$ has degree at most
%five, for otherwise $G-V(C)$ has at least five vertices,
%contrary to~\cite[Lemma~3.1]{Tho5torus}, because $k\le7$.
%This shows that  $G-V(C)$ has a  vertex $u$ of degree five.
%Since $G$ has no separating triangle we deduce that
%the vertex $u$ has two neighbors $u_1$, $u_2$, which are not adjacent
%and not on $C$. Let $J$ be obtained from $G$ by deleting $u$,
%identifying $u_1$ and $u_2$, and deleting all resulting parallel
%edges. Since $m\le2$ the graph $J$ satisfies none of (i)--(vi),
%and hence the coloring $c$ can be extended to a $5$-coloring of
%$J$ by the induction hypothesis.
%It follows that the coloring $c$ can be extended to $G$,
%as desired.
%Thus we may assume that $m\ge3$.

Since (i) does not hold, the coloring $c$ extends to a $5$-coloring $c'$
of the graph $G':=G[V(C)\cup\{v\}]$.
Let $D$ be a facial cycle of $G'$ other than $C$, and let
$H$ be the subgraph of $G$ consisting of $D$ and all vertices and
edges drawn in the disk bounded by $D$.
If $c'$ extends to $H$ for every choice of $D$, then $c$ extends
to $G$, and the lemma holds.
We may therefore assume that $D$ was chosen so that $c'$ does
not extend to $H$.
By the induction hypothesis $H$ and $D$ satisfy one of (i)--(vi).

If $H$ and $D$ satisfy (i), then there is a vertex $w\in V(H)-V(D)$
that sees five distinct colors on $D$.
Thus $w$ has at least four neighbors on $C$, and hence $m\ge4$.
It follows that every bounded face of the graph $G[V(C)\cup\{v,w\}]$
has size at most four, and hence $V(G)=V(C)\cup\{v,w\}$ by ($*$).
Since (i) and (ii) do not hold for $G$, we deduce that $c$ can be extended
to a $5$-coloring of $G$, as desired.

If $H$ and $D$ satisfy (ii), then there are adjacent vertices
$v_1,v_2\in V(H)-V(D)$ that see the same four colors on $D$.
It follows that $m\ge3$, and similarly as in the previous paragraph
we deduce that $V(G)=V(C)\cup\{v,v_1,v_2\}$.
It follows  that $c$ can be extended to a $5$-coloring of $G$:
if both $v_1$ and $v_2$ are adjacent to $v$ we use that
$G$ does not satisfy (i), (ii), or (iii); otherwise we use that
$G$ does not satisfy (i), (ii), or (iv).

If $H$ and $D$ satisfy (iii), then there are three pairwise adjacent vertices
of $v_1,v_2,v_3\in V(H)-V(D)$ that see the same three colors on $D$.
It follows in the same way as above that $V(G)=V(C)\cup\{v,v_1,v_2,v_3\}$.
If $v$ sees at most three colors on $C$, then $c$ extends to
a $5$-coloring of $G$, 
because there are at least two choices for $c'(v)$.
Thus we may assume that $v$ sees at least
four colors. It follows that $m=4$, because $k\le7$.
Since $G$ does not satisfy (v) or (vi) we deduce that $c$
extends to a $5$-coloring of $G$.

If $H$ and $D$ satisfy (iv), then there are three  vertices of
%$v_1,v_2,v_3\in 
$V(H)-V(D)$ forming the first subgraph in Figure~\ref{fig:criticalgraphs}.
But at least one of these vertices has four neighbors on $C$,
and hence $m\ge4$, contrary to $k\le7$.

Finally, if  $H$ and $D$ satisfy (v) or (vi), then $H$ has
a subgraph isomorphic to the second or third graph depicted in
Figure~\ref{fig:criticalgraphs}, 
and the restriction of $c'$ to $D$ is uniquely determined
(up to a permutation of colors). 
Since $D$ has length seven, it follows that $m\le3$,
and hence $c'(v)$ can be changed to a different value, contrary to
the fact that the restriction of $c'$ to $D$ is uniquely determined.~\qed

The following lemma is shown in~\cite{Sas6reg}.

\begin{lemma} 
\mylabel{6regular-5colorable}
All $6$-regular graphs embeddable on the Klein bottle are $5$-colorable.
\end{lemma}

The next lemma is an adaptation of~\cite[Lemma~5.2]{Tho5torus} for
the Klein bottle.

\begin{lemma}
\mylabel{thomlem5.2}
Let $G$ be isomorphic to $C_3 + C_5$, let $S$ be a cycle in $G$ of
length three with vertex-set  $\{z_0, z_1, z_2\}$,
and let $u_1$ be a vertex in $G\backslash V(S)$ adjacent to $z_0$. 
Let $G'$ be obtained from $G$ by splitting $z_0$ into two nonadjacent 
vertices $x$ and $y$ such that $u_1$ and at most one more vertex $u_0$ 
in $G'$ is adjacent to both $x$ and $y$ 
and such that $yz_1z_2x$ is a path in $G'$. 
Let $G''$ be obtained from $G'$ by adding a vertex $v_0$ and joining 
$v_0$ to $x, y, u_1, z_1, z_2$. 
If $G''$ is not $5$-colorable and can be drawn in the Klein bottle,
then either $G'\setminus x$ or $G'\setminus y$ has a subgraph isomorphic to $C_3 + C_5$ or $G''$ is isomorphic to $\numbernine$.
\end{lemma}

\proof
We follow the argument of~\cite[Lemma~5.2]{Tho5torus}.
%except that instead of invoking~\cite[Proposition~2.3]{Tho5torus}
%on two occasions we use the fact that in those cases the graph $G''$
%is isomorphic to~$\numbernine$.~\qed
%Suppose for a contradiction that $G''$ 
%has no subgraph isomorphic to $ \numbernine $.
If one of $x, y$ has the same neighbors in $G'$ as $z_0$
does in $G$, say $x$, then $G'\setminus y$ has a subgraph isomorphic to $C_3 + C_5$,
as desired.  
%Then $G' \not \supseteq \numbernine$ as $G'$ contains
%only nine vertices, and $G' \not \supseteq C_3 + C_5$ by inspection.
Thus we can assume that $z_0$ has  two neighbors in $G$
such that one is a neighbor in $G'$ of $x$ but not $y$ and the other is a
neighbor in $G'$ of $y$ but not $x$.

The vertices $x,y$ have degree at least five in $G''$, for if
say $y$ had degree at most four, then $G''\backslash y\backslash v_0$
would not be $5$-colorable (because $G''$ is not), 
and yet it is a proper subgraph
of $C_3 + C_5$, a contradiction.
It follows that $z_0$ has degree at least six in $G$.  
%Suppose that
%$x$ had degree at most four in $G''$.  Then $G'' - \{x, v_0\}$ is a
%proper subgraph of $C_3 + C_5$ as $y$ is not a neighbor of one of
%x's neighbors.  Since $G'' - \{x, v_0\}$ is a proper subgraph, we
%can properly 5-color it, and extend this to a 5-coloring of $G''$ by
%coloring $v_0$, then $x$.
Let $G$ consist of a 5-cycle $p_1p_2p_3p_4p_5p_1$ and a $3$-cycle
$q_1q_2q_3q_1$ and the $15$ edges $p_iq_j$ where $1 \leq i \leq 3, 1
\leq j \leq 5$.  Since the degree of $z_0$ in $G$ is at least 6,
we have $z_0 \in \{q_1, q_2, q_3\}$.
The remainder of the proof is an analysis based on which vertices
are $z_0, z_1, z_2$.  

First suppose that $z_0,z_1,z_2$ are $q_3,
q_1, q_2$, respectively.  If both $u_0$ and $u_1$ are in $\{p_1,
p_2, p_3,\allowbreak p_4, p_5\}$, then we can color $y, z_1, z_2, x$ with $2,
1, 2, 1$, respectively.  We can color the remaining vertices with
colors $3, 4, 5$ as the remaining vertices are $v_0$ and a 5-cycle,
and in this case $v_0$ is only adjacent to one of the vertices of
the 5-cycle.  If $u_1 = p_1$ and $u_0 = z_1$, then we color $y, z_1,
z_2, x, u_1$ by $2, 1, 2, 3, 4$, respectively.  Since 
some neighbor of $z_0$ in $G$ is not a neighbor of $y$ in $G''$,
%$y$ has degree at least five in $G''$, 
some vertex in $\{p_2, p_3, p_4, p_5\}$ can
obtain color 3 and the remaining vertices may be colored with colors
4 and 5.

Now consider the case where $z_0, z_1, z_2$ are $q_1, p_1, p_2$,
respectively and $u_0$ is not in $\{z_1, z_2\}$.  Color $y, z_1,
z_2, x, u_0, u_1$ by $2, 1, 2, 1, 3, 4$, respectively.  We can
extend this to a 5-coloring of $G''$, coloring $v_0$ last, except 
(up to symmetry) in the following three cases.  If $u_0
= q_2$ and $u_1 = p_4$, color $q_3$ by the same color as $x$ or $y$
and recolor either $z_1$ or $z_2$ by 4 and color the remaining
vertices color 5.  If $u_0 = p_3$ and $u_1 = p_4$, then color $q_3$
by color 1 or 2 and recolor $z_1$ or $z_2$ color 4.  Then we can
color $p_5, q_2$ with colors 3 and 5 respectively.  If $u_1 = p_3$
and $u_0 = p_5$, color $q_3$ by 1 or 2 and recolor one of $z_1, z_2$
by 3 and recolor $p_3, p_4, p_5, q_2$ by $4, 3, 4, 5$, respectively.

Now suppose that $z_0, z_1, z_2$ are $q_1, p_1, p_2$, respectively
and $u_0$ is in $\{p_1, p_2\}$.  Without loss of generality let $u_0
= p_2$.  Suppose that $u_1 \in \{p_3, p_4, p_5\}$.  Then we can color
$y, z_1, z_2, x$ by $2, 4, 3, 1$ and can color $u_1$ by $3$, except
 when $u_1 = p_3$, in which case we color $u_1$ by 4. Next,
color one of $q_2, q_3$ color 1 or 2.  If both $q_2, q_3$ can be
colored 1, 2 then the rest of the coloring follows.  So we can
assume that $q_2, q_3$ are colored by $2, 5$, respectively and both
$q_2, q_3$ are adjacent to $x$.  (The argument is analogous if $q_2,
q_3$ are both adjacent to $y$.)  Since $y$ has degree at least four
in $G'$, at least one vertex in $\{p_3, p_4, p_5\} \setminus
\{u_1\}$ is joined to $y$ and is colored 1.  With possibly a
swapping of the colors of $z_1$ and $z_2$, we can now complete the
5-coloring.

Suppose that $z_0, z_1, z_2$ are $q_1, p_1, p_2$, respectively and
$u_0 = p_2$ and $u_1 = q_2$.  Color $y, z_1, z_2, x, q_2$ colors
$2,1,3,1,4$, respectively.  If $q_3$ can be colored $2$, then color
$p_3, p_4, p_5, v_0$ colors $5, 3, 5, 5$, respectively.  So we may assume that
$q_3$ is adjacent to $y$.  Then color $q_3$ by 5.  If we can color
$\{p_3, p_4, p_5\}$ by colors $\{1, 2, 3\}$, then color $v_0$ with
5. If not, then $p_3, p_4$ are adjacent to the same vertex in $\{x,
y\}$.  Since $x$ has degree at least four in $G'$ and only $u_0,u_1$
are adjacent to both $x,y$, that vertex must
be $x$. We may assume that $p_5$ is adjacent to $y$ since otherwise we color
$p_3, p_4, p_5$ by $2, 3, 2$, respectively.
It follows that $G''$ is isomorphic to $\numbernine$
by an isomorpism that maps the vertices 
$z_1, y, q_3, z_2, q_2, p_5, p_4, p_3, x, v_0$
to the vertices of $\numbernine$ in order, where the vertices of $\numbernine$
are numbered by reference to Figure~\ref{fig:allels}, starting at 
top left and moving horizontally to the right one row at a time.
%from  left to right and then top-down, order given 
%by Figure~\ref{fig:allels}.

%XXXX

%isomorphism, we first will label the vertices of the figure of
%$\numbernine$.  We will label them from top to bottom, breaking ties
%from left to right.  So the vertex on the top left is labeled $a$,
%the vertex on the top right is labeled $b$.  The next row of
%vertices is labeled $c$ and $d$, respectively.  The vertices in the
%third row from the top is labeled $e$, and the vertices on the
%bottom row are labeled $f, g, h, i, j$ from left to right.  With
%this notation, we have the following correspondence: $a \equiv z_1,
%b \equiv y, c \equiv q_3, d \equiv z_2, e \equiv q_2, f \equiv p_5,
%g \equiv p_4, h \equiv p_3, i \equiv x, j \equiv v_0$.  Thus $G''$
%is isomorphic to $\numbernine$.

Now, consider the case when $z_0, z_1, z_2$ are $q_1, q_2, p_1$,
respectively.  If $u_0 \not \in \{z_1, z_2\}$, then color $y, z_1,
z_2, x, p_2, p_3, p_4, p_5, q_3$ by $2, 1, 2, 1, 3, 4, 3, 4, 5$,
respectively.  If $u_0 = p_1$, color $y, z_1, z_2, x$ by $2, 1, 3
,1$, respectively.  If $q_3$ is not adjacent to $y$, then color
$q_3$ by $2$ and the vertices $p_2, p_3, p_4, p_5$ colors 4 and 5.  If
$q_3$ is adjacent to $y$, color $q_3$ by 5. Since $x$ has degree at
least four in $G'$, some vertex in $\{p_2, \ldots, p_5\}$ can be
colored 2.  The other vertices in this set could then be colored
with colors 3 and 4.  Thus assume that $u_0 = q_2 = z_1$.  Color $y,
z_1, z_2, x, u_1$ by $2,3, 2, 1, 4$ and we now will try and extend
this coloring.  If $q_3$ can be colored 1, then color $p_2, p_3,
p_4, p_5$ by colors 4 and 5.  So we assume that $q_3$ is adjacent to
$x$.  If $u_1 = p_3$, then recolor $z_2$ by color 4 and color $q_3$
by 2.  Since $y$ also has degree at least four in $G'$, it must be
adjacent to at least one of $p_4, p_5$, which we color 1.  The
remaining vertices of $\{p_1, \ldots, p_5\}$ are colored 5.  If $u_1
= q_3$, then we color one of $p_2$ or $p_5$ color 1 if possible and
complete the coloring by using 5 for two vertices in $\{p_2, p_3,
p_4, p_5\}$.  Now assume that both $p_2$ and $p_5$ are joined to
$x$. Since $y$ has degree at least four in $G'$ it follows that $y$
is adjacent to $p_3$ and $p_4$.  We now claim that $G''$ is not
embeddable on the Klein bottle.  Notice that if an embedding of this
graph exists, it must be that it is a triangulation as it has 10
vertices and 30 edges.  Consider the induced embeddings of 
$G'' \backslash p_2, G''\backslash p_5$ and $G'' \backslash v_0$, 
respectively.  The face of $G'' \backslash
p_2$ containing $p_2$ is bounded by a Hamiltonian cycle of $N_{G''}(p_2)$.
There exist similarly constructed
Hamiltonian cycles in $N_{G''}(p_5)$ and $N_{G''}(v_0)$.  
However, each of these
cycles contains the edge $xp_1$.  This would mean that $xp_1$ is
part of three facial triangles, a contradiction.

 Finally, consider
the subcase where $z_0, z_1, z_2, u_0, u_1$ are $q_1, q_2, p_1, q_2,
p_2$, respectively.  Color $y, z_1, z_2, x, u_1, q_3$ by $2, 3, 2,
1, 4 ,5$,  respectively.  We may assume that $q_3$ is adjacent to
$x$ else we can recolor $q_3$ by 1 and complete the coloring.  Also,
we can assume that $p_5$ is adjacent to $x$ else we color $p_5,
p_4$, by $1, 4$ and complete the coloring. Color $p_5$ by 4.  The
coloring can be completed unless $p_3$ and $p_4$ are both adjacent
to the same vertex in $\{x, y\}$.  Since $y$ must have degree at
least four in $G'$, it follows that $p_3$ and $p_4$ are adjacent to
$y$.  It follows that $G''$ has a subgraph isomorphic to $\numbernine$
by an isomorphism that maps
$x, z_2 = p_1 ,  q_3, p_2 = u_1,  q_2 =z_1 = u_0,  p_5, p_4,  p_3,  y, v_0$
to the vertices of $\numbernine$ in order, using the same numbering of
the vertices of $\numbernine$ as above.
%. Using the notation
%from above we have the following correspondence: $a \equiv x, b
%\equiv z_2 = p_1 , c \equiv q_3, d \equiv p_2 = u_1, e \equiv q_2 =
%z_1 = u_0, f \equiv p_5, g \equiv p_4, h \equiv p_3, i \equiv y, j
%\equiv v_0$. 
Thus $G''$ is isomorphic to $\numbernine$.~\qed

We also need a minor variation of the previous lemma, a case not
treated in~\cite{Tho5torus}.

\begin{lemma}
\mylabel{thomlem5.2b}
Let $G$ be isomorphic to $C_3 + C_5$, let $S$ be a cycle in $G$ of
length three with vertex-set  $\{z_0, z_1, z_2\}$,
and let $u_1$ be a vertex in $G\backslash V(S)$ adjacent to $z_0$. 
Let $G'$ be obtained from $G$ by adding an edge between two nonadjacent 
vertices neither of which is $z_0$,
and then splitting $z_0$ into two nonadjacent vertices $x$ and $y$ 
such that $u_1$ is the only vertex in $G'$ that is adjacent to both 
$x$ and $y$ and such that $yz_1z_2x$ is a path in $G'$.
Let $G''$ be obtained from $G'$ by adding a vertex $v_0$ and joining 
$v_0$ to $x, y, u_1, z_1, z_2$. 
If $G''$ is not $5$-colorable and can be drawn in the Klein bottle, then either $G'\setminus x$ or $G'\setminus y$ has a subgraph isomorphic to either $C_3 + C_5$ or $K_6$.
\end{lemma}

\proof
%Suppose for a contradiction that $G''$ 
%has no subgraph isomorphic to $C_3 + C_5$.
If one of $x, y$ has the same neighbors in $G'$ as $z_0$
does in $G$, say $x$, then $G'\setminus y$ has a subgraph isomorphic to $C_3 + C_5$,
as desired.
Thus we can assume that $z_0$ has  two neighbors in $G$
such that one is a neighbor in $G'$ of $x$ but not $y$ and the other is a
neighbor in $G'$ of $y$ but not $x$.
We may assume that the vertices $x,y$ have degree at least five in $G''$, 
for if
say $y$ had degree at most four, then $G''\backslash y\backslash v_0=G'\setminus y$
would not be $5$-colorable (because $G''$ is not), yet this is a proper subgraph
of $C_3 + C_5$ plus an additional edge, and hence by Lemma~\ref{smallcrit}
must contain either $C_3+C_5$ or $K_6$ as a subgraph, as desired.
Moreover, the sum of the degrees of $x$ and $y$ in $G''$
is at most $10$ since $z_0$ has degree at most seven in $G$. 
Thus, $z_0$ must have degree seven in $G$ while $x$ and $y$
must have degree five in $G''$.

Let $G$ consist of a 5-cycle $p_1p_2p_3p_4p_5p_1$ and a $3$-cycle
$q_1q_2q_3q_1$ and the $15$ edges $p_iq_j$ where $1 \leq i \leq 3, 1
\leq j \leq 5$.  Since the degree of $z_0$ in $G$ is seven, we have $z_0 \in \{q_1, q_2, q_3\}$.
Without loss of generality, let $z_0=q_1$. Moreover, in $G'$, there is an edge between two of the $p$'s
that are not adjacent in $G$. Without loss of generality, suppose that this edge is $p_1p_3$.

As $u_1$ is the only vertex in $G'$ adjacent to both $x$ and $y$, we have that $x$ and $z_1$ are not adjacent.
Consider the graph $G_{xz_1}$ obtained from $G''$ by deleting $v_0$,
identifying $x$ and $z_1$ into a new vertex $w$, and deleting parallel edges. 
Now $G_{xz_1}$ must not be $5$-colorable, as otherwise we could color $G''$. 
Now $G_{xz_1}$ must contain a $6$-critical subgraph $H$. 
As $y$ has degree at most four in $G_{zx_1}$, $y$ is not in $H$. 
Thus $|V(H)| \le 7$. 
By Lemma~\ref{smallcrit} we find that $H$ is isomorphic to $K_6$. 
The vertex $w$ must be in $H$ as otherwise $G\backslash x$ 
would contain $K_6$ as a proper subgraph, a contradiction. 
The remaining five vertices of $H$ induce a $K_5$. 
So these vertices must be $q_2, q_3, p_1, p_2, p_3$. 
Hence $z_1$ must be one of $p_4$ or $p_5$.

A similar argument shows that $y$ and $z_2$ are not adjacent and that 
the analogously defined graph $G_{yz_2}$ must contain a subgraph $H'$ 
isomorphic to $K_6$ with vertices $q_2, q_3, p_1, p_2, p_3$ and the new
vertex of $G_{yz_2}$. 
Hence $z_2$ must be one of $p_4$ or $p_5$. 
Without loss of generality, suppose that $z_1=p_4$ and $z_2=p_5$. 
As there are edges between $w$ and $p_1, p_2$, the edges $xp_1$ and $xp_2$ 
must be present in $G'$. 
Similarly, the edges $yp_3$ and $yp_2$ must be in $G'$. Hence $u_1=p_2$. 
Finally, as $x$ and $y$ have degree four in $G'$ and exactly one of $p_4=z_1$
and $p_5=z_2$ is adjacent to $x$ and exactly one is adjacent to $y$, 
we may assume without loss of generality that $x$ is adjacent to 
$q_2$ and $y$ is adjacent to $q_3$. 

It is straightforward to color $G''$. Color $q_2$ and $y$ with color $5$; 
color $q_3$ and $x$ with color $4$. Color $p_2$ and $p_4$ with color $1$. 
Color $p_3$ and $p_5$ with color $3$. Color $p_1$ and $v_0$ with color $2$.
This $5$-coloring of $G''$ contradicts the hypothesis of the lemma.~\qed

We also need an adaptation of~\cite[Lemma~5.3]{Tho5torus} for
the Klein bottle. We leave the similar proof to the reader.

\begin{lemma}
\mylabel{thomlem5.3}
Let $G$ be isomorphic to $K_2 + H_7$, let $S$ be a cycle in $G$ of
length three with vertex-set  $\{z_0, z_1, z_2\}$,
and let $u_1$ be a vertex in $G\backslash V(S)$ adjacent to $z_0$.
Let $G'$ be obtained from $G$ by splitting $z_0$ into two nonadjacent
vertices $x$ and $y$ such that $u_1$ and at most one more vertex $u_0$
in $G'$ is joined to both $x$ and $y$ and such that $yz_1z_2x$ is a path in $G'$
.
Let $G''$ be obtained from $G'$ by adding a vertex $v_0$ and joining
$v_0$ to $x, y, u_1, z_1, z_2$.
If $G''$ is not $5$-colorable and can be drawn in the Klein bottle,
then $G'\setminus x$ or $G'\setminus y$ has a subgraph isomorphic to  $K_2+H_7$.
\end{lemma}

We also need a similar variation of the previous lemma to handle a case
not treated in~\cite{Tho5torus}.

\begin{lemma}
\mylabel{thomlem5.3b}
Let $G$ be isomorphic to $K_2 + H_7$, let $S$ be a cycle in $G$ of
length three with vertex-set  $\{z_0, z_1, z_2\}$,
and let $u_1$ be a vertex in $G\backslash V(S)$ adjacent to $z_0$. 
Let $G'$ be obtained from $G$ by adding an edge between two nonadjacent 
vertices neither of which is $z_0$,
and then splitting $z_0$ into two nonadjacent vertices $x$ and $y$ such 
that $u_1$ is the only vertex in $G'$ that is adjacent to both $x$ and $y$ 
and such that $yz_1z_2x$ is a path in $G'$.
Let $G''$ be obtained from $G'$ by adding a vertex $v_0$ and joining 
$v_0$ to $x, y, u_1, z_1, z_2$. 
If $G''$ is not $5$-colorable and can be drawn in the Klein bottle, then either $G'\setminus x$ or $G'\setminus y$ has a subgraph isomorphic to either $K_2 + H_7$ or $K_6$.
\end{lemma}

\proof
%Suppose for a contradiction that $G''$ 
%has no subgraph isomorphic to $K_2+H_7$.
If one of $x, y$ has the same neighbors in $G'$ as $z_0$
does in $G$, say $x$, then $G'\setminus y$ has a subgraph isomorphic to 
$K_2+H_7$, as desired.
Thus we can assume that $z_0$ has two neighbors in $G$
such that one is a neighbor in $G'$ of $x$ but not $y$ and the other is a
neighbor in $G'$ of $y$ but not $x$.
The vertices $x,y$ have degree at least five in $G''$, for if
say $y$ had degree at most four, then 
$G''\backslash y\backslash v_0 =G'\setminus y$
would not be $5$-colorable (because $G''$ is not), 
and yet this is a proper subgraph
of $K_2+H_7$ plus an additional edge and by Lemma~\ref{smallcrit} 
must contain $K_2+H_7$, $C_3+C_5$ or $K_6$ as a subgraph, a contradiction. 
%(Note that it cannot contain $C_3+C_5$ as such cannot be obtained
%from a proper subgraph of $K_2+H_7$ by adding an edge, not that this really matters??)

Hence $x$ and $y$ have degree at least four in $G'$ and so $z_0$ has degree at least seven in $G$. Moreover, the sum of the degrees of $x$ and $y$ is at most $9$ since $z_0$ has degree at most eight in $G$. Thus, $z_0$ must have degree eight in $G$. 
Without loss of generality we may assume that $x$ has degree five and 
$y$ has degree four in $G'$.

We label $K_2+H_7$ as follows. The two degree eight vertices are $q_1,q_2$. 
The degree six vertex is $p_1$. 
The degree fives are $p_2,p_3,p_4, p_5, p_6, p_7$, 
where $p_2,p_3,p_4$ and $p_5,p_6,p_7$ are triangles, $p_4p_5$ is an edge, 
and $p_2,p_3,p_6,p_7$ are adjacent to $p_1$.   
Since the degree of $z_0$ in $G$ is eight, we have $z_0 \in \{q_1, q_2\}$.
Without loss of generality, let $z_0=q_1$. 
Moreover, in $G'$, there is an edge between two of the $p$'s
that are not adjacent in $G$. 

As $u_1$ is the only vertex in $G'$ adjacent to both $x$ and $y$, we have that $x$ and $z_1$ are not adjacent.
Consider the graph $G_{xz_1}$ obtained from $G''$ by deleting $v_0$,
identifying $x$ and $z_1$ into a new vertex $w$, and deleting parallel edges. 
Now $G_{xz_1}$ must not be $5$-colorable, as otherwise we could $5$-color $G''$.
Thus $G_{xz_1}$ must contain a $6$-critical subgraph $H$. 
As $y$ has degree at most four in $G_{zx_1}$, $y$ is not in $H$. 
Thus $|V(H)| \le 8$. By  Lemma~\ref{smallcrit} we find that $H$ is 
isomorphic to $K_6$ or $C_3 + C_5$. 
The vertex $w$ must be in $H$ as otherwise $G''$ would contain a 
proper subgraph that is not $5$-colorable, a contradiction. 

Let $J=G'\setminus \{x,y,z_1\}$. If $H$ is isomorphic to $C_3+C_5$, 
then $J$ must contain a subgraph isomorphic to $K_2+C_5$,
because $q_2$ and $p_1$ are the only vertices of 
$G\backslash z_0=G\backslash q_1$
that could have degree at least six.
Thus there must be two degree six vertices in $J$. These must be $q_2$ and $p_1$. The other five vertices must be neighbors of $p_1$ and yet must form a $C_5$. This is impossible. So $H$ must be isomorphic to $K_6$. Now $J$ must contain $K_5$ as a subgraph. This can only happen if one of the edges $p_1p_4$ or $p_1p_5$ is present in $G'$. Without loss of generality suppose that $p_1p_4$ is present in $G'$. 
Then $H$ must consist of the vertices $w, q_2, p_1, p_2, p_3, p_4$. 
So $z_1$ must be one of $p_5, p_6, p_7$. 
It follows that $x$ is adjacent to $p_2$ and $p_3$.

Similarly, as $u_1$ is the only vertex in $G'$ adjacent to both $x$ and $y$, we have that $y$ and $z_2$ are not adjacent.
Consider the graph $G_{yz_2}$ obtained from $G''$ by deleting $v_0$, 
identifying $y$ and $z_2$ into a new vertex $w'$, and deleting parallel edges. 
Now $G_{yz_2}$ must not be $5$-colorable, as otherwise we could $5$-color $G''$.
Now $G_{yz_2}$ must contain a $6$-critical subgraph $H'$. 
Thus $|V(H')| \le 9$. By Lemma~\ref{smallcrit} we find that $H'$ is 
isomorphic to $K_6$, $C_3 + C_5$, or $K_2+H_7$. The vertex $w'$ must be 
in $H$ as otherwise $G''$ would contain a proper subgraph that is not $5$-colorable, a contradiction. 

Suppose that $x$ is not $H'$. 
The previous argument for $H$ shows that $H'$ is isomorphic to $K_6$, 
that $H'$ consists of $w', q_2,p_1,p_2,p_3,p_4$ and that $y$ is adjacent to $p_2$ and $p_3$. But then there are two vertices, $p_2$ and $p_3$, adjacent to both $x$ and $y$, a contradiction. 

So  $x$ is in $H'$. Now the neighbors of $x$ must be in $H'$. 
Specifically, $p_2$ and $p_3$ are in $H'$. 
Note that $p_2$ and $p_3$ are not equal to $z_2$ as they are not 
adjacent to $p_5$, $p_6$ or $p_7$. 
Meanwhile, at least one of $p_2,p_3$ is not adjacent to $y$. 
Without loss of generality, suppose that $p_2$ is not adjacent to $y$. 
Now $p_2$ has degree five in $G'$ and hence degree at most five in $G_{yz_2}$. 
Thus the neighbors in $G'$ of $p_2$ and all edges incident in $G'$ 
with $p_2$ must be in $H'$. 

If $H'$ is isomorphic to $K_6$, then it follows that $x$ must be adjacent 
to all of $H\setminus w'$ as well as $z_2$. 
That is, $x$ must be adjacent to all the neighbors of $p_2$, namely $q_2,p_1,p_3,p_4$. Now $G'$ contains $K_6$ as a subgraph, a contradiction. If $H'$ is isomorphic to $C_3+C_5$, then $xp_2p_3$ is a triangle in $H'$. Thus one of these vertices must have degree seven in $H'$. However, $x$ and $p_2$ have degree five in $G'$ while $p_3$ has degree at most six, a contradiction. 

Thus $H'$ is isomorphic to $K_2+H_7$. As $H'$ has nine vertices, $z_1$ must be in $H'$ and have degree five. As $z_1$ is adjacent to $y$ but not adjacent to $x$, $z_1$ has degree five in $G'$.  However, $z_1$ is adjacent to $z_2$. So $z_1$ has degree four in $G_{yz_2}$ and so has degree at most four in $H'$, a contradiction.~\qed

\begin{lemma}
\mylabel{2curve}
Let $G$ be a graph drawn in the Klein bottle, and let $c,d\in V(G)$
be such that 
%$G\backslash\{c,d\}$ is connected, 
$G\backslash c$ does not embed in the projective plane,
and $G$ does not embed in the torus.
Then every closed curve in the Klein bottle intersecting $G$ in
a subset of $\{c,d\}$ separates the Klein bottle.
\end{lemma}

\proof
Let $\phi$ be a closed curve in the Klein bottle intersecting $G$ in
a subset of $\{c,d\}$, and suppose for a contradiction that it
does not separate the Klein bottle.
Then $\phi$ is either one-sided or two-sided.
If  $\phi$ is one-sided, then it intersects $G\backslash c$ in at most
one vertex, and hence the Klein bottle drawing of $G\backslash c$ can
be converted into a drawing of  $G\backslash c$ in the projective plane,
a contradiction.
Thus $\phi$ is two-sided, but then the drawing of $G$ can be
converted into a drawing of $G$ in the torus, again a contradiction.~\qed

\begin{figure}
\centering
\includegraphics[scale=1.2]{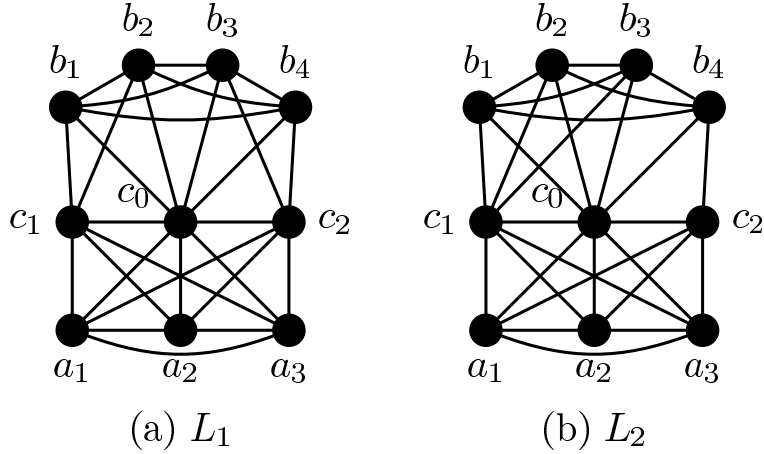}
\caption{The graphs $L_1$ and $L_2$ with their vertices labeled}
\label{fig:l1l2}
\end{figure}

\begin{lemma} 
\mylabel{facial-cycle-lemmas-10}
Let $G$ be $\tengirl$ or $\tenguy$ with its vertices numbered as in
Figure~\ref{fig:l1l2}, and let it be drawn in the Klein bottle.
Then 
\myitem{(i)}every face is bounded by a triangle, except for exactly one, which
is bounded by a cycle of length five with vertices 
$c_1, a_i, c_2, b_j, b_k$ in order for some indices $i,j,k$, and
\myitem{(ii)} for $i=0,1,2$ the vertices $a_1,a_2,a_3$ appear
consecutively in the cyclic order around $c_i$
(but not necessarily in the order listed), and so do the neighbors of $c_i$
that belong to $\{b_1,b_2,b_3,b_4\}$.
\end{lemma}

\proof
Let $i\in\{1,2\}$. There are indices $j,k$ such that $a_j$ and $b_k$
are both adjacent to $c_i$ and are next to each other
in the cyclic order around $c_i$.
Let $f_i$ be the face incident with both the edges $c_ia_j$ and $c_ib_k$.
We claim that the walk bounding $f_i$  
includes at most one occurrence of $c_i$ and no occurrence of $c_0$.
Indeed, otherwise we can  construct a simple closed curve either passing through
$f_i$ and intersecting $G$ in $c_i$ only (if $c_i$ occurs at least twice
in the boundary walk of $f_i$), or passing through
$f_i$ and a neighborhood of the edge $c_ic_0$ and intersecting $G$
in $c_i$ and $c_0$ (if $c_0$ occurs in the boundary walk of $f_i$).
By Lemma~\ref{2curve} this simple closed curve separates the Klein bottle.
It follows from the construction that it also separates $G$, contrary
to the fact that $G\backslash\{c_i,c_0\}$ is connected.
This proves our claim that the walk bounding $f_i$ includes at most 
one occurrence of $c_i$ and no occurrence of $c_0$.

Since the boundary of $f_i$ includes a subwalk from $a_j$ to $b_k$ 
that does not use $c_i$,
we deduce that $c_{3-i}$ belongs to the facial walk bounding $f_i$.
But the neighbors of $c_1$ and $c_2$ in $\{b_1,b_2,b_3,b_4\}$ are
disjoint, and hence $f_i$ has length at least five.
By Euler's formula $f_1=f_2$, this face has length exactly five, and
every other face is bounded by a triangle.
This proves (i).
Statement (ii) also follows, for otherwise there would be another face
with the same properties as $f_1=f_2$, and yet we have already shown
that this face is unique.~\qed

\begin{figure}
\centering
\includegraphics[scale=1.2]{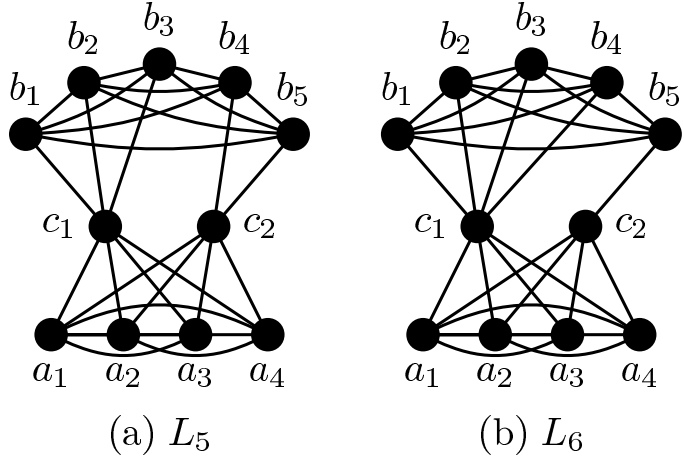}
\caption{The graphs $L_5$ and $L_6$ with their vertices labeled}
\label{fig:l5l6}
\end{figure}

\begin{lemma}
\mylabel{facial-cycle-lemmas-11}
Let $G$ be $\elevengirl$ or $\elevenguy$ with its vertices numbered as in
Figure~\ref{fig:l5l6}, and let it be drawn in the Klein bottle.
Then
\myitem{(i)}every face is bounded by a triangle, except for exactly two, which
are bounded by  cycles $C_1,C_2$ of length five, each with vertices
$c_1, a_i, c_2, b_j, b_k$ in order for some indices $i,j,k$, 
\myitem{(ii)} if $G=\elevengirl$, then $C_1\cap C_2$ consists of the
vertices $c_1,c_2$, and if $G=\elevenguy$, then $C_1\cap C_2$ consists of the
vertices $c_1,c_2,b_5$ and the edge $c_2b_5$, and
\myitem{(iii)} for $i=1,2$ the vertices $a_1,a_2,a_3,a_4$ appear
consecutively in the cyclic order around $c_i$
(but not necessarily in the order listed), and so do the neighbors of $c_i$
that belong to $\{b_1,b_2,b_3,b_4,b_5\}$.
\end{lemma}

\proof
 The proof is similar to the proof of Lemma~\ref{facial-cycle-lemmas-10}.
% and is omitted.
%Again, notice for $i = 1, 2$,
%Let $i\in\{1,2\}$.
There are distinct pairs $(j_1,k_1)$ and $(j_2,k_2)$ of indices such that 
$a_{j_i}$ and $b_{k_i}$ are both adjacent
to $c_1$ and are next to each other in the cyclic order around $c_1$.  
Let $f_i$ be the face incident with both $c_1a_{j_i}$ and $c_1b_{k_i}$.
We claim that the walk bounding $f_i$
includes at most one occurrence of $c_1$.  
%We will, as in Lemma~\ref{facial-cycle-lemmas-10} construct 
For if not, then there is a simple closed curve $\phi$
that passes through $f_i$ and intersects $G$ in $c_1$ only.
%at least twice in the boundary walk of $f_i$.  
But since $L_5$ and $L_6$ are not
embeddable in the torus and $L_5 \setminus c_1$ and $L_6
\setminus c_1$ are not embeddable in the projective plane, it follows from
Lemma~\ref{2curve}  that $\phi$ separates the Klein bottle.  
By construction, $\phi$ also
separates $G$, a contradiction, as $G\setminus c_1$ is connected.  
This proves our claim that the walk bounding $f_i$
includes at most one occurrence of $c_1$.
Thus the walk bounding $f_i$ includes $c_2$, and it follows similarly
that $c_2$ occurs in that walk at most once.
We deduce that $f_1$ and $f_2$ are distinct and have length at least five. 
Euler's formula implies that $f_1,f_2$ have length exactly five,
and that every other face is bounded by a triangle.
It follows that conditions (i), (ii) and (iii) hold.~\qed

%So the walk bounding $f_i$ contains at most one
%occurrence of each $c_i$.
%For each of $L_5$ and $L_6$ there are four pairs of edges that must
%be parts of faces, call them $f_1, f_2$.  These edges are, for
%$L_5$, $b_1c_1, c_1a_1; b_3c_1, c_1a_4; b_4c_2, c_2a_1; b_5c_2,
%c_2a_4$, and for $L_6$, $b_1c_1, c_1a_1; b_4c_1, c_1a_4; b_5c_2,
%c_2a_1; b_5c_2, c_2a_4$.  Since $c_1$ and $c_2$ are each used at
%most once in each face, there are at least two faces.  In addition,
%the neighbors of $c_1, c_2$ in $\{b_1, b_2, b_3, b_4, b_5\}$ are
%disjoint, so each of $f_1$ and $f_2$ have length at least five.
%However, by Euler's formula there are at most two faces of length
%exactly five, and every other face is bounded by a triangle.  This
%proves (i).
%Statement $(ii)$ follows for $L_5$ since edges $b_4c_2$ and $b_5c_2$
%are in different faces and there are no multiple edges in $L_5$.
%Similarly $(ii)$ follows for $L_6$ because edge $b_5c_2$ and vertex
%$c_1$ must be in each of the two 5-faces of $L_6$.  Statement (iii)
%also follows, for otherwise there would be another non-triangular
%face other than $f_1$ and $f_2$, and yet we have already shown that
%there are at most two of these 5-faces.

\section{Reducing to $K_6$}
\mylabel{reducek6}
%%%%%%%%%%%%%%%%%%%%%%%%%%%%%%%%%%%

If $v$ is a vertex of a graph $G$, then we denote by $N_G(v)$, or simply
$N(v)$ if the graph can be understood from the context, 
the open neighborhood of the vertex $v$; that is, the subgraph of
$G$ induced by the neighbors of $v$. 
Sometimes we will use $N(v)$ to mean the vertex-set of this subgraph.
We say that a vertex $v$ in a graph $G$ embedded in a surface
has a {\em wheel neighborhood} if the neighbors of $v$ form a cycle $C$ in
the order determined by the embedding, and the cycle $C$ is null-homotopic.
(The cycle $C$ need not be induced.)

Let $G_0$ be a graph drawn in the Klein bottle such that 
$G_0$ is not $5$-colorable and has no subgraph isomorphic to any of the
graphs listed in Theorem~\ref{main}.
Let a vertex $v_0\in V(G_0)$ of degree exactly five
be chosen so that
each of the following conditions hold subject to all previous conditions:
\myitem{(i)} $|V(G_0)|$ is minimum,
\myitem{(ii)} the clique number of $N(v_0)$, the neighborhood of $v_0$, is maximum,
\myitem{(iii)} the number of largest complete subgraphs in $N(v_0)$
is maximum,
\myitem{(iv)} the number of edges in $N(v_0)$ is maximum,
\myitem{(v)} $|E(G_0)|$ is minimum,
\myitem{(vi)} the number of homotopically-trivial triangles containing $v_0$ is maximum.

\noindent
In those circumstances we say that
the pair $(G_0, v_0)$ is an {\em optimal pair}.
Given an optimal pair $(G_0, v_0)$ we say that a
pair of vertices  $v_1,v_2$ is  an $\textit{identifiable pair}$ 
if $v_1$ and $v_2$ are non-adjacent neighbors of $v_0$.
If  $v_1,v_2$ is an identifiable pair, then we define $G_{v_1v_2}$ to be the
graph obtained from $G_0$ by deleting
all edges incident with $v_0$ except $v_0v_1$ and $v_0v_2$,
contracting the edges $v_0v_1$ and $v_0v_2$
into a new vertex $z_0$, %\marginpar{Used to be $w$}
and deleting all resulting parallel edges.
This also defines a drawing of  $G_{v_1v_2}$ in the Klein bottle.
%If $G_{xy}$ is $5$-colorable, then so is $G_0$, as is easily seen.
%Thus $G_{xy}$ has a $6$-critical subgraph, say $J$.

We now introduce notation that will be used throughout the rest of the paper.
Let $G_0'$ be obtained from $G_0$ by deleting all those edges
that got deleted during the construction of $G_{v_1v_2}$.
That means all edges incident with
$v_0$ except $v_0v_1$ and $v_0v_2$ and all those edges of $G_0$ that
got deleted because they became parallel to another edge.
Thus if a vertex $v$ of $G_0$ is adjacent to both $v_1$ and $v_2$, then
$G_0'$ will include exactly one of the edges $vv_1$, $vv_2$.
Thus the edges of $G_0'\backslash v_0$ may be identified 
with the edges of $G_{v_1v_2}$,
and in what follows we will make use of this identification.
Now if $J$ is a subgraph of $G_{v_1v_2}$ with $z_0\in V(J)$, then
let $\hat J$ be the corresponding subgraph of $G_0'$; that is,
$\hat J$ has vertex-set $\{v_0,v_1,v_2\}\cup V(J)-\{z_0\}$ and
edge-set $\{v_0v_1,v_0v_2\}\cup E(J)$.
Let $\hat{R_1}$ and $\hat{R_2}$ be the two faces of $\hat J$ incident
with $v_0$, and let $R_1,R_2$ be the corresponding two faces of $J$.
We call $R_1,R_2$ the {\em hinges} of $J$.
Finally, let $\hat R$ be the face of $\hat J\backslash v_0$ containing
$v_0$.

%\newpage
%%%%%%%%XXXXX

\begin{lemma} 
\mylabel{G'notc3c5k2h7}
Let  $(G_0, v_0)$ be an optimal pair, and let $v_1,v_2$ be an identifiable pair.
Then $G_{v_1v_2}$ has no subgraph isomorphic to $C_3 + C_5$ or $K_2 + H_7$.
\end{lemma}
\proof
%This follows by using the argument of~\cite[Theorem 6.1, Claim (9)]{Tho5torus},
%using Lemma~\ref{thomlem5.2} instead of~\cite[Lemma~5.2]{Tho5torus}.~\qed
%We follow the argument of~\cite[Theorem 6.1, Claim (9)]{Tho5torus}.
Suppose for a contradiction that there exists a subgraph
$J$ of $G_{v_1v_2}$ such that $J=C_3+ C_5$ 
%is isomorphic to  $C_3+ C_5$ 
or $J=K_2 + H_7$. 
Let us recall that $z_0$ is the vertex of $G_0$ that arises from
the identification of $v_1$ and $v_2$.
Since $J$ is not $5$-colorable the choice of $G_0$ implies that $z_0\in V(J)$.
Thus we apply the notation introduced prior to this lemma.
Let $R_1,R_2$ be the hinges of $J$, let $\hat R_1$ be bounded by
the walk $v_1u_1u_2\cdots u_kv_2v_0$,
and let $\hat R_2$ be bounded by the walk
$v_2z_1z_2 \cdots z_mv_1v_0$.
Then $k,m\ge2$.
We may assume that $k \leq m$
%\{v_1, v_2\} \cap \{u_1,\ldots, u_k\} = \emptyset$ 
and that $G_0$ is drawn on the Klein
bottle such that $k + m$ is minimized. 
%{\bf IS IT USED?}
Since $|E(J)|=3|V(J)|-1$ it follows that
$J$ has exactly one face bounded by a
$4$-cycle and all other faces are bounded by  $3$-cycles.  
So $ k = 2 $ and $m \leq 3$.  
Furthermore, if $m=3$, then all faces of $\hat J$ other than 
$\hat R_1$ and $\hat R_2$ are triangles; otherwise at most one
face other than $\hat R_1$ and $\hat R_2$ is bounded by a cycle of length four.
%We will refer to such a face as an {\em exceptional face}.
%
It follows that $z_1\ne u_2$, for otherwise the cycle
$z_1z_2\ldots z_m v_1u_1$ of $G_0$ has length at most five and bounds a disk
containing $v_0$ and $v_2$, contrary to Lemma~\ref{7cycle-5coloring}.
Similarly, $u_1\ne z_m$.
Since $J$ has no parallel edges we deduce that $z_1\ne u_1$
and $u_2\ne z_m$.
It follows that the vertices $v_1,v_2,u_1,u_2,z_1,z_m$ are pairwise distinct.
However, if $m=3$, then possibly $z_2\in\{u_1,u_2\}$.
Finally, all vertices of $G_0$ are
either in $\hat J$ or inside one of the faces $\hat R_1,\hat R_2$ of $J$ by
Lemma~\ref{7cycle-5coloring}.

Next we claim that $z_0$ has degree at least six. 
Indeed, otherwise $z_0$ is contained in
the open disk bounded by a walk $w$ of $J$ of length at most six
(because $J$ has at most one face that is not a triangle).
But $W$ is also a walk in $G_0$, and the disk it bounds includes
$v_0,v_1,v_2$. 
But $v_1$ is not adjacent to $v_2$, contrary to Lemma~\ref{7cycle-5coloring}.
This proves our claim that $z_0$ has degree at least six.

We now make a couple of remarks about vertices of degree five in $J$.
If $J=C_3+ C_5$, then $J$ has five vertices of degree five,
and the neighborhood of each is isomorphic to $K_5^-$.
If $J=K_2 + H_7$, then $J$ has six vertices of degree five;
four of them have neighborhoods isomorphic to $K_5^-$
and the remaining two have neighborhoods isomorphic to $K_5-E(P_3)$.

Let us say a vertex $v$ of degree five in $J$ is {\em good} if its neighborhood in
$J$ has the property that there are at least two triangles disjoint from
any given vertex.
Thus $J$ either has five good vertices, or it has exactly four,
and they induce a matching of size two.
It follows from the definition of optimal pair that if $N(v_0)$ has
at most one triangle, then the degree of
each good vertex of $J$ must be at least six in $G_0$.

Note that if $z_0$ is a vertex of degree six in $K_2+H_7$, then all the vertices of degree five have a $K_4$ in their neighborhood disjoint from $z_0$. Hence if $N(v_0)$ does not contain a $K_4$, each vertex of degree five in $J$ must have degree at least six in $G_0$.

%Another useful fact is the following. We claim that if $v_1, v_2$ are consecutive in the cyclic order around $v_0$ and $N(v_0) \cup \{v_1v_2\}$ does not contain a $K_4$, then $v_1v_2$ is in $E(G_0)$ and $v_0v_1v_2$ is a homotopically-trivial triangle. Suppose not. If $v_1v_2$ is in $E(G_0)$ but $v_0v_1v_2$ is a homotopically-trivial triangle then redraw $v_1v_2$ such that $v_0v_1v_2$ is a homotopically trivial triangle. This new embedded graph contradicts condition (vi) of an optimal pair. So we may assume that $v_1v_2$ is not an edge of $G_0$. Consider $G_0 \cup \{v_1v_2\}$. If this graph does not contain a critical subgraph in our list, then this contradicts condition (iv) of an optimal pair. So we may suppose that $G_0 \cup \{v_1v_2\}$ contains a critical subgraph $H$ in our list. Clearly $H$ contains the edge $v_1v_2$. Moreover $v_0\in V(H)$, as otherwise $G_0-v_0$ has a $5$-coloring in which all the neighbors of $v_0$ receive different colors, as otherwise the coloring extends to $G_0$. This yields a $5$-coloring of $H$, a contradiction. Thus $N(v_0)\subseteq V(H)$ as $v_0$ has degree five. Yet by assumption $N(v_0) \cup \{v_1v_2\}$ does not contain a $K_4$. But then $H$ must be isomorphic to $L_3$ or $L_4$ as these are the only graphs in our list that contain a vertex of degree five with a neighborhood that does not contain a $K_4$. However as $L_3$ and $L_4$ are triangulations, we find that $G_0 = H-\{v_1v_2\}$ and thus $G_0$ is $5$-colorable, a contradiction. This proves the claim.

We now condition on the cases of Lemma~\ref{7cycle-5coloring} for $\hat{R}$. 
Suppose that case (i) holds. First consider the case that $m=3$. 
Let us say that two vertices of $G_0$ are 
{\em adjacent through a face $f$} of $\hat J$ if the edge joining them
lies in $f$.
We condition on the number of edges incident with $v_0$ through $\hat{R_1}$. 
Suppose there are two such edges. Hence $v_0$ is adjacent to $u_1$ and $u_2$ through $\hat{R_1}$. Further suppose that $v_0$ is adjacent to $z_2$ through $\hat{R_2}$. If $z_2$ is adjacent to $z_0$ in $J$, then without loss of generality suppose that $v_1$ is adjacent to $z_2$ but not through $\hat{R_2}$. Redrawing the edge through $\hat{R_2}$ contradicts condition (vi) of an optimal pair.

So we may assume that $z_2$ is not adjacent to $z_0$. 
It follows that $J=K_2+H_7$ and that $z_0$ must be the vertex of degree 
six in $K_2+H_7$, because that is the only vertex of degree at least
six in $C_3+C_5$ or  $K_2+H_7$ that has a non-neighbor.
However, $v_1$ is not adjacent to $u_2$ and $v_2$ is not adjacent to $u_1$. So $N(v_0)$ does not contain a $K_4$. But then the vertices of degree five in $J$ must be a subset of $\{u_1,u_2,z_2\}$, a contradiction.

%If $N(v_0)$ has at most one triangle, then the good vertices must be a subset of $\{u_1,u_2,z_2\}$, a contradiction. Hence, $N(v_0)$ has at least two triangles. As $v_1$ is not adjacent to $u_2$ and $v_2$ is not adjacent to $u_1$, we must have that $z_2$ is adjacent to $u_1$ and $u_2$ in $J$. Note that $z_1$ and $z_3$ are also adjacent to $z_2$ in $J$. Moreover, $v_1$ and $v_2$ each have exactly three neighbors in $J-z_2$. The cyclic order of $z_0$ around $J$ must be $z_3x_1u_1u_2x_2z_1$.  The only possibility is that $x_1$ and $x_2$ are adjacent good vertices of degree five. But then $x_1$ and $x_2$ are adjacent to only two of the four vertices $\{u_1,u_2,z_1,z_3\}$, the ones of degree eight. Thus, there is another face incident with $z_0$ in the embedding of $J$ besides $R_2$ which is not a triangle, a contradiction.

So we may assume without loss of generality that $v_0$ is adjacent to $z_1$ 
through $\hat{R_2}$. 
Now we may apply Lemmas~\ref{thomlem5.2}, \ref{thomlem5.2b}, 
\ref{thomlem5.3} and~\ref{thomlem5.3b} to $G_0+v_1z_1$,
where $G_0+v_1z_1$ denotes the graph obtained from $G_0$ by adding the
edge $v_1z_1$ if $v_1$ is not adjacent to $z_1$ in $G_0$ and
$G_0+v_1z_1=G_0$ otherwise. 
We find that either $G_0 \setminus v_1 \setminus v_0$ or 
$G_0+v_1z_1 \setminus v_2 \setminus v_0$ contains a subgraph $H$ 
isomorphic to $K_6, C_3+C_5$ or $K_2+H_7$, or $G_0+v_1z_1$ is isomorphic 
to $\numbernine$. In the latter case, $G_0$ is $5$-colorable or
isomorphic to $L_4$, a contradiction. 
In the former case, note that $G_0\setminus v_0$ has a proper $5$-coloring that does not extend to a $5$-coloring of $G_0$ and hence in this coloring all of the neighbors of $v_0$ must receive different colors. This yields a $5$-coloring of $H$, a contradiction.

Suppose $v_0$ is incident with exactly one edge through $\hat{R_1}$. 
Without loss of generality we may assume that $v_0$ is adjacent to 
$u_2$ through $\hat{R_1}$. Suppose that $v_0$ is adjacent to $z_1$ and $z_2$. 
If $z_2$ is adjacent to $z_0$ in $J$, we may apply 
Lemmas~\ref{thomlem5.2}, \ref{thomlem5.2b}, \ref{thomlem5.3} 
and~\ref{thomlem5.3b} to $G_0+v_1u_2$. 
We find that either $G_0 \setminus v_1 \setminus v_0$ or 
$G_0 +v_1u_2 \setminus v_2 \setminus v_0$ 
contains a subgraph $H$ isomorphic to $K_6, C_3+C_5$ or $K_2+H_7$, 
or $G_0 + v_1u_2$ is isomorphic to $\numbernine$. 
In the latter case, $G_0$ is $5$-colorable, a contradiction. In the former case, note that $G_0\setminus v_0$ has a proper $5$-coloring that does not extend to a $5$-coloring of $G_0$ and hence in this coloring all of the neighbors of $v_0$ must receive different colors. This yields a $5$-coloring of $H$, a contradiction. So we may assume that $z_2$ is not adjacent to $z_0$ in $J$. As $v_1$ is not adjacent to $z_1$ and $v_2$ is not adjacent to $z_2$, $N(v_0)$ does not contain a $K_4$. But then the verties of degree five in $J$ would have to be a subset of $\{z_1,z_2,u_2\}$, a contradiction. 

If $v_0$ is adjacent to $z_2$ and $z_3$, then a similar but easier 
argument applies as above.
Let us assume next that $v_0$ is adjacent to $z_1$ and $z_3$. 
Note that $z_0$ must have degree at least seven in $J$ for $v_1$ and 
$v_2$ to have degree at least five in $G_0$ in this case. 
As $z_0$ has degree at most eight in $J$, at least one of $v_1$ or $v_2$ 
has degree five in $G_0$. 
If $v_1$ has degree five, consider $G_{v_2z_3}$, defined as before. 
This graph contains a subgraph $H$ isomorphic to a graph listed in
Theorem~\ref{main}. 
%As $v_0$ has degree four in $G_{v_2z_3}$, 
Since $v_0\not\in V(H)$, the vertex  $v_1$ has degree four in $G_{v_2z_3}$ 
and hence $v_1\not\in V(H)$. 
It follows that the graph obtained from $H$ by deleting the new vertex
of $G_{v_2z_3}$ is a proper subgraph of $J\backslash z_0$.
Consequently, $H$ is isomorphic to a proper subgraph of $J$,
a contradiction.
%Thus $H$ is a subgraph of $J\backslash z_0$ with some additional edges 
%all incident with the same vertex $z_3$. 
%If $J=C_3+C_5$, then $H=K_6$ but this cannot happen as $C_2+C_5$ does not have a $K_5$ disjoint from a given vertex. If $J=K_2+H_7$, then $H=C_3+C_5$ or $K_6$ but this cannot happen as $K_1+H_7$ does not have a $K_5$ or $C_2+C_5$ disjoint from a given vertex. 
If $v_2$ has degree five in $G_0$, we consider $G_{v_1z_1}$ 
similarly to obtain a contradiction.

Finally suppose $v_0$ is not incident with any edge through $\hat{R_1}$. 
Hence $v_0$ is adjacent to $z_1$ $z_2$, and $z_3$ through $\hat{R_2}$. 
Then same argument as in the preceding paragraph applies.

We may assume that $m=2$. 
We may assume without loss of generality that $v_0$ is adjacent to $u_1$ 
and $u_2$ through $\hat{R_1}$ and to $z_1$ through $\hat{R_2}$. 
Now we may apply Lemmas~\ref{thomlem5.2}, \ref{thomlem5.2b}, 
\ref{thomlem5.3} and~\ref{thomlem5.3b} to $G_0+ v_1z_1$. 
We find that either $G_0 \setminus v_1 \setminus v_0$ or 
$G_0 +v_1z_1 \setminus v_2 \setminus v_0$ contains a subgraph $H$ 
isomorphic to $K_6, C_3+C_5$ or $K_2+H_7$, or $G_0 +v_1z_1$ is 
isomorphic to $\numbernine$. 
In the latter case, $G_0$ is $5$-colorable, a contradiction. In the former case, note that $G_0\setminus v_0$ has a proper $5$-coloring that does not extend to a $5$-coloring of $G_0$ and hence in this coloring all of the neighbors of $v_0$ must receive different colors. This yields a $5$-coloring of $H$, a contradiction.
This concludes the case when (i) of Theorem~\ref{7cycle-5coloring} holds.

For cases (iv)-(vi), we have that $m=3$. Note that case (vi) cannot happen as $v_0$ must be adjacent to $v_1$ and $v_2$, which are distance three on the boundary of $\hat{R}$. For cases (iv) and (v), there are in each case two possibilities, up to symmetry, as to which internal vertex is $v_0$. In all cases, it is easy to check that $N(v_0)$ is triangle-free. All vertices of degree five in $J$ must then have degree six in $G_0$, since their neighborhood in $G_0$ has a triangle and would thus contradict that $(G_0,v_0)$ is an optimal pair. To have higher degree in $G_0$, these vertices must be a subset of $\{u_1,u_2,z_1,z_2,z_3\}$. As $K_2+H_7$ has six vertices of degree five, $J=C_3+C_5$. Furthermore, $u_1,u_2,z_1,z_2,z_3$ are all distinct and induce $C_5$. We now color $G_0$ as follows. We may assume without loss of generality that $z_2$ is adjacent to $v_2$ but not to $v_1$. Color $z_1$ and $z_3$ by 1. Color $z_2$, $v_1$, and $u_2$ by 2. Color $u_1$ and $v_2$ by 3. Finally color the other two vertices of $C_3$ using colors 4 and 5. As only three colors appear on the boundary of $\hat{R}$, this coloring extends to $G_0$ by Lemma~\ref{7cycle-5coloring}, a contradiction.

Suppose that case (iii) happens. Suppose that $m=3$. We may assume without 
loss of generality that $v_0$ is adjacent to $u_1$. If $u_1 \ne z_2$, 
then $N(v_0)$ is triangle-free. Hence the good vertices of $J$ must be 
a subset of $\{u_1,z_1,z_2,z_3\}$, which do not induce a matching, a 
contradiction. 
Thus $u_1=z_2$ and we color $G_0$ by~\cite[Lemma~5.1(a)]{Tho5torus}. 
For $m=2$, case (iii) cannot happen as $v_0$ must be adjacent to $v_1$ 
and $v_2$, which are distance three on the boundary of $\hat{R}$.

So finally we may assume case (ii). 
Let $v_0'$ be the other vertex in the interior of $\hat{R}$. 
Suppose that $m=3$. Further suppose that $z_2\in\{u_1,u_2\}$. 
Note that if $N(v_0)$ has at most one triangle, then all good vertices 
of $J$ must have degree six in $G_0$. 
However, they must be a subset of $\{u_1, u_2, z_1, z_3\}$. 
Hence there are at most four good vertices in $J$ and they do not induce 
a perfect matching, a contradiction. 
So $N(v_0)$ has at least two triangles.

Suppose that one of $v_0$ or $v_0'$ is adjacent to both $z_1$ and $z_3$. 
Now $N(v_0)$ has at most one triangle unless that vertex is $v_0'$ which 
is also adjacent to $z_2$, and $v_0$ is adjacent to one of $z_1$ or $z_3$ 
through $\hat{R_2}$ as well as $z_2$ through $\hat{R}_1$. 
In that case, the hypotheses of~\cite[Lemma~5.1(c)]{Tho5torus} 
%Lemma 5.1(c) of Thomassen 
are satisfied and we can extend that coloring to a coloring of $G_0$ by 
Lemma~\ref{7cycle-5coloring}, a contradiction.

So we may suppose that neither $v_0$ or $v_0'$ is adjacent to both $z_1$ and $z_3$. Thus $v_0'$ must be adjacent to both $v_1$ and $v_2$. Without loss of generality, we may assume that one of $v_0$ or $v_0'$ is adjacent to both $z_1$ and $z_2$ thourgh $\hat{R}_2$. Now $N(v_0)$ will have at most one triangle unless $z_2=u_2$. 
Let $G$ be the graph obtained from $G_0\backslash\{v_0,v_0'\}$ 
by adding the edge $u_1z_1$. It follows that $G$ is not $5$-colorable, because every $5$-coloring of $G$ can be extended to a $5$-coloring of $G_0$ by Lemma~\ref{7cycle-5coloring}. Since $G$ has fewer vertices than $G_0$, it follows that $G$ has a subgraph $G'$ isomorphic to one of the graphs listed in Theorem~\ref{main}. 
But the edge $u_1z_1$ belongs to $G'$, because $G\setminus u_1z_1$ is $5$-colorable. 

On the other hand, we claim that the edge $u_1z_1$ belongs to no facial triangle of $G'$. Indeed, if it did, say it belonged to a facial triangle $u_1z_1q$, then either $qz_1v_2u_2u_1q$ or $qz_1z_2v_1u_1q$ would be a contractible $5$-cycle with more than one vertex in its interior, contradicting Lemma~\ref{7cycle-5coloring}. Thus $u_1z_1$ belongs to no facial triangle of $G'$. But there are only two graphs among those listed in Theorem~\ref{main} that have an embedding with an edge that does not belong to a facial triangle, namely, $K_6$ and $L_6$. But $G'$ has at most 10 vertices, because it is obtained from $J$ by splitting one vertex, and hence $G'$ is isomorphic to $K_6$. We have $u_1,z_1\in V(G')$, but $u_1$ is not adjacent to $v_2$ (in $G_0$ and hence in $G'$) and $z_1$ is not adjacent to $v_1$, because there are no exceptional vertices and no parallel edges. Thus $v_1, v_2 \not\in V(G')$. It follows that $G'$ can be obtained from $J$ by first deleting a vertex of degree at least six (and some other vertices) and then adding an edge. This is impossible because $J=C_3+C_5$ or $J=K_2+H_7$. 

Thus $z_2 \not\in \{u_1,u_2\}$. Suppose that one of $v_0$ or $v_0'$ is adjacent to both $z_1$ and $z_3$. Now $N(v_0)$ is triangle-free. Thus all the vertices of degree five in $J$ must have degree six in $G_0$. As these are subset of $\{u_1,u_2,z_1,z_2,z_3\}$, we find that $J=C_3+C_5$. Moreover $u_1,u_2,z_1,z_2,z_3$ are distinct and induce $C_5$. As in cases (iv) and (v), we may color so that the boundary of $\hat{R}$ only uses colors 1, 2, and 3 and then color $v_0$ and $v_0'$ with colors 4 and 5, a contradiction.

So we may suppose that neither $v_0$ or $v_0'$ is adjacent to both $z_1$ and $z_3$. Thus $v_0'$ must be adjacent to both $v_1$ and $v_2$. Without loss of generality, we may assume that one of $v_0$ or $v_0'$ is adjacent to both $z_1$ and $z_2$ thourgh $\hat{R}_2$. Now $N(v_0)$ has at most one triangle and hence the good vertices of $J$ must have degree six in $G_0$. However, $z_3$ is not adjacent to either $v_0$ or $v_0'$. Thus the good vertices must be a subset of $\{u_1,u_2,z_1,z_2\}$. Hence $J=K_2+H_7$. Color $G_0$ as follows. Consider a $5$-coloring of $G_0\setminus \{v_0, v_0'\}$. If $\{u_1,u_2\}$ and $\{z_1,z_2\}$ do not receive the same pair of colors, then we may extend this coloring to $G_0$ by Lemma~\ref{7cycle-5coloring}. So we may assume they are colored with colors 1 and 2. But then no other vertex of $J-z_0$ must receive colors 1 or 2. By swapping the colors of $u_1$ and $u_2$ if necessary, we may assume that $u_2$ and $z_1$ have the same color. We may now recolor $v_1$ with this color and extend the coloring to $G_0$ by Lemma~\ref{7cycle-5coloring}. 

We may now assume that $m=2$. 
It follows that $v_0$ is adjacent to $u_1$ and $u_2$ and that $v_0'$ is adjacent to $z_1,z_2,v_1,v_2$ and $v_0$. Now $N(v_0)$ has at most one triangle. Moreover $N(v_0)$ is triangle-free unless the edge $v_1u_2$ or the edge $v_2u_1$ is in the interior of the unique facial $4$-cycle in $J$. So let us suppose that $N(v_0)$ has a triangle. Without loss of generality suppose the edge $v_1u_2$ is present. All of the good vertices of $J$ must be degree six in $G_0$. These must be a subset of $\{u_1,u_2,z_1,z_2\}$. Hence, $J=K_2+H_7$ and these vertices induce a perfect matching. Repeating the argument from the above paragraph, color $G_0$ as follows. Consider a $5$-coloring of $G_0\setminus \{v_0, v_0'\}$. If $\{u_1,u_2\}$ and $\{z_1,z_2\}$ do not receive the same pair of colors, then we may extend this coloring to $G_0$ by Lemma~\ref{7cycle-5coloring}. So we may assume they are colored with colors 1 and 2. But then no other vertex of $J-z_0$ must receive colors 1 or 2. By swapping the colors of $u_1,u_2$ if necessary, we may assume that $u_2$ and $z_1$ have the same color. We may now recolor $v_1$ with this color and extend the coloring to $G_0$ by Lemma~\ref{7cycle-5coloring}. 

So $N(v_0)$ is triangle-free. All vertices of degree five in $J$ must be degree six in $G_0$. Such vertices are a subset of $\{u_1,u_2,z_1,z_2,x_1,x_2\}$ where $x_1, x_2$ are the ends of an edge in the interior of the facial $4$-cycle of $J$. 
Let us assume that $J=K_2+H_7$.
There are six vertices of degree five; hence, all of these vertices are 
distinct. As $x_1$ is not adjacent to $x_2$ in $J$, we may assume without 
loss of generality that $x_1$ is a good vertex, while $x_2$ may be good or not.
Color $G_0$ as follows. Consider a $5$-coloring of $G_0\setminus \{v_0, v_0'\}$. If $\{u_1,u_2\}$ and $\{z_1,z_2\}$ do not receive the same pair of colors, then we may extend this coloring to $G_0$ by Lemma~\ref{7cycle-5coloring}. So we may assume they are colored with colors 1 and 2. 

We claim that one of the pairs $\{u_1,u_2\}$ and $\{z_1,z_2\}$ only sees two 
other colors in $J\backslash z_0$. 
Suppose not. If $x_2$ is good, then it follows that the other two vertices of degree five receive the same color but they are adjacent, a contradiction. If $x_2$ is not good, then as the pair which contains two good vertices sees all the colors 3, 4, and 5 then $x_1$ will also see 3, 4, and 5, as well as 1, 2 from the pair which contains one good vertex and one not good vertex. Hence $x_1$ cannot receive a color, a contradiction. Now consider the pair, say $\{u_1, u_2\}$ that only sees 2 other colors, say colors 3 and 4. As $v_1$ and $v_2$ are not colored the same, one of $v_1$ or $v_2$ must not be colored 5. Without loss of generality suppose $v_1$ is not colored 5. Then recolor $u_1$ with color 5 and extend this coloring to $G_0$ by Lemma~\ref{7cycle-5coloring}.

So we may assume that $J=C_3+C_5$. Suppose that at least one of $\{u_1,u_2,z_1,z_2\}$ is not a vertex of degree five in $J$. Then it must be exactly one, say $u_1$. Consider a $5$-coloring of $G_0\setminus \{v_0,v_0'\}$. If $u_1$ is not colored the same as one of $\{z_1,z_2\}$, then this coloring extends to $G_0$ by Lemma~\ref{7cycle-5coloring}. However, as $u_1$ is not a vertex of degree five, it is adjacent to all of $J-z_0$ and hence to $z_1$ and $z_2$, it cannot be colored the same as $z_1$ or $z_2$. Thus we may assume that all of $\{u_1,u_2,z_1,z_2\}$ are vertices of degree five in $J$. 

Consider a $5$-coloring of $G_0\setminus \{v_0,v_0'\}$. Now $\{u_1,u_2\}$ must receive the same colors as $\{z_1,z_2\}$, as otherwise this coloring extends to $G_0$ by Lemma~\ref{7cycle-5coloring}. Now the other vertex of degree five in $J$, call this $x_1$ must receive a third color, say color 3. Meanwhile, the other two vertices of $J-z_0$ must receive new colors, namely, colors 4 and 5. Now if $v_1$ is not adjacent to any vertex of color 1, we may recolor $v_1$ by 1 and extend the coloring to $G_0$.  Similarly with color 2 and the same applies for $v_2$ with colors 1 and 2. So we may assume that $u_1$ and $z_1$ are colored 1 and $u_2$ and $z_2$ are colored 2. Further, we may assume that $u_1$ and $z_2$ are adjacent to $x_1$. As $x_1$ must have degree six in $G_0$ there exists an edge $x_1x_2$ through the the $4$-cycle in $J$. As it is not a parallel edge, $x_2\in \{v_1, v_2, u_2, z_1\}$. Thus at least one of $u_2,z_1$ is not adjacent to $x_1$. We may assume without loss of generality that $u_2$ is not adjacent to $x_1$. Now recolor $u_2$ by 3. If the resulting coloring of $G_0 \setminus v_0,v_0'$ is proper, then we may extend it to $G_0$ by Lemma~\ref{7cycle-5coloring}, a contradiction. Thus $u_2$ must be adjacent to a vertex colored 3. As $u_2$ is not adjacent to $x_1$ nor to $v_1$, that vertex must be $v_2$. So recolor $v_2$ by color 2. The resulting coloring is proper as $v_2$ is not adjacent to $z_2$. This coloring extends to a coloring of $G_0$ by Lemma~\ref{7cycle-5coloring}, a contradiction.~\qed

\begin{lemma}
\mylabel{R1R2touch}
Let $(G_0,v_0)$ be an optimal pair,
let $v_1,v_2$ be an identifiable pair, let $J$ be a subgraph
of $G_{v_1v_2}$ isomorphic to $L_1$, $L_2$, $L_5$ or $L_6$, and
let ${R_1}$, ${R_2}$ be the hinges of $J$. 
If $R_1$ and $R_2$ share a vertex  $u\ne z_0$ and at least one of
them has length three, then the other one has length five and
there exists an index $i\in\{1,2\}$ such that
$\hat{R_1}\cup \hat{R_2}\backslash \{v_0,v_i\}$ is a cycle in $G_0$
that bounds an open disk containing $v_0$ and $v_{i}$.
%In particular, at least one of $R_1$, $R_2$ have length five.
\end{lemma}

\proof By the symmetry we may assume that $R_2$ has length three.
Thus $u$ is adjacent to $z_0$ in $J$.
Since $R_1$ is an induced cycle, the cycles $R_1,R_2$ share the
edge $z_0u$.
Thus $\hat{R_1}$,  $\hat{R_2}$ share the edge $v_iu$ for some $i\in\{1,2\}$,
and the second conclusion follows.
By Lemma~\ref{7cycle-5coloring} the cycle
$\hat{R_1}\cup \hat{R_2}\backslash \{v_0,v_i\}$
has length at least six, and hence $R_1$ has length five, as desired.~\qed

We denote by $K_5^-$ the graph obtained from $K_5$ by deleting an edge,
and by $K_5-P_3$ the graph obtained from $K_5$ by deleting two adjacent edges.

\begin{lemma}
\mylabel{finds}
Let $(G_0,v_0)$ be an optimal pair, 
let $v_1,v_2$ be an identifiable pair, and let $J$ be a subgraph
of $G_{v_1v_2}$ isomorphic to $L_1$, $L_2$, $L_5$ or $L_6$.
Then there exists a vertex $s\in V(G_0)-\{v_0\}$ of degree five
such that 
\myitem{\rm(i)} $N_{G_0}(s)$ has a subgraph isomorphic to $K_5-P_3$, and
\myitem{\rm(ii)} if both  hinges of $J$ have length five, then
$N_{G_0}(s)$ has a subgraph isomorphic to $K_5^-$.
\end{lemma}

\proof
We only prove the first assertion, leaving the second one to the reader.
A proof of the second assertion may be found in~\cite{YerPhD}.
Assume that the notation is as in the paragraph prior to 
Lemma~\ref{G'notc3c5k2h7}, and
suppose first that $J = L_5$.   
Let the vertices of $J$ be numbered as in Figure~\ref{fig:l5l6}.
%\ignore{By Lemma ~\ref{facial-cycle-lemmas-11} and the symmetry
%between $\{a_1, a_2, a_3, a_4\}, \{b_1, b_2, b_3\}$, and $\{b_4,
%b_5\}$ we may assume that the faces around $y$ are $a_1ya_2,
%a_2ya_3, a_3ya_4, a_4yb_5b_3x, b_5yb_4, b_4ya_1xb_1$, in order. we
%know that the split can only occur at $x$.}  
It follows from Lemma~\ref{facial-cycle-lemmas-11}
that the indices
of $a_i$ and $b_j$ can be renumbered so that
the faces of $J$ around $c_1$ are $a_1c_1a_2, a_2c_1a_3, a_3c_1a_4,
 a_4c_1b_3b_5c_2,
b_3c_1b_2, b_2c_1b_1, b_1c_1a_1c_2b_4$, in order. Recall that $z_0$ is the
vertex of $J$ that results from the identification of $v_1$ and $v_2$.   
If $z_0 \neq c_1$, then one of
the vertices $a_2, a_3, b_2$ is not incident with $\hat{R_1}$ or
$\hat{R_2}$, and hence has the same neighbors in $J$ and in $G_0$.
It follows that such a vertex satisfies the conclusion of
the lemma, as desired.  We will use the same argument again later,
whereby we will simply say that a certain vertex satisfies the
conclusion of the lemma.

Thus we may assume that $z_0 = c_1$, 
and since we may assume that no vertex satisfies the
conclusion of the lemma, we deduce that one of $R_1$ and $R_2$ is
the face $a_2c_1a_3$ and the other is $b_1c_1b_2$ or $b_2c_1b_3$.  Thus we
may assume that $R_1$ is $a_2c_1a_3$ and $R_2$ is $b_1c_1b_2$.  We may
assume, by swapping $v_1$ and $v_2$, that the neighbors of $v_1$ in
$\hat{J}$ are $a_1, a_2, v_0, b_1$ and that the neighbors of $v_2$
are $a_3, a_4, b_3, b_2, v_0$.  Hence the face $\hat{R}$ is
$v_1a_2a_3v_2b_2b_1$.  
%Recall that $v_0$ is adjacent to both $v_1$ and $v_2$. 
Now $v_1$ is not adjacent to $a_3$ in $G_0$, for
otherwise $a_2$ satisfies the conclusion of the lemma. We shall
abbreviate this argument by $a_2 \Rightarrow v_1 \not \sim a_3$.
Similarly, we have 
%\ignore{$b_1 \Rightarrow b_2 \not \sim v_1$,}
$b_5 \Rightarrow b_3 \not \sim c_2$ and $b_3 \Rightarrow v_2 \not \sim
b_5$.  We shall define 
a $5$-coloring $c$ of $\hat{J} \setminus v_0$.
Let $c(a_1) = c(v_2)= c(b_5) = 1, c(a_2) = c(b_1) = 2,
c(a_3) = c(v_1) = 3, c(a_4) = 4,$ and $c(c_2) = c(b_3) = 5$.  
Assume first that $b_4$ is adjacent to $a_1$. Then $b_2$ is not adjacent to
$v_1$, for otherwise $b_1$ satisfies the conclusion of the lemma.
Furthermore, there is no vertex of $G$ in the face of $\hat J$
bounded by $b_1v_1a_1c_2b_4$.
In that case  we let $c(b_4)= 4$ and  $c(b_2) = 3$.
If  $b_4$ is not adjacent to $a_1$, then we let
$c(b_4)= 3$ and $c(b_2) = 4$.
In either case it follows from Lemma~\ref{7cycle-5coloring}
and the fact that $v_0$ is adjacent to $v_1$ and $v_2$
 that $c$ extends to a 5-coloring of $G_0$,
a contradiction. 
This completes the case $J = L_5$.

%If $b_4$
%is adjacent to $a_1$, and $v_1$ is adjacent to $b_2$, then  $c(b_4)
%= 3, c(b_2) = 4$ and $c(b_5) = 1$.  If $b_4$ is not adjacent to
%$a_1$, and regardless of whether $v_1$ is adjacent to $b_2$, let
%$c(b_4) = 1, c(b_2) = 4,$ and $c(b_5) = 3.$ (This is a modified
%argument, please check.)
% \ignore{Now if\right] \right\rbrace \right( \right( 
%$b_4 $ is adjacent to $a_1$, we define $c(b_4) = 4, c(b_2) = 3$ and
%$c(b_5) = 1$, and otherwise we define $c(b_4) = 1, c(b_2) = 4$ and
%$c(b_5) = 3$. Further, if $v_1$ is adjacent to $b_2$, we now have
%that $c(b_4) = 3, c(b_2) = 4$ and $c(b_5) = 1$.} 

If $J = L_6$ we proceed analogously.  
By Lemma~\ref{facial-cycle-lemmas-11} we may assume that the
%
%First if $z_0 \neq x$,
%then one of the vertices satisfies the conclusion of the lemma. The
faces around $c_1$ are $a_1c_1a_2$, $a_2c_1a_3$, $a_3c_1a_4$, $a_4c_1b_4b_5c_2$,
$b_4c_1b_3$ $b_3c_1b_2$, $b_2c_1b_1$ and $b_1c_1a_1c_2b_5$. 
If $z_0\ne c_1$, or if  one of $R_1$, $R_2$ is not $a_1c_1a_2$ or
$b_2c_1b_3$, then one of $a_2,a_3,b_2,b_3$ satisfies the conclusion of
the lemma.
%Note that the only
%place where we can split vertex $x$ is such that the one of $R_1$
%and $R_2$ is $a_1xa_2$ and the other is $b_2xb_3$. If any other
%faces of length three were $R_1$ or $R_2$, then one of $a_2, a_3,
%b_2$ or $b_3$ would satisfy the conclusion of the lemma.
Thus we may assume that $R_1$ is $a_2c_1a_3$ and
$R_2$ is $b_2c_1b_3$.  We may also assume, by swapping $v_1$ and $v_2$
that the neighbors of $v_1$ in $\hat{J}$ are $a_1, a_2, v_0, b_1$
and $b_2$ and the neighbors of $v_2$ in $\hat{J}$ are $a_3, a_4,
b_3, b_4, $ and $v_0$.  Now $a_1 \Rightarrow v_1 \not \sim y$,  $b_4
\Rightarrow v_2 \not \sim b_5 $,  $a_3 \Rightarrow a_2 \not \sim
v_2$, and $b_2 \Rightarrow b_3 \not \sim v_1 $.  With these
constraints in mind and recalling that $v_0$ is adjacent to $v_1$
and $v_2$, consider the following coloring: $c(a_4) = c(b_1) = 1,
c(a_1) = c(b_2) = 2, c(b_3) = c(v_1) = c(c_2) = 3, c(a_3) = c(b_4) =
4$ and $c(b_5) = c(a_2) = c(v_2) = 5$. 
%Also $c(v_3) = 2$ and $c(v_0) = 4$. (remove this sentence?) 
It follows from Lemma
~\ref{7cycle-5coloring} 
%(7-cycle, remove this comment when combining this document) 
that $c$ extends to a 5-coloring of $G_0$, a
contradiction. This completes the case $J = L_6$.

We now consider the case $J=L_1$.
By Lemma~\ref{facial-cycle-lemmas-10} exactly one face of $J$, say $F$,
is bounded by a cycle of length five, and the remaining faces are bounded 
by triangles. Furthermore, we may assume, by swapping $b_1,b_2$,
and by permuting $a_1,a_2,a_3$ that the faces around $c_1$ in order are
$F, b_2c_1b_1,b_1c_1c_0,c_0c_1a_1,a_3c_1a_1,a_2c_1a_3$.
By swapping $b_3,b_4$ we may assume that the faces around $c_2$ are
$F,b_3c_2b_4,b_4c_2c_0,c_0c_2a_\alpha,a_\beta c_2a_\alpha, a_\gamma c_2a_\beta$
for some distinct indices $\alpha,\beta,\gamma\in\{1,2,3\}$.
Thus the face $F$ is bounded by the cycle $c_1a_2c_2b_3b_2$, and hence
$\gamma=2$. Since $a_1c_0c_1, c_1c_0b_1,b_4c_0c_2$ and $c_2c_0a_\alpha$
are faces of $J$ we deduce that the faces around $c_0$ in order are
$a_1c_0c_1,c_1c_0b_1,b_1c_0b_i,b_ic_0b_j,b_jc_0b_4,b_4c_0c_2,c_2c_0a_\alpha,
a_\alpha c_0a_\delta,a_\delta c_0a_1$ for some integers $i,j,\delta$
with $\{i,j\}=\{2,3\}$ and $\delta\in\{2,3\}-\{\alpha\}$.
Since $\gamma=2$ we have $\alpha\ne2$, and hence $\alpha=3$ and $\delta=2$.

Now if $z_0\ne c_0$, then one of the vertices $a_1,a_2,a_3,b_1,b_2,b_3,b_4$
satisfies the conclusion of the lemma, and hence we may assume that  $z_0= c_0$.
Furthermore, it is not hard to see that one of the above vertices satisfies
the conclusion of the lemma unless one of $R_1$, $R_2$ is
$a_1c_0a_2$ or $a_2c_0a_3$ and the other is one of
$b_1c_0b_i$, $b_ic_0b_j$, $b_jc_0b_4$.
Thus by symmetry we may assume that $R_1$ is $a_1c_0a_2$
and that $R_2$ is one of
$b_1c_0b_i$, $b_ic_0b_j$, $b_jc_0b_4$.

We may assume that in $\hat J$ the vertex $v_1$ is adjacent to $c_1$
and $v_2$ is adjacent to $c_2$.
We see that $a_3\Rightarrow c_1 \not \sim c_2$ and 
$a_3 \Rightarrow a_1 \not \sim v_2$.
Furthermore, if $R_2$ is the face $b_1c_0b_i$, then
$b_4 \Rightarrow b_1 \not \sim v_2$,
and if $R_2$ is the face $b_1c_0b_i$, then
$b_1 \Rightarrow v_1 \not \sim b_4$.
Let $c$ be the coloring of $\hat{J}\backslash v_0$ defined by
$c(b_1) = c(v_2) = 1$, $c(b_i) = c(a_1) = 2$,
$c(b_j) = c(v_1) = c(a_3) = 3$, $c(b_4) = c(a_2) = 4$, and $c(x) =
c(y) = 5$,
and let $c'$ be obtained from $c$ by changing the colors of the
vertices $v_1,v_2,a_2$ to $4,2,1$, respectively.
It follows from Lemma~\ref{7cycle-5coloring} 
by examining the three cases for $R_2$ separately that one of $c,c'$
extends to a $5$-coloring of $G$, a contradiction.
This completes the case $G=L_1$.

Finally, let $J=L_2$. We proceed similarly as above, using
Lemma~\ref{facial-cycle-lemmas-10}.
Let $F$ be the unique face of $J$ of size five.
By renumbering $a_1,a_2,a_3$ and $b_1,b_2,b_3$ we may assume
that the faces around $c_1$ are 
$F,b_3c_1b_2,b_2c_1b_1,b_1c_1c_0,c_0c_1a_1,a_1c_1a_3,a_3c_1a_2$.
Then the faces around $c_2$ are
$F,b_4c_2c_0,c_0c_2a_\alpha,a_\alpha c_2a_\beta,a_\beta c_2a_\gamma$
for some distinct integers $\alpha,\beta,\gamma\in\{1,2,3\}$.
It follows that $\gamma=2$ and that $F$ is bounded by $c_1b_3b_4c_2a_2$.
Since $b_1c_1c_0,c_0c_1a_1,b_4c_2c_0,c_0c_2a_\alpha$ are faces of
$J$ we deduce that $\alpha\ne1$ (and hence $\alpha=3$ and $\beta=1$) and that
the cyclic order of the neighbors of $c_2$ around $c_2$ is
$c_1b_1b_ib_jb_4c_2a_3a_2a_1$ for some distinct integers $i,j\in\{2,3\}$. 
(Recall that all faces incident with $c_0$ are triangles.)
Since $b_4$ is adjacent to $b_3$ in the boundary of $F$ we deduce that
$i=3$ and $j=2$.

Similarly as above, it is easy to see that some $a_i$ or $b_j$ satisfies
the conclusion of the lemma, unless $z_0\in\{c_0,c_1\}$.
Suppose first that $z_0=c_1$.
We may assume that $R_1$ is $b_1b_2c_1$ and $R_2$ is $a_1a_3c_1$,
for otherwise some vertex satisfies the conclusion of the lemma.
We may assume that $v_1$ is adjacent to $a_2,a_3,b_2,b_3$.
We have $a_2\Rightarrow v_1\not\sim c_2$,
$a_1\Rightarrow a_3\not\sim v_2$ and
$b_2\Rightarrow v_1\not\sim b_1$.
Let $c(a_2) = c(b_2) = 1, c(a_3) = c(b_4) = c(v_2) = 2, c(a_1) = c(b_3)
= 3, c(v_1) = c(b_1) = c(c_2) = 4, $ and $c(c_0) = 5$.  It follows from
Lemma ~\ref{7cycle-5coloring} 
that $c$ extends to a 5-coloring of $G_0$,
a contradiction.
Thus we may assume that $z_0=c_0$.
Similarly as above we may assume that
$R_1$ is $b_1b_3c_0$ or $b_3b_2c_0$ and that $R_2$ is $a_1a_2c_0$
or $a_2a_3c_0$.
We may assume that $v_1$ is adjacent to $a_1$ and $b_1$.
If  $R_2$ is $a_1a_2c_0$, then we have
$a_3\Rightarrow c_1\not\sim c_2$ and
$a_3\Rightarrow a_1\not\sim v_2$.
If  $R_2$ is  $a_2a_3c_0$, then 
$a_1\Rightarrow c_1\not\sim c_2$ and
$a_1\Rightarrow a_3\not\sim v_2$.
If $R_1$ is $b_1b_3c_0$, then 
$b_2\Rightarrow b_1\not\sim v_2$.
Let 
$c(a_1)=c(b_1)=c(v_2)=1, c(b_3)=2, c(a_2)=c(b_2)=3, c(a_3)=c(b_4)=c(v_1)=4$
and $c(c_1)=c(c_2)=5$.
It follows from
Lemma ~\ref{7cycle-5coloring}
that $c$ extends to a 5-coloring of $G_0$,
a contradiction.~\qed

% junk

%junk

\begin{lemma}
\mylabel{L1L2L5L6}
Let $(G_0,v_0)$ be an optimal pair,
let $v_1,v_2$ be an identifiable pair, and let $J$ be a subgraph
of $G_{v_1v_2}$. Then $J$ is not isomorphic to $L_1$, $L_2$, $L_5$ or $L_6$.
\end{lemma}

\proof Let $G,v_0,v_1,v_2$ and $J$ be as stated, and suppose for
a contradiction that $J$ is isomorphic to $L_1$, $L_2$, $L_5$ or $L_6$.
Let $R_1,R_2$ be the hinges of $J$, and let $\hat J$, $\hat{R_1}$
and $\hat{R_1}$ be as prior to Lemma~\ref{G'notc3c5k2h7}.
 From Lemma~\ref{finds} and  conditions (ii)--(iv) in the definition of
an optimal pair we deduce that

\claim{1} {$N_{G_0}(v_0)$ has a subgraph isomorphic to $K_5-P_3$,}

and

\claim{2} 
  {if both $R_1$ and $R_2$ have length five, then $v_1,v_2$ is the
   only non-adjacent pair of vertices in $N_{G_0}(v_0)$.}

Let $v_3,v_4,v_5$ be the remaining neighbors of $v_0$ in $G_0$.
If at least two of them belong to the interior of $\hat{R_1}$ or 
$\hat{R_2}$, then
they belong to the interior of the same face, say $R_1$, by (1).
But then $\hat{R_1}$ is bounded by a cycle of length seven, 
and that again contradicts (1) by inspecting the outcomes of
Lemma~\ref{7cycle-5coloring}.
Thus at most one of $v_3,v_4,v_5$ belongs to the interior of $\hat{R_1}$
or $\hat{R_2}$.

From the symmetry we may assume that the edges $v_0v_4$ and $v_0v_5$
belong to the face $\hat{R_1}$. 
We may also assume that $v_5$ belongs to the boundary of $\hat{R_1}$,
and that if $v_4$ does not belong to the boundary of $\hat{R_1}$,
then the edge $v_0v_3$ belongs to $\hat{R_2}$.
We claim that $v_4$ belongs to the boundary of $\hat{R_1}$.
To prove this suppose to the contradiction that $v_4$ belongs to the
interior of $\hat{R_1}$. Then one of the edges $v_1v_4$, $v_2v_4$ does not
belong to $G_0$, and so we may assume $v_2v_4$ does not.
By (1) $v_1,v_2$ and $v_2,v_4$ are the only non-adjacent pairs of
vertices in $N_{G_0}(v_0)$, and by (2) at least one of $R_1,R_2$ has
length three. It follows that $v_3$ belongs to the boundary of $\hat{R_1}$,
and the choice of $v_4,v_5$ implies that the edge $v_0v_3$ lies in the face
$\hat{R_2}$. Thus $v_3$ belongs to the boundary of $\hat{R_2}$.
By Lemma~\ref{R1R2touch} there exists an index $i\in\{1,2\}$ such that
the cycle $R_1\cup R_2\backslash \{v_0,v_i\}$
bounds a disk containing $v_0,v_i$ in its interior. By shortcutting
this cycle through $v_0$ we obtain a cycle of $G_0$ of length at most four
bounding a disk that contains the vertex $v_i$ in its interior,
contrary to Lemma~\ref{7cycle-5coloring}.
This proves our claim that $v_4$ belongs to the boundary of $\hat{R_1}$.
We may assume that $v_0,v_1,v_4,v_5,v_2$ occur on the boundary of
$\hat{R_1}$ in the order listed.

%If $v_3$ belongs to the interior of $\hat{R_1}$, then $v_4, v_5$ belong
%to the boundary of $R_1$ by (1); and if none of $v_3,v_4,v_5$ belong to the
%interior of $\hat{R_1}$ or $\hat{R_2}$, then some two belong to the
%same face. Thus we may assume that $v_4, v_5$ belong to the boundary of 
%$\hat{R_1}$.

%We claim that $v_4,v_5$ do not belong to the boundary of $\hat{R_2}$.
%For suppose that one of them does. 
%Then by Lemma~\ref{R1R2touch} both $R_1$ and $R_2$ have length five,
%and hence $v_1v_5,v_2v_4\in E(G_0)$ by (2). Both those edges must lie
%in $\hat{R_2}$, because they cannot be in $\hat{R_1}$ ($v_0$ is adjacent
%to $v_4$ and $v_5$) or any other face (all other faces are triangles).
%Thus $R_1$ and $R_2$ share the vertices $z_0,v_4,v_5$, and so $J=L_6$.
%It follows by inspecting the two faces of $J_6$ of length five
% that the edges $v_1v_5,v_2v_4$ cannot be
%simultaneously added to $\hat{R_2}$, because they would have to cross.
%That contradiction proves our claim that 
%$v_4,v_5$ do not belong to the boundary of $\hat{R_2}$.

Let $e\in E(G_0)$ have ends either $v_1,v_5$, or $v_2,v_4$.
Then $e\not\in E(\hat J)$, because the boundary of $\hat{R_1}$ is
an induced cycle of $\hat J$.
Moreover, $e$ does not belong to the face $\hat{R_1}$, because the edges
$v_0v_4,v_0v_5$ belong to that face.
Thus $e$ belongs to $\hat{R_2}$ or a face of $\hat J$ of length five.
We claim that $e$ does not belong to $\hat{R_2}$.
To prove the claim suppose to the contrary that it does, and from
the symmetry we may assume that $e=v_2v_4$.
We now argue that not both $R_1,R_2$ are pentagons.
Indeed, otherwise $v_1$ is adjacent to $v_5$ by (2), and the edge
$v_1v_5$ belongs to $\hat{R_2}$, because there is no other face of length
at least five to contain it.
In particular, $v_4,v_5$ belong to the boundary of $\hat{R_2}$,
and because the edges $v_1v_5$, $v_2v_4$ do not cross inside $\hat{R_2}$,
the vertices $v_1,v_0,v_2,v_4,v_5$ occur on the boundary of $\hat{R_2}$
in the order listed.
It now follows by inspecting the $5$-cycles of $L_5$ and $L_6$
that this is impossible.
Thus not both $R_1,R_2$ are pentagons.
By Lemma~\ref{R1R2touch} the cycle
%there is an index $i\in\{1,2\}$ such that
$\hat{R_1}\cup \hat{R_2}\backslash \{v_0,v_1\}$  bounds
a disk with $v_0,v_1$ in its interior. 
By shortcutting this cycle using the chord $v_2v_4$ we obtain
a cycle in $G_0$ of length at most five bounding a disk with at least
two vertices in its interior, contrary to Lemma~\ref{7cycle-5coloring}.
This proves our claim that $v_1v_5$ and $v_2v_4$ do not lie in the
face $\hat{R_2}$.

By (1) and the symmetry we may assume that $v_2v_4\in E(G_0)$,
and hence the edge 
%Since the boundary of $\hat{R_1}$ is an induced cycle in $\hat J$
%we deduce that $v_2v_4\not\in E(\hat J)$, and hence the edge
$v_2v_4$ belongs to a face $\hat F$ of $\hat J$ such that 
$\hat F\ne \hat{R_1},\hat{R_2}$.
Let $F$ be the corresponding face of $J$.
Since $F$ is bounded by an induced cycle, we deduce that $v_4$ is
not adjacent to $z_0$ in $J$. 
Consequently, $R_1$ has length five.
Thus $R_1$ and $F$ have length five, and all other faces of $J$,
including $R_2$, are triangles.
In particular, $J=L_5$ or $J=L_6$,
and $v_1,v_5$ are not adjacent in $G_0$ (because no face of $\hat J$
can contain the edge $v_1v_5$).
By (1) $v_1,v_2$ and $v_1,v_5$ are the only non-adjacent pairs of
vertices in $N_{G_0}(v_0)$.
Condition (1) also implies that $v_3$ belongs to the boundary of $\hat{R_2}$.
Using that and the fact that $v_3$ is adjacent to $v_1$ and $v_2$ in $G_0$,
it now follows that there exists a vertex of $G_0\backslash v_0$ 
whose neighborhood 
in $G_0$ has a subgraph isomorphic to $K_5-P_3$.
%, contrary to the fact that $(G_0,v_0)$ is an optimal pair.
Finding such a vertex requires a case analysis reminiscent of but simpler
than  the proof of Lemma~\ref{finds}. We omit further details.
The existence of such a vertex contradicts the fact that 
$(G_0,v_0)$ is an optimal pair.~\qed

\begin{lemma}
\mylabel{L3L4}
Let $(G_0, v_0)$ be an optimal pair, let $v_1, v_2$ be an
identifiable pair, and let $J$ be a subgraph of $G_{v_1v_2}$.  Then
$J$ is not isomorphic to $L_3$ or $L_4$.
\end{lemma}

\proof Let $G_0, v_0, v_1, v_2$ and $J$ be as stated, and suppose
for a contradiction that $J$ is isomorphic to $L_3$ or $L_4$.  Let
$R_1, R_2$ be the hinges of $J$, and let $\hat{J}, \hat{R_1},
\hat{R_2}$ be as prior to Lemma~\ref{G'notc3c5k2h7}.  
Since by Euler's formula
$J$ triangulates the Klein bottle, we deduce that the faces $
\hat{R_1}, \hat{R_2}$ have size five, and every other face of
$\hat{J}$ is a triangle.  Let the boundaries of $\hat{R_1}$ and
$\hat{R_2}$ be $v_1v_0v_2a_1b_1$ and $v_1v_0v_2cb_l$, respectively.
Let the neighbors of $v_1$ in $\hat J$ in cyclic order be $v_0,b_1, b_2, \ldots,
b_l$, and let the neighbors of $v_2$ in $\hat J$ be 
$v_0,a_1, a_2, \ldots, a_k, c$.
Then $\deg_J(z_0) = k + l + 1$.  Since $J$ has no parallel edges the
vertices $a_1, a_2, \ldots, a_k, c, b_l, b_{l-1}, \ldots, b_1$ are
distinct, and since $J$ is a triangulation they form a cycle, say
$C$, in the order listed.  Since $v_1$ is not adjacent to $v_2$ in $G_0$,
Lemma~\ref{7cycle-5coloring} implies that $|V(C)| \geq 7$.

Let us assume that $|V(C)| = 7$.  Then $z_0$ has degree seven, and
hence $J = L_4$, because $L_3$ has no vertices of degree seven.  Let
$X$ be the set of neighbors of $z_0$ in $J$.  By inspecting the
graph obtained from $L_4$ by deleting a vertex of degree seven, we
find that for every $x \in X$, there exists a $5$-coloring of $J
\setminus z_0$ such that no vertex of $X - \{x\}$ has the same color
as $x$. But this contradicts Lemma~\ref{7cycle-5coloring} 
applied to the subgraph of
$G_0$ consisting of all vertices and edges drawn in the closed disk
bounded by $C$, because $X = V(C)$.  This completes the case when
$|V(C)| = 7$.

Since $L_3$ and $L_4$ have no vertices of degree eight, it follows
that $|V(C)| = 9$, and hence $z_0$ is the unique vertex of $J$ of
degree nine.  From the symmetry between $v_1$ and $v_2$, we may
assume that $\deg_{\hat{J}}(v_1) \leq 5$; in other words $l \leq 4$.
The graph $J$ is 6-critical.  Since $z_0$ is adjacent to every other
vertex of $J$, we deduce that $J \setminus z_0 \setminus x$ is
$4$-colorable for every vertex $x \in V(J) - \{z_0\}$,
and hence
%.  We deduce the following:
\begin{enumerate}
\item[(1)]  {\it for every vertex  $x \in V(J)-\{z_0\}$, 
the graph $J \setminus z_0$ has a 5-coloring 
such that $x$ is the only vertex colored $5$.}
\end{enumerate}

From Lemma~\ref{7cycle-5coloring} 
applied to the boundary of the face $\hat{R}$ of
$\hat{J} \setminus v_0$, we deduce that one of $\hat{R_1},
\hat{R_2}$ contains no vertex of $G_0$ in its interior, and the
other contains at most one.  Since $v_0$ has degree five, we may
assume from the symmetry between $\hat{R_1}$ and $\hat{R_2}$ that
$v_0$ is adjacent to $a_1$ and $b_1$ (and hence $\hat{R_1}$ includes
no vertices of $G_0$ in its interior).  We claim that $l = 4$ and
$v_1$ is adjacent to $c$.  To prove the claim suppose to the
contrary that either $l \leq 3$ or $v_1$ is not adjacent to $c$.
Then $\deg_{\hat{J}}(v_1) \leq 5$.  By (1) there exists a coloring of
$J \setminus z_0 = \hat{J} \setminus \{v_0, v_1, v_2 \}$ such that
$b_1$ is the only vertex colored 5.  We give $v_2$ the color $5$,
then we color $v_1$, then we color the unique vertex in the interior
of $\hat{R_2}$ if there is one, and finally color $v_0$.  The last
three steps are possible, because each vertex being colored sees at
most four distinct colors.  Thus we obtain a $5$-coloring of $G_0$, a
contradiction.  This proves our claim that $l = 4$ and $v_1$ is
adjacent to $c$.  It follows that $k = 4$ and $V(G_0) = \{v_0, v_1,
v_2, a_1, a_2, a_3, a_4, b_1, b_2, b_3, b_4, c \}$.  We have
$\deg_{G_0}(v_1) =  \deg_{G_0}(v_2) = 6$, and since $\deg_J(c) \leq
\deg_{G_0}(c) - 2$, we deduce that $\deg_{G_0}(c) \geq 7$.  Thus we
have shown that

\begin{enumerate}
\item[(2)]   { \it if $x_1, x_2, x_3, x_4, x_5$ are the
neighbors of $v_0$ in $G_0$ listed in their cyclic order around $v_0$, 
the vertex $x_1$
is not adjacent to $x_3$ in $G_0$ and $G_{x_1, x_3}$ has a subgraph
isomorphic to $L_3$ or $L_4$, then $\deg_{G_0}(x_1) = \deg_{G_0}(x_3)
= 6 $ and $\deg_{G_0}(x_2) \geq 7$.}
\end{enumerate}

It also follows that $v_1$ is not adjacent to $a_1$ in $G_0$ and
that $v_2$ is not adjacent to $b_1$ in $G_0$.  Not both $G_{v_1a_1}$
and $G_{v_2b_1}$ have a subgraph isomorphic to $L_3$ or $L_4$ by
$(2)$, and so from the symmetry we may assume that $G_{v_1a_1}$ does
not.  By the optimality of $(G_0, v_0)$ and Lemmas~\ref{G'notc3c5k2h7}
 and~\ref{L1L2L5L6}, it
follows that $G_{v_1a_1}$ has a subgraph isomorphic to $K_6$.
Thus $G\backslash\{v_0,v_1,v_2\}$ has a subgraph $K$ isomorphic to $K_5$. If
$v_2 \not \in V(K)$, then $V(K) \cup \{z_0\}$ induces a $K_6$
subgraph in $J$, a contradiction.  Thus $v_2 \in V(K)$, and hence
$V(K) = \{v_2, a_2, a_3, a_4, c\}$.  Let $i \in \{3,4\}$.  If $a_1$
is not adjacent to $a_i$ in $G_0$, then we 5-color $G_0$ as follows.
By (1) there is a 5-coloring of $G_0 \setminus \{v_0, v_1, v_2\}$
such that $a_1$ and $a_i$ are the only two vertices colored 5.  We
give $v_1$ color $5$, then color $v_2$ and finally $v_0$.  Similarly
as before, this gives a valid 5-coloring of $G_0$ a contradiction.
Thus, $a_1$ is adjacent to $a_3$ and $a_4$ and hence $a_1$ is not
adjacent to $c$, for otherwise $\{a_1, a_2, a_3, a_4, v_2, c\}$
induces a $K_6$ subgraph in $G_0$.

Since $\deg_{G_0} (v_2) = 6$, it follows from (2) that $G_{ca_1}$ has
no subgraph isomorphic to $L_3$ or $L_4$.  By the optimality of
$(G_0, v_0)$ and Lemmas~\ref{G'notc3c5k2h7} and~\ref{L1L2L5L6}
 it follows that $G_{ca_1}$ has
a subgraph isomorphic to $K_6$.  By an analogous argument as above
we deduce that $\{v_1, b_1, b_2, b_3, b_4\}$ is the vertex-set of a
$K_5$ subgraph of $G_0$.  The existence of the two $K_5$ subgraphs
implies that $a_2, a_3, a_4, b_2, b_3, b_4$ have $K_4$ subgraphs in
their neighborhoods, and the optimality of $(G_0, v_0)$ implies that
$a_2, a_3, a_4, b_1, b_2, b_3$ all have degree at least six in
$G_0$, and hence in $J$.  Thus $a_1, b_1, c$ are the only vertices
of $J$ of degree five.  Thus, $J = L_3$ and $a_1, b_1, c$ are
pairwise adjacent, a contradiction, because we have shown earlier
that $a_1$ is not adjacent to $c$.~\qed

The results of this section may be summarized as follows.

\begin{lemma}
\mylabel{J=K6}
Let $(G_0,v_0)$ be an optimal pair, and let $v_1,v_2$ be an identifiable
pair.  Then $G_{v_1v_2}$ has a subgraph isomorphic to $K_6$.
\end{lemma}

\proof Every $5$-coloring of $G_{v_1v_2}$ may be extended to a $5$-coloring
of $G_0$, and hence $G_{v_1v_2}$ is not $5$-colorable.
By the choice of $G_0$ the graph $G_{v_1v_2}$ has a subgraph isomorphic
to one of the graphs listed in Theorem~\ref{main}.
By Lemmas~\ref{G'notc3c5k2h7}, \ref{L1L2L5L6} and~\ref{L3L4}
that subgraph is $K_6$, as desired.~\qed

\section{Using $K_6$}
\mylabel{usek6}
%%%%%%%%%%%%%%%%%%%%%%%%%%%%%%%%%%%%%%%%%%%%%%%%%%%%%%%%%%%%%%%%%%%%%%%%%%%

\begin{lemma} 
\mylabel{THOM-9}
Let $(G_0,v_0)$ be an optimal pair.  Then 
$G_0$ has at least $10$ vertices, and if it has exactly $10$, then
it has a vertex of degree nine.
\end{lemma}

\proof This follows immediately from Lemma~\ref{smallcrit}.~\qed

% junk

\begin{lemma} 
\mylabel{K5MINUS}
Let $(G_0,v_0)$ be an optimal pair.
Then there are at least two identifiable pairs.
\end{lemma}

\proof
Since $G_0$ has no subgraph isomorphic to $K_6$ there is at least
one identifiable pair.
Suppose for a contradiction that $v_1,v_2$ is the only identifiable
pair. Thus the subgraph of $G_0$ induced by $v_0$ and its neighbors
is isomorphic to $K_6$ with one edge deleted.
By Lemma~\ref{J=K6} the graph $G_0\backslash\{v_0,v_1,v_2\}$ has
a subgraph $K$ isomorphic to $K_5$, and every vertex of $K$ is
adjacent to $v_1$ or $v_2$.
Let $t$ be the number of neighbors of $v_0$ in $V(K)$.
Since $v_0$ has degree five and its neighbors $v_1,v_2$ are not in $K$
it follows that $t\le3$.
If $t=0$, then $G_0$ has a subgraph isomorphic to $L_5$ or $L_6$;
if $t=1$, then  $G_0$ has a subgraph isomorphic to $L_1$ or $L_2$;
if $t=2$, then  $G_0$ has a subgraph isomorphic to $K_2+H_7$; and 
if $t=3$, then  $G_0$ has a subgraph isomorphic to $C_3+C_5$.~\qed

%We say that a vertex $v$ in a graph $G$ embedded in a surface
%has a {\em wheel neighborhood} if the neighbors of $v$ form a cycle in
%the order determined by the embedding.

\begin{lemma} 
\mylabel{wheelnbhd}
Let $(G_0,v_0)$ be an optimal pair.
Then $v_0$ has a wheel neighborhood.
\end{lemma}

\proof
Let us say that a vertex $v\in V(G_0)$ is a {\em fan} if its neighbors
form a cycle in the order determined by the embedding of $G_0$.
We remark that if $v_0$ is a fan and $v_0$ does not have a wheel neighborhood, then the embedding of $G_0$ can be
modified to $G_0'$ so that $v_0$ will have a wheel neighborhood contradicting condition (vi).
Thus it suffices to show that $v_0$ is a fan.
Suppose for a contradiction that  there
exist non-adjacent vertices $a,b \in N(v_0)$ 
that are consecutive in the cyclic order of the
neighbors of $v_0$.  By  condition
(iv) in the definition of an optimal pair, 
the graph $G' = G_0+ab$ has a subgraph
$M$ isomorphic to one of the graphs from  Theorem~\ref{main}.  
Assume, for a contradiction,
that $v_0 \notin V(M)$.  By optimality condition (i), $G_0
\backslash v_0$ has a 5-coloring $c$. 
Since $c$ is not a $5$-coloring of $M$ it follows that $c(a) = c(b)$. 
But then $c$ can be extended to a 5-coloring of $G_0$, a contradiction.  
Thus $v_0 \in V(M)$.  Since $\deg(v_0) = 5$, we get that
$N_{G_0}(v_0)\subseteq V(M)$.  Further note that $a,b$ are adjacent in $M$,
because $M$ is not a subgraph of $G_0$.  
%We proceed be considering the identity of $M$.

First, assume $M = K_6$.  Then $V(M) = \{v_0\} \cup N(v_0)$. This
implies that there is at most one identifiable pair, 
contrary to Lemma~\ref{K5MINUS}.
Second,  assume $M = L_3$  or $L_4$. As each is a triangulation,
Lemma $\ref{7cycle-5coloring}$ implies that $G_0 = M\backslash ab$.  But $M$ is
$6$-critical, so $G_0$ has a $5$-coloring, a contradiction.

Third, assume that $M = C_3 + C_5$ or $K_2 + H_7$. Because $M$ 
is one edge short of being a triangulation,
there is a unique face in $M$ of length four.  As $ab
\in E(M)$, the embedding of $M\backslash ab$ has at most two faces of size strictly
bigger than three, and if it has two, then they both have size four.
Since $G_0$ has at least $10$ vertices by Lemma~\ref{THOM-9},
Lemma $\ref{7cycle-5coloring}$  implies
that $M\backslash ab$ has a face $f$ of size five whose interior includes
a vertex of degree five.
However, $f$ is the only face of  $M\backslash ab$ of size at least four,
and hence it also includes the edge $ab$, but that is impossible.
This completes the case when $M = C_3 + C_5$ or $K_2 + H_7$.

Fourth, suppose that $M$ is either $L_5$ or $L_6$,
and let the notation be as in the proof of Lemma~\ref{finds}.
In particular, every face incident with $a_2$ or $b_2$ is a triangle.
At least one of $a_2,b_2$, say $s$, is not equal to $v_0$ and does not
include both $a,b$ in its neighborhood.
But then the neighborhoods of $s$ in $G$ and in $M$ are the same,
and hence $s$ satisfies conditions (ii)-(iv) in the definition of an optimal pair
by Lemma~\ref{K5MINUS}. 
But $s$ is a fan in $M$, and hence has a wheel neighborhood in some embedding
of $G_0$, contrary to condition (vi) in the definition of optimal pair.

If $M=L_1$, then we apply the argument of the previous paragraph
to the vertices $a_1,b_1,b_4$,
using the notation of Lemma~\ref{finds}.
Finally, suppose that $M$ is $L_2$,
and let the notation be again as in the proof of Lemma~\ref{finds}.
Every face incident with one of the vertices $a_3,b_2$
is a triangle, and at least one of those vertices, say $s$,
has the property that $s\ne v_0$ and if 
the neighborhood of $s$ includes both $a$ and $b$,
then  $a,b$ are not consecutive in the cyclic ordering
around $s$ and
$\{a,b\}\cap\{x,y\}\ne\emptyset$ for every pair of distinct non-adjacent
vertices $x,y\in N_M(v_0)$.
Since $s$ is a fan in $M$ its choice implies
that it is  a fan in $G_0$, and hence  has a wheel neighborhood in some embedding
of $G_0$. Furthermore, in $G_0$ there are at most two pairs of non-adjacent 
vertices in the neighborhood of $s$, and if there are two, then they
are not disjoint.
Thus $s$ satisfies conditions (ii)-(iv) in the definition of an optimal pair
by Lemma~\ref{K5MINUS}, contrary to 
condition (vi) in the definition of an optimal pair.~\qed

A drawing of a graph $G$ in a surface is {\em $2$-cell} if every
face of $G$ is homeomorphic to an open disk.

\begin{lemma} 
\mylabel{crosscapinface}
Let $(G_0,v_0)$ be an optimal pair, and let $v_1,v_2$ be an identifiable
pair, and let $J$ be a subgraph of $G_{v_1v_2}$ isomorphic to $K_6$.
Then the drawing of $J$ in the Klein bottle is $2$-cell.
% there does not exist an
%embedding of $J \cong K_6$ such that a face surrounds a crosscap.
%(Should we just say that the embedding is 2-cell?  Yes, I think.)
\end{lemma}

\proof
Let $v_0,R_1,R_2,\hat R_1,\hat R_2$ be as before,
and suppose for a contradiction that the drawing of $J$ is not $2$-cell.
Since $K_6$ has a unique drawing in the projective 
plane~\cite[page 364]{KocKre},
\nocite{KocKre}
it follows that every face of $J$ is bounded by a triangle, and
exactly one face, say $F$, is homeomorphic
to the M\"obius strip. If $F$ is not $R_1$ or $R_2$, then the boundary
of $F$ is a separating triangle of $G_0$, a contradiction, because
no $6$-critical graph has a separating triangle. Thus we may assume
that $F=R_2$.

Since both $R_1$ and $R_2$ are triangles, and they share at least
one vertex, there exists a vertex $s\in V(J)$ not incident with
$R_1$ or $R_2$. Thus in $\hat J$ all the faces incident with $s$
are triangles, and hence $\deg_{G_0}(s)=\deg_{J}(s)=5$ by  
Lemma~\ref{7cycle-5coloring}.
Furthermore, if $R_1$ and $R_2$ share an edge, then 
$N_{G_0}(s)$
has a subgraph isomorphic to $K_5^-$, the complete graph on
five vertices with one edge deleted.
This implies, by the
optimality of $(G_0,v_0)$, that 
$N_{G_0}(v_0)$ has a subgraph isomorphic to $K_5^-$,
contradicting Lemma~\ref{K5MINUS}.

So we may assume that $R_1$ and $R_2$ have no common edge. Let the
facial walk incident with $\hat{R_1}$ be $v_0,v_1,z_1,z_2,v_2,v_0$,
and the facial walk incident with $\hat{R_2}$ be
$v_0,v_1,z_3,z_4,v_2,v_0$.  Notice, from the embedding of $J$, that the
$z_i$ are distinct.  Also notice that $s$ is the lone vertex in
$\hat J$ not incident with either $\hat{R_1}$ or $\hat{R_2}$, and
$N_{G_0}(s)$ includes no two disjoint pairs of non-adjacent vertices.
This implies, by the optimality of
$(G_0,v_0)$, that 
$N_{G_0}(v_0)$ includes no two disjoint pairs of non-adjacent vertices.
We shall refer to this as the DP property.

Let $N(v_0)=\{v_1,v_2,v_3,v_4,v_5\}$.  Assume that some neighbor of
$v_0$, say $v_3$, belongs to $\hat{R_1}$.  By Lemma
$\ref{7cycle-5coloring}$, $v_3$ is adjacent to all vertices incident
with $\hat{R_1}$.  Thus $v_4$ and $v_5$ belong to the closure of $R_2$.
In either case, $v_3$ and $v_4$ are not adjacent in $G_0$.  Since $v_1$
and $v_2$ are also not adjacent, this contradicts the DP property.
% to the previous paragraph.

Since $v_1$ is not adjacent to $v_2$ in $G_0$ it follows from 
Lemma~\ref{wheelnbhd} that at least one of $v_3,v_4,v_5$ belongs 
to the closure of $\hat R_1$. Thus there remain two cases,
depending on whether one or two of those vertices are incident with $\hat R_1$.
If it is two vertices, then we may assume without loss of
generality that $v_3=z_1$ and $v_4=z_2$.  As $z_1$ and $z_2$ are
not incident to $\hat R_2$, $v_3 ,v_2$ and $v_4 ,v_1$ are not
adjacent in $G_0$, contrary to the DP property.
Thus we may assume that $v_3=z_1$ and $v_4$ and $v_5$ belong to the
closure of $\hat R_2$.
By the DP property $v_3, v_4$ and
$v_3,v_5$ are adjacent in $G_0$.  Thus, without loss of generality, $v_4=z_3$ and
$v_5=z_4$.  
Furthermore, 
%since $N_{G_0}(v_0)$ includes no two disjoint pairs of  non-adjacent vertices, 
it follows from the DP property that either $v_1,v_5$ or  $v_2,v_4$
are adjacent in $G_0$.
%, and so we may assume the former.
Thus the subgraph $L$ of $\hat J$ consisting of the vertices
$v_0,v_1,v_2,v_4,v_5$ and edges between them that belong to the
closure of $\hat R_2$ has five vertices and at least eight edges.
We can regard $L$ as drawn in the M\"obius band with the cycle
$v_1v_0v_2v_5v_4$ forming the boundary of the  M\"obius band.
As such the graph $L$ has at least three faces.
% One of them is bounded by the cycle  $v_1v_0,v_2,v_5,v_4$, 
%and of the other three faces, at most one has length 
%at least five. 
Since the sum of the lengths of the faces is at least $11$, at most one of
them has length at least five.
That face of $L$ includes at most one vertex of $G_0$
by Lemma~\ref{7cycle-5coloring}, and the other faces of $L$ include none.
Thus $G_0$ has at most nine vertices, contrary to Lemma~\ref{THOM-9}.~\qed

\begin{lemma}
\mylabel{6faceexists}
Let $(G_0, v_0)$ be an optimal pair,  let $v_1,v_2$ be an
identifiable pair, and let $J$ be a subgraph of $G_{v_1v_2}$ isomorphic to
$K_6$. Then some face of $J$ has length six.
\end{lemma}

\proof
Let $\tilde{J}$ denote the graph consisting of $\hat{J}$
and edges of $G_0$ not in $\hat{J}$ from $v_1$ or $v_2$ to
the boundary of $\hat R_1$ or $\hat R_2$ that are drawn inside 
$\hat R_1$ or $\hat R_2$.
Let $\tilde{R_1}$ be the face in $\tilde{J}$ that
contains $v_0$ and is contained in $R_1$, and similarly for
$\tilde{R_2}$. 
We assume for a contradiction that no face of $J$ has length six.
By Lemma~\ref{crosscapinface} the embedding of $J$ is $2$-cell,
and so, by Euler's formula, all faces of $J$ are bounded by triangles,
except for either three faces of length four,
or one face of length four and one face of length five.
Each face of $\tilde J$ other than $\tilde{R_1}$ and $\tilde{R_2}$
will be called {\em special} if it has length at least four.
Thus there are at most three special faces, and if there are exactly three,
then they have length exactly four.

  Let us denote the vertices on the boundary
of $\tilde{R_1}$ as $v_1$, $v_0$, $v_2$, $u_1, \ldots, u_k$ in order, 
and let the vertices on the boundary of $\tilde{R_2}$ be
$v_2$, $v_0$, $v_1$, $z_1, \ldots, z_l$ in order.
Note that $u_1,u_2, \ldots,u_k$ are pairwise distinct, and similarly for
$z_1,z_2, \ldots,z_l$.
A special face of length five may include a vertex of $G_0$ in its
interior; such vertex will be called {\em special}.
It follows that there is at most one special vertex.
An edge of $G_0$ is called {\em special} if it has both ends in
$\hat J\backslash v_0$, but does not belong to $\hat J$, and is not
$v_1z_1$ or $v_2z_1$ if $l=1$, and is not $v_1u_1$ or $v_2u_1$ if $k=1$.
It follows that every special edge is incident with $v_1$ or $v_2$.
Furthermore, the multigraph obtained from $G_0$ by deleting all
vertices in the faces $\tilde R_1$ and $\tilde R_2$ and contracting
the edges $v_0v_1$ and $v_0v_2$ has $J$ as a spanning subgraph,
and each special edge belongs to a face of $J$ of length at least four.
It follows that there are at most three special edges.
Furthermore, if there is a special vertex, then there is at most one
special edge, and each increase of $k$ or $l$ above the value of two
decreases the number of special edges by one.

%Let $\tilde{J}$ denote the graph consisting of $\hat{J}$
%and edges of $G_0$ not in $\hat{J}$ from $v_1$ or $v_2$ to 
%the boundary of $R_1$ or $R_2$ that are drawn inside $R_1$ or $R_2$. 
%Let $\tilde{R_1}$ be the face in $\tilde{J}$ that
%contains $v_0$ and is contained in $R_1$. Similarly for
%$\tilde{R_2}$.

Since $R_1$ and $R_2$ have length three, four, or five, we deduce
that $k,l\in\{1,2,3,4\}$.
The graph $\hat J\backslash\{v_0,v_1,v_2\}=J\backslash z_0$ is
isomorphic to $K_5$, and $u_1,u_2, \ldots,u_k$ are its distinct vertices;
let $u_{k+1}, \ldots,u_5$ be the remaining vertices of this graph.
It follows that if $c$ is a $5$-coloring of $\tilde J$ and
$c(u_i)=c(z_j)$, then $u_i=z_j$.
We will refer to this property as {\em injectivity}.
From the symmetry we may assume that $k\ge l$.
Since $J$ has at most one face of length five, it follows that $l\le 3$.
%$|\tilde{R_1}| \ge|\tilde{R_2}|$. 
We distinguish three cases depending
 on the value of $l$.
%is either 4, 5, or 6, because $J$ has at most one face of length five.
%We will assume to a contradiction that there
%exists an identifiable pair $(v_1, v_2)$ such that $G_{v_1 v_2}$
%contains a subgraph $J$ isomorphic to $K_6$, embedded such that the
%embedding of $J$ has either three 4-cycles or one 4-cycle and one
%5-cycle.

%{\bf IMPROVE}
%In any 5-coloring of a subgraph of $G_0$ that contains the induced
%subgraph $J \setminus \{z_0\}$, which is isomorphic to $K_5$, and we
%will refer to the vertex of $J \setminus \{z_0\}$ that gets colored
%$i$ as $k_i$.

\medskip
\noindent
{\bf Case 1: $l=1$}

\noindent
%Note that $|\tilde{R_1}| \le 7$.  
By Lemma~\ref{wheelnbhd} the vertex $v_0$ is adjacent to $z_1$.
Also notice then that 
$v_1z_1v_2u_1u_2 \ldots u_k$ is a null-homotopic walk $W$ of length 
at most seven.
Since $G_0$ is $6$-critical, the graph $G\backslash v_0$ has a $5$-coloring,
say $c$.
By Lemma~$\ref{7cycle-5coloring}$ applied to the subgraph $L$ of $G_0$ drawn 
in the disk bounded by $W$ and the coloring $c$, 
the graph $L$ satisfies one of (i)--(vi) of that lemma.
We discuss those cases separately.

Case (i): There are eight vertices in $\tilde{J}$ and none in the
interior of $\tilde{R_1}$ and $\tilde{R_2}$, and at most one special vertex. 
%If the embedding of $J$
%has a 5-cycle which is not $R_1$ or $R_2$, then there may be a
%vertex inside that 5-cycle.  Nevertheless, 
Thus $|V(G_0)| \le 9$,
contradicting Lemma $\ref{THOM-9}$.

Case (ii): As before $|V(G_0)| \le 9$, a contradiction, unless there
exists a special vertex $v_0'$.
%embedding of $J$ has a 5-cycle which is not $R_1$ or $R_2$ and there
%exists a vertex $v_0'$ inside that 5-cycle. 
This implies that
$|\tilde{R_1}| = |\hat{R_1}| = 6$.
Without loss of generality suppose
$v_0$ is adjacent to $u_3, v_1, z_1, v_2$ and a
 vertex $v_3$ which is adjacent to $v_0,v_2, u_1, u_2, u_3$.
Notice that $v_0'$ must have degree five in $G_0$ and its neighborhood
must contain a subgraph isomorphic to $K_5-P_3$, since four of its
neighbors are in $J \setminus z_0$ and thus form a clique.
Meanwhile the neighborhood of $v_0$ is missing the edges
$v_1v_2$, $v_1v_3$, and $v_2u_3$. 
The last one does not belong to $\tilde J$, does not lie in $\tilde{R_1}$
(because we have already described the graph therein), and is not special,
because all special edges have been accounted for.
%a different face of  $\tilde J$, because all such faces are bounded by triangles.
Thus the pair $(G_0,v_0')$ contradicts the optimality of $(G_0,v_0)$.

Case (iii): 
The graph $L\backslash W$ consists of three pairwise adjacent vertices,
and $v_0$ is one of them. Let $v_3,v_4$ be the remaining two.
%Let us recall that $u_1,u_2, \ldots,u_k$ are pairwise adjacent in $G_0$,
%and hence the coloring $c$ assigns them pairwise different colors.
%Yet there exist three distinct colors such that each of $v_0,v_3,v_4$
%sees a vertex of $W$ of each of those three colors.
%This can happen in only two ways, up to symmetry that exchanges
By Lemma~\ref{wheelnbhd} we may assume, using the symmetry that exchanges
$v_1,v_4,u_1,u_2$ with $v_2,v_3,u_k,u_{k-1}$, that
$v_3$ has neighbors $v_0,v_2,u_1,u_2,v_4$
and $v_4$ is adjacent to $v_1,v_0,v_3,u_2$ and either $u_3$ or $u_4$.
%\myitem{(a)} $k\in\{3,4\}$, $v_3$ has neighbors $v_0,v_2,u_1,u_2,v_4$
%and $v_4$ is adjacent to $v_1,v_0,v_3,u_2,u_3$,
%\myitem{(b)} $k=4$, $v_3$ has neighbors $v_0,v_2,u_1,u_2,v_4$
%and $v_4$ is adjacent to $v_1,v_0,v_3,u_2,u_4$.
\noindent 
In either case $z_1$ and $u_2$ are colored the same, and hence
they are equal by injectivity.
To be able to treat both cases simultaneously, we swap $u_3$ and $u_4$
if necessary; thus we may assume that $v_4$ is adjacent to $u_3$.
We can do this, because we will no longer use 
the order of $u_1,u_2, \ldots,u_k$ for the duration of case (iii).
The vertex  $v_1$ is adjacent to $u_2,u_3,u_4,u_5$, for otherwise
its color can be changed, in which case the coloring $c$ could be extended
to $L$, contrary to the fact that $G_0$ has no $5$-coloring.
Similarly, $v_2$ is adjacent to $u_1,u_2,u_4,u_5$.
It follows that $G_0$ has a subgraph isomorphic to $L_3$, a contradiction.
%$v_0,v_1,\ldots,v_4,u_1,u_2,\ldots,u_5$ is isomorphic to $L_3$,
To describe the isomorphism, the vertices corresponding to the top row
of vertices in Figure~\ref{fig:allels}(c) in left-to-right order are
$u_1,u_4,u_5,u_3$,  the vertices corresponding to the middle row are
$v_3,v_2,u_2=z_1,v_1,v_4$, and the bottom vertex is $v_0$.
This completes case (iii).

Cases (iv)-(vi): 
We have $k=4$. Hence $R_1$ has length five, and therefore there is
at most one special edge.
%every face of $\tilde J$ is bounded by a triagle,
%except for $\tilde{R_1}$ and one other face that has length four.
%Thus at most one edge of $E(G_0)-E(\tilde J)$ does not belong to
%$\tilde{R_1}$ or $\tilde{R_2}$.
%It follows from the description of cases (iv)-(vi)
%in Lemma~$\ref{7cycle-5coloring}$ that $\tilde{R_1}$ includes no edge
%of $G_0$ with one end  $v_1$ and the other in $\{u_1,u_2,u_3\}$,
%and no edge with one end $v_2$ and the other end in   $\{u_2,u_3,u_4\}$.
%Since in $\hat J$ not both $v_1,v_2$ are adjacent to $u_i$ for $i=1,2,3,4$,
%we deduce that the same is true in $\tilde J$, possibly with one exception.
Consequently, one of $v_1$, $v_2$ is not adjacent in $\tilde J$
to at least two vertices
among $u_1,u_2,u_3,u_4$. 
Since  every face of $\tilde J$ except $\tilde{R_1}$ and one other face of length
four is bounded by a triangle this implies that in the coloring $c$, one of
 $v_1$, $v_2$ sees at most three colors.
From the symmetry we may assume that $v_2$ has this property.
Thus $c(v_2)$ may be changed to a different color.

By using this fact and examining the cases (iv)-(vi) of 
Lemma~$\ref{7cycle-5coloring}$ we deduce
that $L$ is isomorphic to the graph of case (iv).
Let the vertices of $L$ be numbered as in Figure~\ref{fig:criticalgraphs}(iv).
It further follows that $v_2=x_4$ or $v_2=x_5$, and so from the symmetry
we may assume the former. Since $z_1$ has a unique neighbor in $L\backslash W$
we deduce that $z_1=x_3$, $v_1=x_2$, $u_4=x_1$ and so on.
Notice that $x_8$ has degree five in $G_0$ and that its neighborhood
is isomorphic to $K_5-P_3$. Meanwhile, the neighborhood of $v_0$ is
certainly missing the edges $v_1v_2$ and $v_1x_9$. Now if 
$x_3\ne x_5$ then $x_3$ is not adjacent to $x_9$ and
$N(v_0)$ is missing at least three edges, a contradiction to the optimality
of $(G_0,v_0)$, given the existence of $x_8$. 
So $x_3=x_5$, but then the edges $x_3v_2$, $x_5v_2$ are actually the
same edge, because $\tilde J$ does not have parallel edges.
It follows that $v_2$ has degree at most four in $G_0$,  a contradiction.

%belongs to
%$E(G_0)-E(\tilde J)$, and hence edges must
%be parallel edges instead of being the same edge, otherwise the
%7-cycle is really a 5-cycle.  So $z_0$ has three edges to the
%$x_3=x_5$ vertex; that is, we added two edges, using $R_2$ for one
%of them and the 4-cycle for the other (since $R_1$ is the
%five-cycle).   But then $G_0$ must already be a triangulation
%without the edge $v_2x_1$.  So the neighborhood is missing three
%edges, again a contradiction.

\medskip
\noindent{\bf Case 2: $l=2$}

\noindent
By Lemma~\ref{wheelnbhd} either   $v_0$ is adjacent to both $z_1$
and $z_2$, in which case we define $\overline v_0:=v_0$, 
or there exists a vertex $\overline v_0$ in $\tilde{R_2}$
adjacent to $v_0,v_1,v_2,z_1,z_2$.
Let $W$ denote the walk
$v_1\overline v_0v_2 u_1 \ldots u_k$ of length at most seven,
and let $X$ denote the set of vertices of $G_0$ drawn in the open disk
bounded by $W$.
We claim that $X\ne\emptyset$.
This is clear if $\overline v_0\ne v_0$,
and so we may assume that $\overline v_0= v_0$. But then $X=\emptyset$
implies $|V(G_0)|\le9$,
contrary to Lemma~\ref{THOM-9}. Thus $X\ne\emptyset$.
Let $x\in X$ have the fewest number of neighbors on $W$.
Since $G_0$ is $6$-critical, the graph $G_0\backslash x$ 
has a $5$-coloring, say $c$.
By Lemma~$\ref{7cycle-5coloring}$ applied to the subgraph $L$ of $G_0$ drawn
in the disk bounded by $W$ and the coloring $c$,
the graph $L$ and coloring $c$ satisfy one of (i)--(vi) of that lemma.
%Again, we discuss those cases separately.
%Since  $G_0\backslash (X\cup\{v_1,v_2\}$ is isomorphic
%to $J\backslash z_0$, which, in turn, is isomorphic to $K_5$,
%we deduce that for every coloring $d$ of $G_0\backslash X$, if $d(u_i)=d(z_j)$,
%then $u_i=z_j$.
%We will refer to this as {\em injectivity}.

Suppose first that $L$ and $c$ satisfy (i).
Then $|X|=1$ by the choice of $x$.
As before $|V(G_0)| \le 9$, contradicting Lemma~$\ref{THOM-9}$,
unless there is a special vertex.
Hence $k\le3$. If $k=3$,
then $R_1$ has length four, and all special faces have been accounted for.
In particular, $\tilde J=\hat J$.
The fact that the coloring $c$ cannot be extended to $L$ implies that
$\{c(z_1),c(z_2)\}\subseteq\{c(u_1),c(u_2),c(u_3)\}$, and hence
$\{z_1,z_2\}\subseteq\{u_1,u_2,u_3\}$ by injectivity. 
Thus $u_1$ or $u_3$ is equal to
one of $z_1,z_2$.
Since there are no special edges, either $u_1v_2$ and $z_2v_2$, 
or $u_kv_1$ and $z_1v_1$ are the same edge, but then $v_1$ or $v_2$ has
degree at most four, a contradiction.
% and hence $J$ has a parallel edge, a contradiction.
If $k=2$ we reach the same conclusion, using the fact that in that case
there is at most one special edge.
%keeping in mind that 
%each edge of $E(\tilde J)-E(\hat J)$ belongs to an special face, and hence
%there is at most one such edge.
It follows that  $L$ and $c$ do not satisfy (i).

Next we dispose of the case $k\le 3$.
To that end assume that $k\le 3$.
Then $W$ has length at most six.
Thus $L$ and $c$ satisfy either (ii) or (iii) of Lemma~\ref{7cycle-5coloring},
and so $W$ has length exactly six and $k=3$.
In particular, $R_1$ has length four, and so there is either at most
one special vertex, or at most two special edges, and not both.
It follows that either $c(v_1)$ or $c(v_2)$ can be changed without
affecting the colors of the other vertices of $G_0\backslash X$.
That implies that  $L$ and $c$ satisfy (ii).
Let $v_3$ be the unique neighbor of $\overline v_0$ in $X$, and let
$v_4$ be the other vertex of $X$.
From the symmetry we may assume that $v_3$ is adjacent to
$\overline v_0,v_1,v_2,u_1,v_4$, and $v_4$ is adjacent to $v_1,v_3,u_1,u_2,u_3$.
By considering the walk $u_1u_2u_3v_1z_1z_2v_2$ and the subgraph drawn
in the disk it bounds, and by applying Lemma~\ref{7cycle-5coloring}
to this graph and the coloring $c$ we deduce that 
$|\{c(u_1),c(u_2),c(u_3)\}\cap \{c(z_1),c(z_2)\}|=1$.
That implies $|\{u_1,u_2,u_3\}\cap \{z_1,z_2\}|=1$ by injectivity,
and so we may assume that $u_5$ is not equal to $z_1$ or $z_2$.
It follows that the neighborhood of $u_5$ has a subgraph isomorphic to
$K_5-P_3$.
 However, the neighborhood of $\overline v_0$ is
missing $v_1v_2$ and at least one of the edges $v_3z_1$ and $v_3z_2$, 
contrary to the optimality of $(G_0,v_0)$ if 
$v_0=\overline v_0$. 
Similarly, the neighborhood of $v_3$ is missing
$v_1v_2$ and $\overline v_0v_4$, a contradiction if $v_0=v_3$.
This completes the case $k\le3$.

Thus we may assume that $k=4$. 
It follows that $R_1$ has length five, and hence there is at most one
special edge. 
Let $i\in\{1,2\}$. If $v_i$ is adjacent to both $z_1$ and $z_2$, then one
of the edges $v_iz_1$, $v_iz_2$ is special.
It follows that in $G_0$, either $v_1$ is not adjacent to $z_2$,
or $v_2$ is not adjacent to $z_1$.
But $z_2$ is the only neighbor of $v_1$ in $G_0\backslash X$
colored $c(z_2)$, because $G_0\backslash (X\cup\{v_0,v_1,v_2\}$ is isomorphic
to $J\backslash z_0$, which, in turn, is isomorphic to $K_5$.
Thus there is a (proper) $5$-coloring $c_1$ of $G_0\backslash X$
obtained by changing the color of at most one of the vertices $v_1,v_2$
such that either
$c_1(v_1)=c_1(z_2)$ or $c_1(v_2)=c_1(z_1)$.
%and $c_1$ is obtained from $c$ by changing the color of at most one
%vertex, namely either $v_1$ or $v_2$.
Now $c_1(\overline v_0)$ can be changed to another color, thus yielding
a coloring $c_2$ of $G_0\backslash X$.

If $L$ and $c$ satisfy one of the cases (iii)-(vi),
then one of the colorings $c_1,c_2$ extends into $L$, a contradiction.
Thus $L$ and $c$ satisfy (ii) of 
Lemma~\ref{7cycle-5coloring}.
Let $v_3\in X$ be the unique vertex of $X$ adjacent to $\overline v_0$, and let 
$v_4$ be the other vertex in $X$.
If both $v_3$ and $v_4$ have degree five in $G_0$, then one of the
colorings $c_1,c_2$ extends into $L$, a contradiction. Thus one of
$v_3$, $v_4$ has degree five, and the other has degree six.
It follows that $v_3$ is adjacent to $v_1$, $v_2$, and either
$u_1$ or $u_4$, and so from the
symmetry we may assume it is adjacent to $u_1$.
If $c_1(v_1)=c_1(u_1)$, then we can extend one of the colorings $c_1,c_2$ 
into $L$ by first coloring $v_4$ and then $v_3$.
Thus $c_1(v_1)\ne c_1(u_1)$.
If $v_4$ is not adjacent to $u_1$, then we can extend $c_1$ or $c_2$
by giving $v_4$ the color $c_1(u_1)$, and then coloring $v_3$.
Thus $v_4$ is adjacent to $v_1$.
If $v_4$ has degree five, then its neighbors are $u_1,u_2,u_3,u_4,v_3$,
and the neighbors of $v_3$ are $\overline v_0,v_1,v_2,u_1,u_4,v_4$.
Let $d$ a $5$-coloring of $G_0\backslash \overline v_0$.
Since the coloring $d$ cannot be extended to $\overline v_0$, it follows
that the neighbors of $\overline v_0$ receive different colors.
Now similarly as in the contruction of $c_1$ above, we can change either
the color of $v_1$, or the color of $v_2$. The resulting coloring then
extends to $\overline v_0$, a contradiction.
This completes the case when $v_4$ has degree five, and hence $v_4$ has
degree six.
It follows that the neighbors of $v_4$ are $u_1,u_2,u_3,u_4,v_1,v_3$
and the neighbors of $v_3$ are $\overline v_0,v_1,v_2,u_1,v_4$.
Let $d_1$ be a $5$-coloring of
the graph $G_0\backslash\{\overline v_0,v_3\}$.
Since the coloring $d_1$ does not extend into $\overline v_0,v_3$,
we deduce that $\{d_1(z_1),d_1(z_2)\}=\{d_1(v_4),d_1(u_1)\}$.
%Given that $G_0\backslash (X\cup\{v_1,v_2\}$ is isomorphic to $K_5$,
By injectivity this implies that $u_1=z_1$ or $u_1=z_2$.
If $u_1=z_2$, then one of the edges $v_2u_1$, $v_2z_2$ is special,
because they cannot be the same edge, given that $v_2$ has degree at least
five in $G_0$.
Thus all special edges have been accounted for, and so $z_1$ is
not adjacent to $u_1$.
Thus $d_1(v_1)$ can be changed to $d_1(u_1)$, and the new coloring
extends to all of $G_0$, a contradiction.
%But $u_1=z_2$ implies that the edges $v_2u_1$ and $v_2z_2$ are actually
%the same edge, which in turn implies that the degree of $v_2$ in $G_0$
%is at most four, a contradiction.
Thus  $u_1=z_1$. It follows that $G_0$ is isomorphic to $L_3$.
First of all, the vertex $v_1$ is not adjacent to both $u_2$ and $u_3$,
for otherwise the vertices $v_1,v_4,u_1,u_2,u_3,u_4$ form a $K_6$ subgraph
in $G_0$.
If $v_1$ is adjacent to neither  $u_2$ nor $u_3$, then $v_2$ is adjacent
to these vertices, and an isomorphism between $G_0$ and $L_3$ is given
by mapping the vertices in the top row in Figure~\ref{fig:allels}(c),
in left-to-right order, to $u_4,u_2,u_3,u_5$, the middle row to
$v_1,v_4,u_1=z_1,v_2,\overline v_0$ and the bottom vertex to $v_3$.
If $v_1$ is adjacent to exactly one of  $u_2$, $u_3$, then due to
the symmetry in the forthcoming argument we may assume that $v_1$
is adjacent to $u_3$, and hence $v_2$ is adjacent to $u_2$.
Then an isomorphism is given by mapping the top row to
$v_4,u_4,u_3,u_2$, the middle row to $v_3,v_1,u_1=z_1,u_5,v_2$,
and mapping the bottom vertex to $\overline v_0$.
This completes the case $l=2$.

\medskip
\noindent{\bf Case 3: $l=3$}

\noindent
Lemma~\ref{wheelnbhd} implies that $v_0$ has at most one neighbor
among $\{z_1,z_2,z_3,u_1,u_2,\ldots,u_k\}$,
and such neighbor must be $u_1$, $u_k$, $z_1$, or $z_3$.

We claim that either $v_0$ is adjacent to $z_1$ or $z_3$, or
$k=3$ and $v_0$ is adjacent to $u_1$ or $u_3$.
To prove this claim let us assume that $v_0$ has no neighbor among
$\{z_1,z_2,z_3\}$.
Let $C$ be the cycle $v_1z_1z_2z_3v_2v_0$, and let $X$ denote the
set of vertices of $G_0$ drawn in the open disk bounded by $C$.
We have  $X\ne\emptyset$ by Lemma~\ref{wheelnbhd}.
Let $c$ be a coloring of $G\backslash X$, and let $L$ denote
the subgraph of $G_0$ consisting of all vertices and edges drawn in
the closed disk bounded by $C$.
%induced by $V(C)\cup X$.
By Lemma~\ref{7cycle-5coloring} the graph $L$ satisfies one of
the conditions (i), (ii), (iii) of that lemma.
The vertices $z_1$ and $z_3$ are adjacent, because the graph obtained
from $\hat J$ by deleting $v_0,v_1$, $v_2$ and the vertices drawn in the
faces $\tilde R_1$ or $\tilde R_2$ is isomorphic to $K_5$.
We may also assume, by the symmetry between $v_1$ and $v_2$, that
$v_1$ is adjacent to $z_2$.
We claim that we may assume that the neighborhood of $v_0$ is a $5$-cycle.
This is clear if $v_0$ has no neighbor in $\{u_1,u_2,u_3,u_4\}$, and
so we may assume that $v_0$ is adjacent to $u_1$.
Then we may assume that $k=4$, for otherwise the claim we are proving holds.
Thus there is no special edge.
By Lemma~\ref{wheelnbhd} there exists a vertex inside $\tilde R_1$
adjacent to $v_0,v_1,u_1$.
Since there is no special edge the vertex $v_1$ is not adjacent to $u_1$,
and $u_1$ is not adjacent to $z_1$, because $v_2$ has degree at least five
in $G_0$.
It follows that the neighborhood of $v_0$ is indeed a $5$-cycle.
If $|X|\ge2$, then 
%the neighborhood of $v_0$ in $G_0$ is a $5$-cycle, while 
there exists a vertex in $X$
whose neighborhood has a subgraph that is a $5$-cycle plus at least
one additional edge, namely $z_1z_3$ or $v_1z_2$.
That contradicts the optimality of $(G_0,v_0)$.
Thus $|X|=1$. Let $x$ denote the unique element of $X$,
and let us assume first that $k=4$.
Then there are no special edges, and so $v_1$ is not adjacent to $z_3$
and $v_2$ is not adjacent to $z_1$.
Let $C'$ denote the cycle $v_1xv_2u_1u_2u_3u_4$, and let $X'$ be the set
of vertices of $G_0$ drawn in the open disk bounded by $C'$.
Then $G_0\backslash (X'\cup\{x\})$ has a $5$-coloring $c'$ such that $c'(v_1)=c'(z_3)$
and $c'(v_2)=c'(z_1)$. 
Then $c'$ can be extended to $x$ in at least two different ways.
By Lemmas~\ref{7cycle-5coloring} and~\ref{wheelnbhd}
the coloring $c'$ can be extended to
all of $G_0$, unless (up to symmetry betweeen $v_1$ and $v_2$)
$v_0$ is adjacent to $u_1$, there exists a vertex $y$ adjacent to 
$u_1,u_2,u_3,u_4$ and $c'(v_1)=c'(u_5)$.
But $v_1$ is not adjacent to $u_1$ (because $v_2$ is and there are
no special edges), and hence the color of $v_1$ can be changed to
$c'(u_1)$, and the resulting coloring extends to all of $G_0$,
a contradiction.
This completes the case $k=4$.
Thus $k=3$, and so there is at most one special edge.
Let $c''$ be a $5$-coloring of  $G_0\backslash X'$. It follows that
the color of at least one of the vertices $v_1,v_2$ can be changed
to a different color, without affecting the colors of the other vertices
of $G\backslash X'$.
It follows from Lemma~\ref{7cycle-5coloring} that $|X'|\le 2$.
That, in turn, implies that $v_0$ is adjacent to $u_1$ or $u_3$, and hence
proves our claim from the beginning of this paragraph.

Thus we may assume that $v_0$ is adjacent to $z_3$.
By Lemma~\ref{wheelnbhd} there exists a vertex $v_3$ adjacent to 
$v_0,v_1,z_1,z_2,z_3$ and a vertex
%By another application of  Lemma~\ref{wheelnbhd} there exists a vertex
$v_4$ in $\tilde R_1$ that is adjacent to $v_0,v_1,v_2$.
The neighborhood of $v_3$ includes the edge $z_1z_3$, and so
by the optimality of $(G_0,v_0)$ the neighborhood of $v_0$ includes the
edge $v_4z_3$. Thus $z_3\in\{u_1,u_2,u_3,u_4\}$.
Assume first that $k=4$. Then there are no special edges, and hence
$z_3\ne u_4$. Next we deduce that $z_3\ne u_1$, for otherwise $v_2u_1$
and $v_2z_3$ are the same edge, which implies (given that $z_3=u_1$
is adjacent to $v_4$) that $v_2$ has degree at most three, a contradiction.
Thus $z_3\in\{u_2,u_3\}$.
Let $Y$ consist of $v_0$ and all vertices 
in $\tilde R_1$ or $\tilde R_2$.
Since $z_3$ is adjacent to $v_4$ we deduce that $|Y|\le 4$.
Since there are no special edges, $z_3$ is not adjacent to $v_1$,
and $v_2$ is not adjacent to $u_4$.
Thus $G_0\backslash Y$ has a coloring $d$ such that $d(v_1)=d(z_3)$
and $d(v_2)=d(u_4)$.
Since $z_3\in\{u_2,u_3\}$ this coloring can be extended to the vertices drawn in
 $\tilde R_1$, and since $d(v_1)=d(z_3)$ it can be further extended to
$v_0$ and $v_3$, a contradiction. 

Thus $k=3$. 
Let $W$ denote the walk $v_1v_3z_3v_2u_1u_2u_3$, and
let $d'$ be a $5$-coloring of $G_0\backslash (Y-\{v_3\})$.
We now apply Lemma~\ref{7cycle-5coloring} to the graph drawn in the
closed disk bounded by $W$ and coloring $d'$, 
and note that either the color of each of
$v_1$, $v_2$ can be changed to a different color, independently of each
other and independently of the colors of other vertices, except possibly
$v_3$, or the color of one of $v_1$, $v_2$ can be changed to two
different values.
In either case, one of the resulting colorings extends to $G_0$, a 
contradiction.~\qed

\begin{lemma}
\mylabel{6walk}
Let $(G_0,v_0)$ be an optimal pair,  let $v_1,v_2$ be an identifiable
pair, and let $J$ be a subgraph of  $G_{v_1v_2}$ 
isomorphic to $K_6$.
Then the drawing of $J$ in the Klein bottle 
does not have a facial walk of length six.
\end{lemma}

\begin{figure}
\centering
\includegraphics[scale=.5]{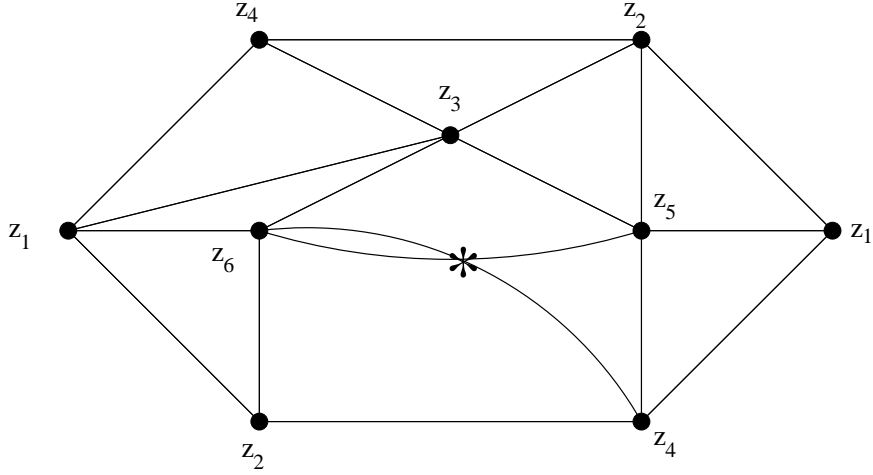}
\caption{An embedding of $K_6$ with a facial walk on five vertices}
\label{fig:5walk}
\end{figure}

\proof
Suppose for a contradiction that there exists a subgraph $J$ of $G_{v_1v_2}$ 
isomorphic to $K_6$ such that the drawing of $J$ in the Klein bottle has a face
$F_0$ bounded by a walk $W$ of length six. 
Let the vertices of $J$ be $z_1,z_2,\ldots,z_6$.
Since $K_7$ cannot be embedded in the Klein bottle, it follows that
$W$ has a repeated vertex.
If $W$ has exactly one repeated vertex, then (since $J$ is simple)
we may assume that the vertices on $W$ are $z_6,z_2,z_4,z_6,z_3,z_5$, 
in order.
There exists a closed curve $\phi$ passing through $z_6$ and otherwise
confined to $F_0$ such that there is an edge of $J$ on either side of $\phi$
in a neighborhood of $z_6$.
The curve $\phi$ cannot be separating, because $G_0\backslash z_0$ is
connected, and it cannot be $2$-sided, because $G_0\backslash z_0$ is not planar.
It follows that $\phi$ is $1$-sided.
By Euler's formula every face of $J$ other than $F_0$ is bounded by
a triangle.
It follows that the triangles $z_4z_5z_6$, $z_1z_6z_3$, and
$z_1z_6z_2$ bound faces of $J$.
Furthermore, either $z_3z_5z_2$ or $z_3z_5z_4$ is a face, but since
$J$ is simple we deduce that it is the former.
It follows that $z_1z_3z_4$, $z_2z_3z_4$, $z_1z_2z_5$ and
$z_1z_4z_5$ are faces of $J$, and those are all the faces of $J$.
The drawing of $J$ is depicted in Figure~\ref{fig:5walk},
where diagonally opposite vertices and edges are identified, and the
asterisk indicates another cross-cap.

\begin{figure}
\centering
\includegraphics[scale=.5]{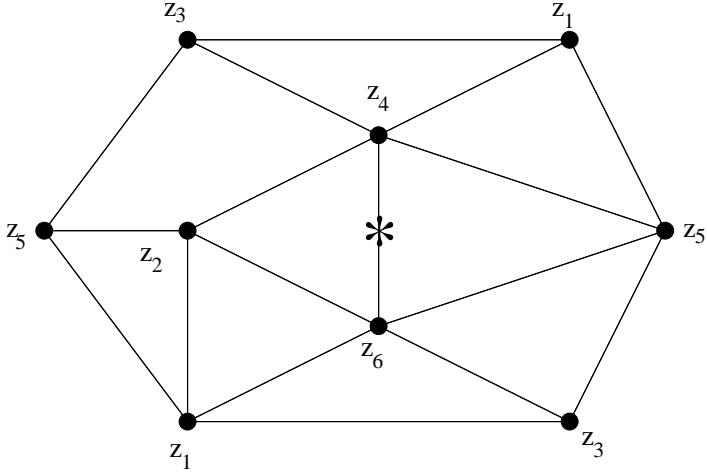}
\caption{An embedding of $K_6$ with a facial walk on four vertices}
\label{fig:4walk}
\end{figure}

Similarly, if $W$ has at least two repeated vertices, then it has
exactly two, and we may assume that the vertices of $W$ are
$z_6z_5z_4z_6z_2z_4$.
Similarly as in the previous paragraph, the embedding is now
uniquely determined, and is depicted in Figure~\ref{fig:4walk}.

%It follows that the vertex-set of $W$ consists of either four or five
%vertices, and for most of the proof we will have to consider those cases
%separately.
%We introduce the following notation.
%If $W$ has five vertices, then
%let $z_6\in V(J)$ be the repeated vertex in
%$W$, let $z_1\in V(J)$ be the vertex not on $W$,
%and label the remaining vertices of $J$ such that the cyclic order in $J$ around
%$z_6$ is $z_1$, $z_2$, $z_4$, $z_5$, $z_3$.
%If $W$ has four vertices, then
%let $z_6$ be a repeated vertex in $W$, let
%$z_4$ be the other repeated vertex, and label the remaining vertices of $J$
%such that the cyclic order around $z_6$ in $J$ is $z_4$, $z_5$, $z_3$, $z_1$,
%$z_2$.
In either case let $R_1$ and $R_2$ be the hinges of $J$, and
let $F_{ijk}$ denote the facial triangle
incident with $z_i,z_j,z_k$ if it exists.
We should note that specifying the hinges does not uniquely determine
the graph $\hat J$, because the face $F_0$ has multiple incidences 
with some vertices.
For instance, if $W$ has five vertices, $z_0=z_6$, and $R_1=F_0$,
then it is not clear whether the split occurs in the ``angle" between
the edges $z_3z_6$ and $z_4z_6$, or in the angle between $z_5z_6$ and
$z_2z_6$. To overcome this ambiguity we will write $R_1=F_{364}$ in
the former case, and  $R_1=F_{265}$ in the latter case.
Notice that this is just a notational device; there is no face bounded
by $z_3z_6z_5$ or  $z_2z_6z_4$.
We proceed in a series of claims.

\begin{enumerate}
\item[(1)]   { \it 
Not both $R_1$ and $R_2$ are bounded by triangles.
}
\end{enumerate}

\noindent
To prove (1) suppose for a contradiction that  $R_1$ and $R_2$ are both facial
triangles. Let us recall that $z_0$ is the vertex of $G_{v_1v_2}$
that results from the identification of $v_1$ and $v_2$.
Suppose first that $R_1$ and $R_2$ are consecutive in the cyclic order
around $z_0$. Then $v_0$ and one of $v_1$ or $v_2$ is in the interior
of a $4$-cycle in $G_0$, contrary to Lemma~\ref{7cycle-5coloring}. 
Similarly, if the cyclic
order around $z_0$ has $R_1$ followed by a facial triangle, followed
by $R_2$, then there would be two vertices in the interior of a $5$-cycle
in $G_0$, contrary to Lemma~\ref{7cycle-5coloring}.
In addition, if the cyclic order has $R_1$,
followed by two facial triangles, followed by $R_2$, then there are
two vertices inside a $6$-cycle. Hence, we are in either case (ii)
or (iii) of Lemma~\ref{7cycle-5coloring}. 
However, the boundary has five vertices
that form a clique. So $5$-color all but the interior of this $6$-walk
(using that $G_0$ is $5$-critical);
the boundary must have five colors, contrary to Lemma~\ref{7cycle-5coloring}. 
We conclude that
$R_1$ and $R_2$ must have $F_0$ in between them in the cyclic order
around $z_0$, on both sides.
In particular, $W$ has five vertices.

Thus the only case remaining is that $z_0=z_6$, where $J$ is
embedded with a facial $6$-walk on five vertices. Suppose without
loss of generality that $R_1 = F_{126}$ and $R_2=F_{456}$, and that  $v_2$
is adjacent to $z_1,z_3$ and $z_4$. Then the faces of the subgraph 
induced by $v_1, v_2, z_1, z_2, z_3, z_4, z_5$
are all triangles but perhaps for two six-cycles: 
$v_1, z_2, z_1, v_2, z_4, z_5$ and $v_1, z_5, z_3, v_2, z_4, z_2$. 
Since $v_0$ is adjacent to $v_1$ and $v_2$ it follows from  
Lemma~\ref{7cycle-5coloring}
that the only vertex in $G_0$ in the interior of the first six-cycle is $v_0$.  
Hence there must be at least two vertices in the interior of the other six-cycle, 
else $|V(G_0)| \le 9$, a contradiction. 
Thus we are in either case (ii) or (iii) of Lemma~\ref{7cycle-5coloring}. 
Note that the disk bounded by the second cycle includes no chord.
%this implies there are no chords through the second six-cycle. 
So $v_1$ is not adjacent to $z_3$. Now if $v_1$ is not adjacent to $z_1$, 
we color $G_0$ as follows. Let the color of $z_i$ be $i$. 
Color $v_1$ with color 1. Then color $v_0$ and $v_2$, and extend the coloring 
to the interior of the second six-cycle by Lemma~\ref{7cycle-5coloring}. 
Hence we may assume that $v_1$ is adjacent to $z_1$. 
But then $v_0$ is adjacent to $z_1, z_4, z_5$ while $v_1$ is not adjacent to $z_4$.
Now $v_1$ may be colored either 3 or 4. One of these options extends to the 
interior of the second six-cycle after we color $v_1$, $v_0$, $v_2$ in that order.
This proves claim (1).
\medskip

%{\bf THIS PARAGRAPH NEEDS TO BE FIXED}
%Thus the only case remaining is that $z_0$ is $z_6$ where $J$ is
%embedded with a facial $6$-walk on five vertices. Suppose without
%loss of generality that $R_1 = F_{126}$ and $R_2=F_{456}$ where $v_2$
%is adjacent to $z_1$,$z_3$ and $z_4$. (to be finished later, should
%be easy)...
%This proves claim (1).
%\medskip

In light of (1) we may assume that $R_1=F_0$.

\begin{enumerate}
\item[(2)]   { \it
If $R_2$ is bounded by a triangle, then it is not consecutive with
$F_0$ in the cyclic order around $z_0$ in $J$.
}
\end{enumerate}

\noindent To prove (2) suppose for a contradiction that
$R_2$ is bounded by a triangle and that it is consecutive with 
$F_0$ in the cyclic order around $z_0$ in $J$.
%Now we may assume that $z_0\in \partial (F_0)$ and without loss of
%generality that $R_1 = F_0$. Suppose $R_2$ is a facial triangle consecutive
% with $F_0$ in the cyclic order around $z_0$. 
It follows that one of $v_1,v_2$ has degree two in $\hat{J}$,
and so we may assume that it is $v_1$ and that its neighbors are
$v_0$ and  $z_j$. 
Consider the subgraph $\hat{J}\setminus \{v_0,v_1\}$.
All of its faces are triangles but for a $7$-walk. We $5$-color this
subgraph, which is isomorphic to $K_6$ minus an edge. 
Thus $v_2$ must receive the same color as $z_j$. 
Since this subgraph only has six vertices, the
interior of the $7$-walk must be as in case (v) or (vi) of 
Lemma~\ref{7cycle-5coloring}, for
otherwise there would be at most nine vertices in $G_0$, contrary to
Lemma~\ref{THOM-9}.
Consider the edge $z_0z_j$ in $J$, which must be on
the boundary of $F_0$. Now the vertex or vertices not on the boundary
of $F_0$ must be on the boundary of $R_2$, for otherwise the $7$-walk
would only have four colors and we could extend the $5$-coloring
to its interior, a contradiction. Since $R_2$ is a facial triangle
this means that either $z_0$ or $z_j$ is $z_6$ and that
$W$ has five vertices.
However, then the color of $z_0$
and $z_j$ appears three times on the boundary of the $7$-walk. So
the $5$-coloring may also be extended, a contradiction.
This proves (2).
\medskip

By an {\em s-vertex} we mean a vertex $s\in V(G_0)$ of degree five such that
$N_{G_0}(s)$ has a subgraph isomorphic to $K_5-P_3$.
If $G_0$ has an s-vertex, then the optimality of $(G_0,v_0)$ implies
that $N_{G_0}(v_0)$ does not include two disjoint pairs of
non-adjacent vertices.

\begin{enumerate}
\item[(3)]   { \it
Let $R_2$ be bounded by a triangle; then $\hat R_2$ is bounded by 
a pentagon, say $v_0v_1r_1r_2v_2$.
Assume further that $G_0$ has an s-vertex. Then either
\myitemitem{\rm(a)} $\hat R_2$ includes a unique vertex $v$ of $G$,
and $v$ is adjacent to $v_0,r_1,r_2$ and all neighbors of $v_0$ other than $v$, or
\myitemitem{\rm(b)} $v_0$ is adjacent to $r_1,r_2$, and $r_1,r_2$ are adjacent
to the neighbor of $v_0$ other than $v_1,v_2,r_1,r_2$, or
\myitemitem{\rm(c)} $v_0,v_1,v_2$ are all adjacent to $r_i$ for some $i\in\{1,2\}$,
and $r_i$ is adjacent to the two neighbors of $v_0$ other than $v_1,v_2,r_i$.
}
\end{enumerate}

\noindent
To prove (3) we first notice that $\hat R_2$ includes at most one vertex
of $G_0$ by Lemma~\ref{7cycle-5coloring}.
If it includes exactly one vertex, then (a) holds by the existence of
an s-vertex, and the optimality of $(G_0,v_0)$.
If $\hat R_2$ includes no vertex of $G_0$, then by
Lemma~\ref{wheelnbhd} either $v_0$ is adjacent to both $r_1$ and $r_2$,
or $v_0,v_1,v_2$ are all adjacent to $r_i$ for some $i\in\{1,2\}$.
We deduce from the existence of an s-vertex and the optimality of $(G_0,v_0)$
that either (b) or (c) holds. 
This proves (3).
\medskip

\begin{enumerate}
\item[(4)]   { \it
The walk $W$ has five vertices.
}
\end{enumerate}

\noindent
To prove (4) we suppose for a contradiction that $W$ has four vertices.
Suppose first that $z_0=z_2$.
Then by (2) and the symmetry we may assume that $R_2=F_{125}$.
It follows that $z_3$ is an s-vertex, and so we may apply (3).
But (a) does not hold, because in that case $v_0$ 
has four neighbors in $\hat R_1$ or on its boundary, and not all of them
can be adjacent to the neighbor of $v_0$ in $\hat R_2$.
%is not adjacent to $z_1$.
If (b) holds, then $v_0$ is adjacent to $z_1$ and $z_5$, and $v$ is adjacent
to $z_1$, where $v$ is the neighbor of $v_0$ other than $v_1,v_2,z_1,z_5$.
Now $v\ne z_5$, because otherwise both $\hat R_1$ and $\hat R_2$
include an edge joining $v_0$ and $z_5$, contrary to the fact that $G_0$
is simple. Since $v$ is adjacent to $z_1$ we deduce that $v=z_4$ or $v=z_6$.
In either case Lemma~\ref{wheelnbhd} implies that $v_1$ or $v_2$ has degree
at most four, a contradiction. 

Thus we may assume that (c) holds, and so $v_0,v_1,v_2$ are all
adjacent to $z_1$ or $z_5$.
In the former case we can change notation so that $R_2=F_{126}$,
contrary to (2). Thus $v_0,v_1,v_2$ are all adjacent to  $z_5$.
Let $v_1$ be adjacent to $z_3,z_4,z_5$; then $v_2$ is adjacent to
$z_1,z_5,z_6$. 
Let the vertices $v_2,z_5,v_1,v_4,v_5$ form the wheel neighborhood of $v_0$,
in order. Since an s-vertex exists, the optimality of  $(G_0,v_0)$ implies
that either  $v_1$ is adjacent to $v_5$,
or $v_2$ is adjacent to $v_4$, or both.
We may assume from the symmetry  that $v_1$ is adjacent to $v_5$.
Since $v_5$ is adjacent to $z_5$ by (c), we deduce that $v_5=z_4$ or
$v_5=z_6$, because $v_5\ne z_5$ for the same reason as above.
If $v_5=z_6$, then $v_2z_6$ and $v_2v_5$ are the same edge, and it
follows from Lemma~\ref{wheelnbhd} that $v_2$ has degree at most four.
Thus $v_5=z_4$. It follows that $v_2$ is adjacent to $z_4$, and hence
the neighborhood of $z_1$ has a subgraph isomorphic to $K_5^-$,
contrary to Lemma~\ref{K5MINUS}. 
This completes the case  $z_0=z_2$.

Thus by symmetry we may assume that $z_0=z_4$.
Again by symmetry we may assume that $R_1=F_{246}$ and 
$R_2$ is either $F_{134}$ or $F_{145}$.
Assume first that $R_2=F_{145}$.
Let $v_1$ be adjacent to $z_1,z_2,z_3$.
Then $z_3$ is an s-vertex, and so we may use (3).
If (a) holds, and $v$ is as in (a), then it is not possible for $v$ to
be adjacent to all neighbors of $v_0$ other than $v$, a contradiction.
If (b) holds, then $v_2$ is not adjacent to $z_1$, and hence
$v_1$ is adjacent to $z_5$, by the optimality of $(G_0,v_0)$, because
an s-vertex exists.
Thus the neighborhood of $z_3$ in $G_0$ has a subgraph isomorphic to
$K_5^-$, contrary to Lemma~\ref{K5MINUS}.
Thus (c) holds.
If $v_0,v_1,v_2$ are adjacent to $z_5$, then 
$N_{G_0}(z_3)$ has a subgraph isomorphic to $K_5^-$, contrary to
Lemma~\ref{K5MINUS}.
%we can change our notation so that $R_2=F_{145}$, contrary to (2).
Hence  $v_0,v_1,v_2$ are adjacent to $z_1$.
By (c) the vertex $z_1$ is adjacent to $v_4,v_5$, the two neighbors of
$v_0$ other than $v_1,v_2,z_1$.
It follows that 
 $\{v_4,v_5\} \subseteq \{z_2,z_5,z_6\}$. 
However, if $v_0$ is adjacent to $z_2$, then $v_1$ would be
of degree at most four in $G_0$, a contradiction. Thus $v_0$ is adjacent to
$z_5$ and $z_6$; hence $v_1$ is adjacent to $z_5$ 
by Lemma $\ref{wheelnbhd}$. Now the graph
has eight vertices and perhaps one more inside the $5$-cycle
$v_1z_2z_6v_2z_5$. Hence $G_0$ has at most nine vertices,
contrary to Lemma~\ref{THOM-9}.
This completes the case  $R_2=F_{145}$.

We may therefore assume that  $R_2=F_{246}$.
From the symmetry we may assume that $v_1$ is adjacent to $z_2$ and $z_3$.
If $\hat R_2$ includes an edge incident with $v_1$ or $v_2$, then
Lemma~\ref{wheelnbhd} implies that $v_0,v_1,v_2$ are all adjacent to
$z_1$ or $z_3$. 
Then we may change our notation so that either $R_2=F_{145}$ or 
$R_2=F_{234}$. In the former case we get a contradiction by the result of the
previous paragraph, and in the latter case we get a contradiction by (2).
Thus  $\hat R_2$ includes no edge incident with $v_1$ or $v_2$.
Hence either $v_0$ is adjacent to $z_1$ and $z_3$, or $v_0$ is adjacent
to an internal vertex $v_3$ of degree five which is adjacent to $z_1$
and $z_3$. In either case there is a vertex of degree five in $G_0$
adjacent to $v_1$, $z_3$, $z_1$, and $v_2$. For this vertex, $z_3,v_2$
is an identifiable pair. Note that $G_{z_3v_2}$ is not $5$-colorable.
We $5$-color the vertices $z_1$, $z_2$, $v_2=z_3$, $z_5$, $z_6$ so that 
each gets a unique color.  
Then this coloring extends to $G_{z_3v_2}$, unless we are in case (ii) of 
Lemma~\ref{7cycle-5coloring} for the following walk on six vertices: 
$z_5, v_2=z_3, z_6, z_2, v_2=z_3, z_6$ in 
$G_{z_3v_2}[\{z_1,z_2,v_2=z_3,z_4,z_5,z_6 \}]$.
This implies that there are two adjacent vertices $w_1$ and $w_2$ such that, 
in $G_0$, $w_1$ is adjacent to $z_2$, $z_6$,
$v_2$, and $z_5$, while $w_2$ is adjacent to $z_6$,
$z_5$, $z_2$, and one of $v_2$, $z_3$. But then the subgraph induced by
the eight vertices: $z_1, z_2, z_3, z_5, z_6, v_2, w_1, w_2$ has all
facial triangles except for perhaps one $5$-cycle. Yet there can be at
most one vertex in the interior of that $5$-cycle. Thus $G_0$ has
at most nine vertices, a contradiction.
This proves (4).

\begin{enumerate}
\item[(5)]   { \it
$z_0\ne z_2,z_3$.
}
\end{enumerate}

\noindent

We may assume to a contradiction that $z_0=z_2$ since the case where $z_0=z_3$ is symmetric.
By (2) $R_2= F_{125}$ or $F_{235}$. Suppose 
first that some edge of $G_0$ is incident with $v_1$ or $v_2$
and lies inside $\hat{R_2}$.
Then $v_0$, $v_1$, and $v_2$ are all adjacent to $z_5$, for otherwise
we may change our notation so that $\hat R_2=F_{126}$, contrary to (2). 
Let $v_4$ and $v_5$ be neighbors of $v_0$ such that the cyclic order around
$v_0$ is $v_1$, $z_5$, $v_2$, $v_5$, $v_4$. 
Now notice that $z_1$ is degree five in
$G_0$ and $N_{G_0}(z_1)$ has a subgraph isomorphic to $K_5-P_3$.  
Since $N_{G_0}(z_0)$ is missing the
edge $v_1v_2$, one of the edges $v_1v_5$, $v_2v_4$ must be present
or $z_1$ would contradict the choice of $v_0$. This implies that $v_1$
and $v_2$ are both adjacent to $v_j$ for some $j\in \{4,5\}$. Thus the
edges $v_1v_j$, $v_2v_j$ must go to a repeated vertex on the boundary
of $R_1$ or $v_0$ would be in a four-cycle in $G_0$, a contradiction.
Thus $v_j=z_6$ and the edge $v_2z_6$ is already present. The edge
$v_1z_6$ then implies that $z_4$ is degree five in $G_0$ and that
$N_{G_0}(z_4)$ has a subgraph isomorphic to $K_5^-$, 
contrary to Lemma~\ref{K5MINUS}.
Thus $\hat R_2$ includes no edge of $G_0$ incident with $v_1$ or $v_2$.

Now suppose that $R_2 = F_{125}$. 
%We may assume that there are no
%chords in $\hat{R_2}$ from $v_1$ or from $v_2$. 
We may assume that $v_1$ is adjacent to $z_3,z_4,z_5$.
Then either the cyclic
order around $v_0$ is $v_1$, $z_5$, $z_1$, $v_2$, and an unspecified
vertex $v_3$, or $v_0$ is adjacent to a vertex $v_3$ of degree five
with cyclic order: $v_1$, $z_5$, $z_1$, $v_2$, $v_0$. In either case,
$z_1v_1$ is an identifiable pair for a vertex of degree five in $G_0$.
Note that $G_{v_1z_1}$ is not $5$-colorable. We $5$-color the vertices
$v_1=z_1$, $z_3$, $z_4$, $z_5$, $z_6$ of  $G_{v_1z_1}$
such that each gets a unique color.
Since this coloring does not extend to $G_{v_1z_1}$ we deduce from
Lemma~\ref{7cycle-5coloring} 
applied to the  walk  $z_6$, $v_1=z_1$,
$z_4$, $z_6$, $z_3$, $z_5$ on six vertices
that case (i) of that lemma holds.
That implies there exists a vertex $w_1$ in $G_0$ that
 is adjacent to $v_1$, $z_4$, $z_6$, $z_3$ and $z_5$. 
Let $H:=G[\{z_1,z_3,z_4,z_5,z_6,v_1,w_1\}]$.
The edge $w_1z_6$ may be embedded in two different ways.
In one way of embedding the edge the graph $H$
%of $G_0$ induced by the vertices  $z_1,z_3,z_4,z_5,z_6,v_1,w_1$ 
has all faces bounded by
triangles, except for one bounded by a $4$-cycle and one bounded by
a $5$-cycle. But then $G_0$ has at most eight vertices by 
Lemma~\ref{7cycle-5coloring}, contrary to Lemma~\ref{THOM-9}.
It follows that the edge  $w_1z_6$ is embedded in such a way that
all faces of $H$ are bounded by triangles, except for one face
bounded by the walk  $z_6z_1z_5v_1w_1z_5$ of length six.
Since $G_0$ has at least ten vertices by Lemma~\ref{THOM-9}, 
we must be in case (iii) of 
Lemma~\ref{7cycle-5coloring} when applied to said walk.
This can happen in two ways. In the first case there are pairwise adjacent
vertices $a,b,c\in V(G_0)$ such that $a$ is adjacent to $z_1,z_5,z_6$,
the vertex $b$ is adjacent to $z_5,v_1,w_1$ and $c$ is adjacent to $w_1,z_5,z_6$.
Now $G_0$ is isomorphic to $L_4$ by an isomorphism that maps
$z_3$ and $z_4$ to the top two vertices in Figure~\ref{fig:allels}(d)
(in left-to-right order), 
$z_6$ and $w_1$ to the vertices in the second row, 
$z_5$ to the unique vertex of degree nine, 
and $z_1,a,c,b,v_1$ to the last row of vertices in that figure.
In the second case there are pairwise adjacent vertices
 $a,b,c\in V(G_0)$ such that $a$ is adjacent to $z_1,z_5,v_1$,
the vertex $b$ is adjacent to $z_5,v_1,w_1$ and $c$ is adjacent to $z_1,z_5,z_6$.
Now $G_0$ is isomorphic to $L_3$ by an isomorphism that maps
the top row of vertices in Figure~\ref{fig:allels}(c) to
$z_6, z_3,z_4,w_1$ (again in left-to-right order),
the middle row to $c,z_1,z_5,v_1,b$ and the bottom vertex to $a$.
Since either case leads to a contradiction, this completes the case 
$R_2 = F_{125}$.

%It now follows that $G_0$ is isomorphic to $L_3$, as follows. 
%There is one vertex of degree nine, $z_5$.
%There are six vertices of degree six, two of which, $z_3, z_4$, are adjacent to all other vertices of degree
%six. The remaining four induce a path: $z_1z_6w_1v_1$. Finally there are three pairwise adjacent vertices of degree five.
%One is adjacent to $z_1$ and $z_6$, the second to $z_6$ and $w_1$, and the third to $w_1$ and $v_1$.
%This completes the case  $R_2 = F_{125}$.

It follows that $R_2 = F_{235}$. 
We may assume that $v_2$ is adjacent to $z_1,z_5,z_6$.
Then either the cyclic order around
$v_0$ is $v_1$, $z_3$, $z_5$, $v_2$, and an unspecified vertex $v_3$,
or $v_0$ is adjacent to a vertex $v_3$ of degree five with cyclic
order: $v_1$, $z_3$, $z_5$, $v_2$, $v_0$. Note that $z_1$ is degree
five in $G_0$ and $N_{G_0}(z_1)$ has a subgraph isomorphic to $K_5-P_3$. 
Thus in either case,
$v_2z_3$ is an identifiable pair for a vertex of degree five in $G_0$,
for otherwise $N_{G_0}(z_1)$ has a subgraph isomorphic to $K_5^-$, a contradiction.
Note that $G_{v_2z_3}$ is not $5$-colorable. We $5$-color the vertices
$z_1$, $v_2=z_3$, $z_4$, $z_5$, $z_6$ of  $G_{v_2z_3}$
such that each gets a unique color.
Since this coloring does not extend to $G_{v_2z_3}$,
we deduce that the $6$-walk $z_6v_2=z_3z_4z_6v_2=z_3z_5$  
satisfies (ii) of Lemma~\ref{7cycle-5coloring}.
% for the following $6$-walk: $z_6$, $v_2=z_3$,$z_4$, $z_6$, $v_2=z_3$, $z_5$.  
%$w_1$ and $w_2$. 
Thus, in $G_0$,
there exists two adjacent vertices $w_1$ and $w_2$ such that $w_1$ is adjacent
 to $z_4$, $z_6$, $z_3$, and $z_5$, while $w_2$ is adjacent to
  $z_4$, $z_5$, $z_6$ and $v_2$. But then $w_1$ is degree five in $G_0$
   and $N_{G_0}(w_1)$ has a subgraph isomorphic to $K_5^-$, a contradiction.
This proves (5).

\begin{enumerate}
\item[(6)]   { \it
$z_0\ne z_4,z_5$.
}
\end{enumerate}

\noindent
To prove (6) we may assume for a contradiction that 
$z_0=z_4$ since the case where $z_0=z_5$ is symmetric.
Thus $R_2 = F_{134}$ or $F_{145}$ by (2).
Assume first that $R_2 = F_{145}$, and that $\hat{R_2}$ includes no
edges incident with $v_1$ or $v_2$.
%Suppose that there are no chords
%in $\hat{R_2}$ from $v_1$ or from $v_2$. Now suppose that $R_2 = F_{145}$.
Then either the cyclic order around $v_0$ is $v_1$, $z_1$, $z_5$,
$v_2$, and an unspecified vertex $v_3$, or $v_0$ is adjacent to a
vertex $v_3$ of degree five with cyclic order: $v_1$, $z_1$, $z_5$,
$v_2$, $v_0$.  If the edge $v_1z_5$ is present, then in the subgraph of $G_0$
induced by $z_1$, $z_2$, $z_3$, $z_5$, $z_6$ and $v_2$, there is only
one face that is not bounded by a triangle or $4$-cycle---the following walk on 
six vertices: $z_5$, $z_3z_6z_5z_1v_2$. 
Thus there are at most nine vertices in $G_0$ by Lemma~\ref{7cycle-5coloring}, 
contrary to Lemma~\ref{THOM-9}. 
Hence, in either case $v_1z_5$ is an identifiable pair for a vertex of degree five
in $G_0$. Note that
$G_{v_1z_5}$ is not $5$-colorable. We $5$-color the vertices 
$z_1$, $z_2$, $z_3$, $v_1=z_5$, $z_6$ of $G_{v_1z_5}$ such that each gets a 
unique color. Since this $5$-coloring does not extend to a $5$-coloring 
of $G_{v_1z_5}$ we deduce that case (ii) of Lemma~\ref{7cycle-5coloring} holds
for the following walk on six vertices: 
$z_6$, $z_2$, $v_1=z_5$, $z_6$, $z_3$, $v_1=z_5$.
Thus, in $G_0$, there are two adjacent vertices $w_1$ and $w_2$ such that
$w_1$ is adjacent to $z_2$, $z_6$, $z_5$, and $z_3$, while $w_2$ is adjacent to
 $z_2$, $z_6$, $z_3$ and $v_1$. But then $w_1$ is
degree five in $G_0$ and $N_{G_0}(w_1)$ has a subgraph isomorphic to $K_5^-$, 
contrary to Lemma~\ref{K5MINUS}.
This completes the case when  $R_2 = F_{145}$ and $\hat{R_2}$ includes no
edges incident with $v_1$ or $v_2$.

For the next case assume that $R_2 = F_{134}$, and again that  $\hat{R_2}$ 
includes no edges incident with $v_1$ or $v_2$.
Then either the cyclic order around $v_0$
is $v_1$, $z_3$, $z_1$, $v_2$, and an unspecified vertex $v_3$, or
$v_0$ is adjacent to a vertex $v_3$ of degree five with cyclic order:
$v_1$, $z_3$, $z_1$, $v_2$, $v_0$.  
Next we dispose of the case that $v_2$ is adjacent to $z_3$.
In that case we consider the subgraph of $G_0$ 
induced by $z_1$, $z_2$, $z_3$, $z_5$, $z_6$
and $v_2$. There is only one face that is not bounded by a triangle or 
$4$-cycle---the following walk on seven vertices: 
$z_5z_3v_2z_1z_3z_2z_6$.
We $5$-color the subgraph as follows: $c(z_i) = i$ for $i=1,2,3,5$, 
 $c(z_6)=4$, and $c(v_2)=2$ and apply Lemma~\ref{7cycle-5coloring}.
By Lemma~\ref{THOM-9} cases (v) or (vi) of Lemma~\ref{7cycle-5coloring} hold.
%However, case (v) does not hold, given the order of colors on the boundary of
%the $7$-walk. Thus case (vi) holds.
Since $z_2$ and $v_2$ have the same color and $z_3$ is a repeated vertex
it follows from Lemma~\ref{7cycle-5coloring} that $G_0$ has four vertices
$a,b,c,d$ such that $d$ is adjacent to $z_2,z_3,z_5,z_6$, the vertices
$a,b,c$ form a triangle and either $a$ is adjacent to $z_1,v_2,z_3$,
the vertex $b$ is adjacent to $z_1,z_2,z_3$, and $c$ is adjacent to
$z_2,z_3,d$ (case (v) of Lemma~\ref{7cycle-5coloring}), or
 $a$ is adjacent to $z_1,v_2,z_3$,
the vertex $b$ is adjacent to $v_2,z_3,d$, and $c$ is adjacent to
$z_2,z_3,d$ (case (vi) of Lemma~\ref{7cycle-5coloring}).
In the former case $d$ is an s-vertex, and yet $v_0=a$, $c$ is not adjacent
to $z_1$ and $b$ is not adjacent to $v_2$, contrary to the optimality
of $(G_0,v_0)$.
In the latter case $G_0$ is isomorphic to $L_3$ by a mapping that sends
the top row of vertices in Figure~\ref{fig:allels}(c) to
$z_1, z_6, z_5, z_2$ (in left-to-right order),
the middle row to $a,v_2,z_3,d,c$ and the bottom vertex to $b$,
a contradiction.
Thus  $v_2$ is not adjacent to $z_3$, and hence $v_2z_3$ is an identifiable
pair for a vertex of degree five in $G_0$. Note that $G_{v_2z_3}$ is not
$5$-colorable. We $5$-color the vertices $z_1$, $z_2$, $v_2=z_3$, $z_5$, $z_6$
of $G_{v_2z_3}$
such that each gets a unique color. 
Since this coloring not extend to $G_{v_2z_3}$ we deduce that
case (ii) of Lemma~\ref{7cycle-5coloring} holds for the following $6$-walk:
$z_6$, $z_2$, $z_3=v_2$, $z_6$, $z_3=v_2$, $z_5$.  However this would imply that
there are two internal vertices $w_1$ and $w_2$, both adjacent to $z_2$
 and both adjacent to $z_5$. 
But then one of them is not adjacent to $z_3=v_2$, a contradiction.
This completes both cases when  $\hat{R_2}$ 
includes no edges incident with $v_1$ or $v_2$.

We continue  the proof of (6). We have just shown that  $\hat{R_2}$ includes
an edge incident with $v_1$ or $v_2$.
Then $v_0,v_1,v_2$ are all adjacent to $z_1,z_3$ or $z_5$.
However, if they are all adjacent to $z_3$, then we can change notation so
that $R_2=F_{234}$, contrary to (2), and if they are all adjacent to $z_5$,
then we can change notation so that $R_2=F_{456}$, again contrary to (2).
Thus $v_0,v_1,v_2$ are all adjacent to $z_1$.
We may assume that the notation is chosen so that  $v_1$ is adjacent to $z_2$ and
$z_3$ while $v_2$ is adjacent to $z_5$ and $z_6$.
Let $v_4$ and $v_5$ be neighbors of $v_0$ numbered so that the 
cyclic order around $v_0$ is $v_2,z_1,v_1,v_4,v_5$.

Next we claim that  $v_1$ is not adjacent to $z_6$. Suppose it were. 
The triangle $z_2v_1z_6$ is null-homotopic in $G_0$ by Lemma~\ref{7cycle-5coloring}
applied to the $4$-cycle $z_1z_5z_6v_1$.
Now consider the subgraph induced by the vertices $z_1$, $z_2$, $z_3$,
$z_5$, $z_6$, and $v_1$. All of its faces are triangles but for the
$7$-walk $z_1z_5z_6z_3z_5z_6v_1$. We $5$-color
these vertices as follows: $c(z_i) = i$ for $i=1,3,5$, $c(z_6)=4$, and $c(v_1)=5$.
 Now we must be in case (v) or (vi) of Lemma~\ref{7cycle-5coloring}, for otherwise
$|V(G_0)| \le 9$, contrary to Lemma~\ref{THOM-9}. Yet, since the fifth color would
appear three times on the boundary, we can extend this coloring to
all of $G_0$, a contradiction. Thus  $v_1$ is not adjacent to $z_6$.

Now we claim that $v_4,v_5\not\in\{z_1,z_2,\ldots,z_6\}$.
To prove this claim we suppose the contrary. Then $v_0$ is adjacent to
$z_2$, $z_3$, $z_5$ or $z_6$.
If $v_0$ is adjacent to $z_2$, then  $v_1$ has degree at most four in $G_0$. 
If $v_0$ is adjacent to $z_6$, then either $v_2$  is degree
four in $G_0$, a contradiction, or $v_1$ is adjacent to $z_6$, a
contrary to the previous paragraph.
If $v_4=z_3$, then  the $5$-cycle $v_1z_3z_6z_5z_1$ has the vertices
$v_0$ and $v_2$ in its interior, contrary to Lemma~\ref{7cycle-5coloring}.
Let us assume that $v_5=z_3$.
%, then the $5$-cycle $v_0z_3z_6z_5z_1$ has the vertices
%$v_0$ and $v_2$ in its interior, again a contradiction.
%Suppose $v_0$ is adjacent to $z_3$. 
Then $v_2$ is degree
five and $N(v_2)$ is missing at most the edges $v_0z_5$ and $v_0z_6$.
Yet these edges must not be present, for otherwise 
$N(v_2)$ has a subgraph isomorphic to $K_5^-$, contrary to Lemma~\ref{K5MINUS}.
Hence $v_4\not\in\{z_1,z_2,\ldots,z_6\}$, but then it is not
adjacent to $z_1$. Thus $N_{G_0}(v_0)$ includes two disjoint edges. However,
$N_{G_0}(v_2)$ has a subgraph isomorphic to $K_5-P_3$, 
contradicting the optimality of $(G,v_0)$. Thus we may assume that
$v_0$ is adjacent to $z_5$. This implies, by Lemma~\ref{wheelnbhd}, 
that $v_4=z_5$, because 
$v_2$ is already adjacent to $z_5$ and $v_5=z_5$ would imply the
existence of another edge from $v_2$ to $z_5$, not homotopic to the
existing one.
Then the subgraph of $G_0$ induced by $z_1$,
$z_2$, $z_3$, $z_5$, $z_6$, and $v_1$ has only one 
face---a six-walk---that can have vertices in its interior. 
But then there are at most
nine vertices in $G_0$ by Lemma~\ref{7cycle-5coloring}, 
contrary to Lemma~\ref{THOM-9}.
This proves our claim  that $v_4,v_5\not\in\{z_1,z_2,\ldots,z_6\}$.

Continuing with the proof of (6), we note that $v_2$ is
not adjacent to $v_4$, for otherwise $v_5$ is of degree four in $G_0$, a
contradiction. Similarly $v_1$ is not adjacent to $v_5$. Since $z_1$
is not adjacent to $v_4$ or $v_5$, the neighborhood of $v_0$ in $G_0$ is
a cycle of length five.
The vertex $v_2$ is not adjacent to $z_2$, for otherwise the $4$-cycle
$z_2v_2v_0v_1$ includes the vertices $v_4$ and $v_5$ in its interior, 
contrary to Lemma~\ref{7cycle-5coloring}.
Furthermore, the vertex $v_4$ is not adjacent to $z_2$, for otherwise
the neighborhood of $v_1$ in $G_0$ has a subgraph isomorphic to 
a $5$-cycle plus one edge, contrary to the optimality of $(G_0,v_0)$.
We now consider the graph $G_{v_2v_4}$.
It has a subgraph $H$ isomorphic to $K_6$, and the new vertex $w$ of $H$
obtained by identifying $v_2$ and $v_4$ belongs to $H$.
Let $\Delta$ denote the open disk bounded by the walk
$z_1z_5z_6z_3z_5z_6z_2z_3$ of  $G_{v_2v_4}$.
Since $w$ belongs to $\Delta$, all vertices of $H$ belong to the closure
of $\Delta$.
However, $z_2\not\in V(H)$, because $z_2$ is not adjacent to $v_2$ or $v_4$
in $G_0$.
Since $v_1$ is not adjacent to $z_6$ as shown two paragraphs ago,
we deduce that not both $z_6$ and $v_1$ belong to $H$.
That implies that $z_1\not\in V(H)$, because at most six neighbors of $z_1$ 
in $G_{v_2v_4}$ (including $z_2\not\in V(H)$) belong to the closure of $\Delta$.
%has at most six
%neighbors in $G_{v_2v_4}$ in the closure of $\Delta$, including $z_2\not\in V(H)$.
If $v_1\not\in V(H)$, then no edge incident with one of the two occurrences
of $z_3$ on the boundary of $\Delta$ belongs to $H$.
Thus regardless of which of $v_1$, $z_6$ does not belong to $H$,
there is a planar graph $H'$ obtained from $H$ by splitting at 
most two vertices, and a drawing of $H'$ in the unit disk with vertices
$p,q,r,s$ drawn on the boundary in order such that $H$ is obtained
from $H'$ by identifying $p$ with $r$, and $q$ with $s$.
It follows that $H$ can be made planar by deleting one vertex, contrary to
the fact that it is isomorphic to $K_6$.
This proves (6).
\medskip

Since $R_1=F_0$ it follows that $z_0\ne z_1$. Thus $z_0=z_6$ by (5) and (6).

\begin{enumerate}
\item[(7)]   { \it
We may assume that $R_2\ne F_{136}$ and  $R_2\ne F_{126}$.
}
\end{enumerate}

\noindent
To prove (7) we may assume for a contradiction by symmetry that  $R_2=F_{136}$.
Then by (2) we have $R_1=F_{264}$.
We may assume that $v_1$ and $v_2$ are numbered so that $v_1$ is adjacent to
$z_1$ and $z_2$.
We may assume that  $\hat{R_2}$ includes no edge incident with $v_1$ or $v_2$;
for if it includes the edge $v_2z_1$, then we can change notation so
that $R_2=F_{126}$, contrary to (2), and if it includes the edge $v_1z_3$, then
we can change notation and reduce to the case when $R_2=F_0$, 
which is handled below.
Then either the cyclic order
around $v_0$ is $v_1$, $z_1$, $z_3$, $v_2$, and an unspecified vertex
$v_3$, or $v_0$ is adjacent to a vertex $v_3$ of degree five with
cyclic order: $v_1$, $z_1$, $z_3$, $v_2$, $v_0$. In either case, $z_1,v_2$
is an identifiable pair for a vertex of degree five in $G_0$. Note that
$G_{v_2z_1}$ is not $5$-colorable. We $5$-color the vertices 
$z_1=v_2, z_2, z_3, z_4,z_5$ of $G_{v_2z_1}$ such that each gets a unique color. 
Since this coloring does not extend to the rest of $G_{v_2z_1}$
we deduce that case (i) of Lemma~\ref{7cycle-5coloring} holds 
for the following $6$-walk on five vertices:
 $z_1v_2$, $z_2$, $z_4$, $z_1v_2$, $z_3$, $z_5$. 
This implies that there exists a vertex $w_1$
in $G_0$ such that $w_1$ is adjacent to 
$z_2$, $z_4$, $v_2$, $z_3$ and $z_5$ in $G_0$.
In the subgraph of $G_0$ induced by those six vertices and $z_1$,
all the faces are triangles but for the face bounded by the
%following walk on six vertices:
cycle $z_1z_3v_2z_5w_1z_2$.  Since $G_0$ must have at
least ten vertices, we must be in case (iii) of Lemma~\ref{7cycle-5coloring}.
Now $5$-color the subgraph induced by those six vertices and $z_4$
such that $c(z_i)=i$ for $i=1,2,3,5$, $c(w_1)=1$, and $c(v_2)=2$. 
The above-mentioned cycle is colored using four colors, and hence
the $5$-coloring may be extended to $G_0$, a contradiction.
This proves (7).
\medskip

In light of (7) we may assume that both $R_1$ and $R_2$ are equal to $F_0$.
Thus we may assume that $R_1=F_{264}$ and  $R_2=F_{365}$.
We may assume that $v_1$ and $v_2$ are numbered so that $v_1$ is adjacent to
$z_1,z_2$ and $z_3$.
Let the remaining neighbors of $v_0$ be $v_3,v_4,v_5$ numbered so that
the cyclic order around $v_0$ is $v_1,v_3,v_2,v_5,v_4$.
This specifies the cyclic order uniquely up to reversal, and so 
we may assume by symmetry that the cyclic order around $v_1$ (of 
a subset of the neighbors of $v_1$) is
$z_1,z_3,v_3,v_0,v_4,z_2$, where possibly $v_3=z_3$ and $z_2=v_4$.

\begin{enumerate}
\item[(8)]   { \it
The vertex  $v_1$ is not adjacent to $z_4$ or $z_5$.
}
\end{enumerate}

\noindent To prove (8) we note that $z_1$ has degree five in $G_0$
and that its neighborhood has a subgraph isomorphic to $K_5-P_3$.
If $v_1$ was adjacent to  $z_4$ or $z_5$, then the neighborhood of
$z_1$ would have a subgraph isomorphic to $K_5^-$, contrary to
Lemma~\ref{K5MINUS} and the optimality of $(G_0,v_0)$.
This proves (8).

Since $z_1$ has degree five in $G_0$
and its neighborhood has a subgraph isomorphic to $K_5-P_3$,
we deduce from the optimality of $(G_0,v_0)$ and Lemma~\ref{K5MINUS} that
the neighborhood of $v_0$ is isomorphic to  $K_5-P_3$. It follows that

\begin{enumerate}
\item[(9)]   { \it
the vertex  $v_3$ is adjacent to $v_4$ or $v_5$
}
\end{enumerate}

\noindent
and

\begin{enumerate}
\item[(10)]   { \it
either $v_1$ is adjacent to $v_5$, or $v_2$ is adjacent to $v_4$, and not both.
}
\end{enumerate}

\begin{enumerate}
\item[(11)]   { \it
The vertex $v_2$ is  adjacent to  $v_4$.
}
\end{enumerate}

\noindent
To prove (11) suppose for a contradiction that $v_2$ and $v_4$ are not
adjacent. 
We will consider $G_{v_2v_4}$ and its new vertex $w$ formed by
identifying $v_2$ and $v_4$. Let us note that all faces of the subgraph 
of $G_{v_2v_4}$ induced
by $z_1$, $z_2$, $z_3$, $z_4$, $z_5$, $v_1$, $w$
are bounded by triangles except for a face bounded by the $8$-walk
$W_1=v_1wz_5z_3v_1wz_4z_2$. 
Let $D_1$ be the open disk bounded by $W_1$, let
$W_0=v_1v_4v_5v_2z_5z_3v_1v_3v_2z_4z_2$ be a corresponding walk in $G_0$,
and let $D_0$ be the open disk bounded by $W_0$.
By Lemma~\ref{J=K6} the graph $G_{v_2v_4}$ 
has a subgraph $H$ isomorphic to $K_6$. 
Since $G$ has no $K_6$ subgraph it follows that $w\in V(H)$.
If $z_1\in V(H)$, then, since $z_1$ has degree five in $G_0$, all
neighbors of $z_1$ belong to $V(H)$, contrary to (8).
Thus all vertices of $H$ belong to $W_1$ or $D_1$, and by Lemma~\ref{wheelnbhd}
each vertex of $H\backslash w$ (when regarded as a vertex of $G_0$)
belongs to $W_0$ or $D_0$.
Assume for a moment that all but possibly one vertex of $H$ belong to $W_1$.
Then $z_4$ or $z_5$ belongs to $V(H)$, and so $v_1\not\in V(H)$ by (8).
Thus exactly one vertex of $H$, say $w_1$,  belongs to $D_1$ and
$V(H)=\{w,w_1,z_2,z_3,z_4,z_5\}$.
It follows that $v_4\not\in\{w_1,z_2,z_3,z_4,z_5\}$.
%, and $v_4\not\in\{z_4,z_5\}$, because $v_2$ is not adjacent to $v_4$.
Thus $v_4$ is
not adjacent to $z_3$ in $G_0$, because the edge $z_3v_4$ would have to
lie in $D_0$, where it would have to
%the edge $z_3v_4$ lies in $D_0$, and hence $v_4$ is
%not adjacent to $z_3$ in $G_0$, because the edge $z_3v_4$ would have to
cross the path $z_4w_1z_5$.
But $w$ is adjacent to $z_3$ in $H$, and so $v_2$ is adjacent to $z_3$
in $G_0$. 
It follows that the the $4$-cycle $v_1v_0v_2z_3$ is null-homotopic,
for otherwise the edge $v_2z_3$ and path $z_2w_1z_5$ would cross in $D_0$.
%The presense of the vertex $w_1$ shows that 
%the $4$-cycle
%v_1v_0v_2z_3$ is null-homotopic;
We deduce from Lemma~\ref{7cycle-5coloring} applied to the $4$-cycle
$v_1v_0v_2z_3$ that $v_3=z_3$. But $v_3$ is adjacent to $v_4$
by (9), and yet $z_3$ is not adjacent to $v_4$, a contradiction.
This completes the case when at most one vertex of $H$ belongs to~$D$.

Thus at least two vertices of $H$, say $w_1$ and $w_2$ belong to
of $D_1$. Since $W_1$ has exactly two repeated vertices,
the argument used at the end of the proof of (6) shows that $w_1$ and
$w_2$ are the only two vertices of $H$ in $D_1$.
Also, it follows that $w,v_1$, the two repeated vertices of $W_0$,
belong to $H$. Since $v_1$ is in $H$, (8) implies that $z_4,z_5\not\in V(H)$.
It follows that $z_2,z_3\in V(H)$, and consequently $v_4\not\in\{z_2,z_3\}$.
Thus each of $w_1,w_2$ is adjacent in $G_0$ to
$v_1,z_2,z_3$ and to $v_2$ or $v_4$.
It follows from considering the drawing of $G_0$ inside $D_0$ that
one of $w_1,w_2$, say $w_1$, is adjacent to $v_2$ and the
$4$-cycle $v_1v_0v_2w_1$ is null-homotopic.
By Lemma~\ref{7cycle-5coloring} applied to this $4$-cycle 
we deduce that $w_1=v_3$. Thus the edge $v_3v_4$ belongs to $D_0$.
But $w_2\ne v_4$, because $v_4$ is not a vertex of $H$, 
and yet the edge $v_3v_4$ intersects the path $z_3w_2z_2$
inside $D_0$, a contradiction.
This proves (11).
\medskip

%Since $v_4$ is adjacent to $v_3$ by (9) and to
%$v_1$, we deduce that $w_2=v_4$, and hence $w_1=v_3$ and $w_2=v_4$
%are both  adjacent to $v_5$, the former by (9).
%Thus $D\cup\{z_2,z_3\}$ includes a closed disk $D'$ that has $z_2,z_3,v_5$
%on its boundary, and $w_1,w_2$ in its interior, and such that $D'$ includes
%all edges from $w_1$ and $w_2$ to $z_2,z_3,v_5$.
%That is impossible, because upon adding a new vertex in the complement
%of $D'$ and joining it to $z_2,z_3,v_5$ we would obtain a planar drawing
%of $K_{3,3}$.

\begin{enumerate}
\item[(12)]   { \it
The vertex $v_5$ is  adjacent to  $v_1$.
}
\end{enumerate}

\noindent
We prove (12) similarly as the previous claim.
Suppose for  a contradiction that $v_1$ and $v_5$ are not adjacent,
and consider $G_{v_1v_5}$ and its new vertex $w$. The subgraph of $G_{v_1v_5}$
induced by $z_1$, $z_2$, $z_3$, $z_4$, $z_5$, $w$, $v_2$ has all
faces bounded by triangles except for one bounded by the 
$8$-walk $W_1=wv_2z_5z_3wv_2z_4z_2$. 
Let $D_1$ be the open disk bounded by $W_1$, and let $W_0,D_0$ be as in (11).
Similarly as in the proof of (11) the graph $G_{v_1v_5}$ has a subgraph
$H$ isomorphic to $K_6$ with $w\in V(H)$. %and $z_1\not\in V(H)$.
We claim that $z_4\not\in V(H)$. Indeed, if $z_4$ is in $H$, then it is
adjacent to $w$ in $H$; but $z_4$ is not adjacent in $G_0$ to $v_1$ by (8), 
and hence $z_4$ is adjacent to $v_5$ in $G_0$. Yet $v_2$ is adjacent to
$v_4$ by (10). Since $v_4\not\in\{z_4,z_5\}$ by (8), the edges 
$v_2v_4$ and $z_4v_5$ must cross inside $D_0$, a contradiction.
This proves our claim that $z_4\not\in V(H)$.
It follows that $z_1\not\in V(H)$, because $z_1$ has degree five in $G_{v_1v_5}$,
and $z_4$ is one of its neighbors.

If $D_1$ includes at most one vertex of $H$, then $w,v_2,z_2,z_3,z_5\in V(H)$,
and exactly one vertex of $H$, say $w_1$, belongs to $D_1$.
%Thus $v_2$ is adjacent to $z_2$ and $z_3$ in $G_0$. By (9)
Thus $w_1$ is adjacent to $z_2$ and $z_5$ in $G_0$, and that implies that
the edges $v_3v_4$ and $v_3v_5$ do not lie in $D_1$.
Therefore $v_3,v_4,v_5\in\{z_2,z_3,z_4,z_5\}$, but that is impossible,
given the existence of $w_1$. 
This completes the case that $D_1$ includes at most one vertex of $H$.
Thus, similarly as in (11), it follows that $D_1$ includes exactly two
vertices of $H$, say $w_1$ and $w_2$.
Now $V(H)$ includes $w,v_2$ and exactly two of $\{z_2,z_3,z_5\}$.
But it cannot include $z_5$ and $z_3$, because otherwise for some
$j\in\{1,2\}$ the paths $z_5w_jv_2$ and $z_3w_{3-j}v_2$ cross inside $D_0$.
Thus $V(H)$ includes $z_2$ and $z_i$ for some $i\in\{3,5\}$.
Choose $j\in\{1,2\}$ such that $w_j\ne v_3$. 
Then the path $z_2w_jz_i$ is not disjoint from  the edges $v_3v_4$,
$v_3v_5$ (because they  cross inside $D_0$), and so it follows that
$i=3$ and $v_3=z_3$.
Since there is no crossing in $D_0$ and $w_1$ and $w_2$ are adjacent to
$z_2$ and $z_3$, they are not both adjacent to $v_5$. Thus we may
assume that $w_1$ is adjacent to $v_1$. This argument shows, in fact, 
that the cycle $v_1v_0v_2w_1$ is null-homotopic, and so it follows from 
Lemma~\ref{7cycle-5coloring} that $v_3=w_1$, a contradiction, because
$w_1$ lies in $D_1$ and $v_3=z_3$ does not.
This proves (12).
\medskip

Now claims (10), (11), and (12) are contradictory. 
This completes the proof of Lemma~\ref{6walk}.~\qed

%\hrule
%\medskip
%TO BE EDITED
%\input insert.tex

%\medskip
\noindent
{\bf Proof of Theorem~\ref{main}.}
It follows by direct inspection that none of the graphs listed in
Theorem~\ref{main} is $5$-colorable.
Conversely, let $G_0$ be a graph drawn in the Klein bottle that is
not $5$-colorable. We may assume, by taking a subgraph of $G_0$,
 that $G_0$ is $6$-critical.
Then $G_0$ has minimum degree at least five.
By Lemma~\ref{6regular-5colorable} the graph $G_0$ has a vertex of degree
exactly five, and so we may select a vertex $v_0$ of $G_0$ such that
$(G_0,v_0)$ is an optimal pair.
If there is no identifiable pair, then $G_0$ has a $K_6$ subgraph,
as desired. Thus we may select an identifiable pair $v_1,v_2$.
Let $G':=G_{v_1v_2}$. 
By Lemma~\ref{J=K6} the graph $G'$ has a subgraph $H$ isomorphic to $K_6$.
By Lemma~\ref{crosscapinface} the drawing of $H$ is $2$-cell, and
by Lemma~\ref{6faceexists} some face of $H$ has length six, contrary to 
Lemma~\ref{6walk}.~\qed

%By the optimality of $(G_0,v_0)$ the graph
%$G'$ has a subgraph isomorphic to one of the graphs listed in
%Theorem~\ref{main}.
%But it has no subgraph isomorphic to $C_3+C_5$ or $K_2+H_7$ by
%Lemma~\ref{G'notc3c5k2h7}, it has no subgraph isomorphic to 
%$L_1,L_2,L_5$ or $L_6$ by Lemma~\ref{L1L2L5L6}, and it has no
%subgraph isomorphic to $L_3$ or $L_4$ by Lemma~\ref{L3L4}.

\baselineskip 11pt
\vfill
\noindent
This material is based upon work supported by the National Science Foundation
under Grant No.~DMS-0701077.
Any opinions, findings, and conclusions or recommendations expressed in
this material are those of the authors and do not necessarily reflect
the views of the National Science Foundation.


\begin{thebibliography}{99}

\def\JCTB{{\it J.~Combin.\ Theory Ser.\ B}}
\def\CMUC{{\it Comment. Math. Univ. Carol.}}
\def\TAMS{{\it Trans.\ Amer.\ Math.\ Soc.}}
\def\JAMS{{\it J.~Amer.\ Math.\ Soc.}}
\def\PAMS{{\it Proc. Amer. Math. Soc.}}
\def\DM{{\it Discrete Math.}}
\def\CM{{\it Contemporary Math.}}
\def\GC{{\it Graphs and Combin.}}
\def\COM{{\it Combinatorica}}
\def\JGT{{\it J.~Graph Theory}}
\def\JAlgorithms{{\it J.~Algorithms}}
\def\SIAMDM{{\it SIAM J.~Disc.\ Math.}}
\def\CPC{{\it Combinatorics, Probability and Computing}}
\def\EJC{Electron.\ J.~Combin.}

%\input refs.tex
%extract_bib.pl 5colkb.tex >! refs.tex


\bibitem{AlbHut} M.~Albertson and J.~Hutchinson,
The three excluded cases of Dirac's map-color theorem,
 Second International Conference on Combinatorial Mathematics 
(New York, 1978),  pp.~7--17,
{\it Ann.\ New York Acad.\ Sci.} {\bf 319}, New York Acad.\ Sci., 
New York, 1979.

\bibitem{AppHak1}K. Appel and W. Haken, 
Every planar map is four colorable, Part I: discharging,
{\it Illinois J. of Math.} {\bf 21} (1977), 429--490.

\bibitem{AppHakKoc}K. Appel, W. Haken and J. Koch, 
Every planar map is four colorable, Part II: reducibility,
{\it Illinois J. of Math.} {\bf 21} (1977), 491--567.

\bibitem{AppHak89}K. Appel and W. Haken, 
Every planar map is four colorable, 
{\it Contemp. Math.} {\bf 98} (1989).

\bibitem{Dirmap} G. A. Dirac, Map color theorems,
{\it Canad.\ J.~Math.} {\bf 4} (1952), 480--490.

\bibitem{Dircritical} G. A. Dirac, The coloring of maps,
{\it J.~London Math.\ Soc.} {\bf 28} (1953), 476--480.

\bibitem{Eppsubiso} D.~Eppstein, 
Subgraph isomorphism in planar graphs and related problems,
{\it J.~Algor.\ and Appl.} {\bf3} (1999), 1--27.

\bibitem{Fis} S.~Fisk, The nonexistence of colorings,
\JCTB\ {\bf 24} (1978), 247--248.

\bibitem{Galcritical} T.~Gallai, Kritische Graphen I, II,
{\it Publ.\ Math.\ Inst.\ Hungar.\ Acad.\ Sci.} {\bf 8} (1963), 165--192
and 373--395.

\bibitem{GarJoh} M.~R.~Garey and D.~S.~Johnson, 
Computers and intractability. A guide 
to the theory of NP-completeness, W. H. Freeman, San Francisco, 1979.

\bibitem{Heawood} P.~J.~Heawood, Map-color theorem,
{\it Quart.\ J.~Pure Appl.\ Math.} {\bf 24} (1890), 332--338.

\bibitem{KawKraKynLid} K.~Kawarabayashi, D.~Kral, J.~Kyn\v{c}l, and B.~Lidick\'y,
6-critical graphs on the Klein bottle,
submitted.

\bibitem{KocKre}W.~Kocay and Donald L. Kreher, 
Graph Algorithms and Optimization, CRC Press, New York, 2004.

\bibitem{KraMohNakPanSuz} 
D.~Kr\'al', B.~Mohar, A.~Nakamoto, O.~Pangr\'ac and Y.~Suzuki,
Coloring Eulerian triangulations on the Klein bottle,
submitted.

\bibitem{MadDirac} W.~Mader,
$3n-5$ edges do force a subdivision of $K_5$,
manuscript.

\bibitem{Massey} W.~S.~Massey, A basic course in algebraic topology,
Springer-Verlag, New York, 1991.

\bibitem{Mohlinear} B.~Mohar,
A linear time algorithm for embedding graphs in an arbitrary surface,
{\it\SIAMDM} {\bf 12} (1999), 6--26.

\bibitem{MohHajos} B.~Mohar,
Triangulations and the Haj\'os conjecture,
{\it Electron.\ J.~Combin.} {\bf 12} N15 (2005).

\bibitem{MohTho} B.~Mohar and C.~Thomassen, Graphs on surfaces, 
Johns Hopkins University Press, Baltimore, MD, 2001.

\bibitem{Ringel} G.~Ringel, Map Color Theorem,
Springer-Verlag, Berlin, 1974.

\bibitem{RobSanSeyTho4CT}
N.~Robertson, D.~P.~Sanders, P.~D.~Seymour and R.~Thomas,
The four-colour theorem, {\it\JCTB} {\bf 70} (1997), 2-44.

\bibitem{RodZich} V.~R\"odl and J.~Zich,
Triangulations and the Haj\'os conjecture,
{\it\JGT} {\bf 59} (2008), 293--325.

\bibitem{Sas6reg} N.~Sasanuma,
Chromatic numbers of $6$-regular graphs on the Klein bottle,
{\it Austral.\ J.~Combin.} {\bf45} (2009), 73--85.

\bibitem{Tho5torus} C. Thomassen, Five-coloring graphs on the torus,
\JCTB\ {\bf62} (1994), 11--33.

\bibitem{Tho5choose} C.~Thomassen, Every planar graph is 5-choosable,
\JCTB\ {\bf62} (1994), 180--181.

\bibitem{ThoCritical} C. Thomassen, Color-critical graphs on a fixed surface,
\JCTB\ {\bf70} (1997), 67--100.

\bibitem{ThoHajos} C. Thomassen,
Some remarks on Haj\'os' conjecture,
{\it\JCTB} {\bf93} (2005), 95--105.

\bibitem{Toft} B.~Toft,
On critical subgraphs of colour-critical graphs,
{\it\DM} {\bf7} (1974), 377--392.

\bibitem{YerPhD} C.~Yerger, Color-critical graphs on surfaces,
Ph.D.\ dissertation, Georgia Institute of Technology, 2010.

\end{thebibliography}
\end{document}